\documentclass[fontsize=12pt,a4paper,headings=normal,
twoside=false,leqno,parskip=half-,abstract=true]{scrartcl}
\usepackage[utf8]{inputenc}
\usepackage{CJKutf8} 
\usepackage[greek,english]{babel}
\usepackage{csquotes}
\usepackage[official]{eurosym}
\usepackage[hyphens]{url}
\usepackage{hyperref}

\hypersetup{
 pdftitle={complex rational ODEs},
 pdfauthor={Bernold Fiedler},
 colorlinks=true,
 linkcolor=blue,
 citecolor=blue,
 filecolor=blue,
 urlcolor=blue} 

 \usepackage{afterpage}

\usepackage{graphicx}
\usepackage{wrapfig} 
\usepackage[format=plain,labelfont=bf,font=small]{caption}
\usepackage{xcolor}
\usepackage[arrow, matrix, curve]{xy}
\usepackage{float}

\usepackage{caption}
\usepackage{tabulary}
\usepackage{array}
\newcolumntype{N}[1]{>{\centering\arraybackslash}m{#1}}

\usepackage{amsmath,amsthm}
\swapnumbers 
\usepackage{amssymb, eurosym} 

\makeatletter
\newcommand{\tpitchfork}{%
  \vbox{
    \baselineskip\z@skip
    \lineskip-.52ex
    \lineskiplimit\maxdimen
    \m@th
    \ialign{##\crcr\hidewidth\smash{$-$}\hidewidth\crcr$\pitchfork$\crcr}
  }%
}
\makeatother
\usepackage{latexsym}
\usepackage{enumerate}

\usepackage[notref,notcite,color,final 
]{showkeys}
\usepackage[textwidth=2cm,textsize=tiny,backgroundcolor=none,textcolor=red]{todonotes}

\definecolor{refkey}{rgb}{1,0,0}
\definecolor{labelkey}{rgb}{1,0,0}

\usepackage{cancel}

\usepackage{tikz}

  \mathchardef\ordinarycolon\mathcode`\:
  \mathcode`\:=\string"8000
  \begingroup \catcode`\:=\active
    \gdef:{\mathrel{\mathop\ordinarycolon}}
  \endgroup

\theoremstyle{plain}
\newtheorem{thm}{Theorem}[section]
\newtheorem{lem}[thm]{Lemma}
\newtheorem{prop}[thm]{Proposition}
\newtheorem{cor}[thm]{Corollary}

\newtheorem{defi}[thm]{Definition}

\newcommand\eps{\varepsilon}
\newcommand\mi{\mathrm{i}}
\renewcommand\theta{\vartheta}
\renewcommand\rho{\varrho}
\renewcommand\phi{\varphi}
\renewcommand\Re{\mathrm{Re}\,}
\renewcommand\Im{\mathrm{Im}\,}

\begin{document}

\title{\LARGE{Real-time blow-up and connection graphs\\ 
\medskip
of rational vector fields on the Riemann sphere }}
\vspace{1cm}
{\subtitle{}
	\vspace{1ex}
	{}}\vspace{1ex}

\author{
 \\
\emph{-- In memoriam Professor Masaya Yamaguti --}\\
{~}\\
Bernold Fiedler*
\\
\vspace{2cm}}

\date{\small{version of \today}}
\maketitle
\thispagestyle{empty}

\vfill

*\\
Institut für Mathematik\\
Freie Universität Berlin\\
Arnimallee 3\\ 
14195 Berlin, Germany


\newpage
\pagestyle{plain}
\pagenumbering{roman}
\setcounter{page}{1}

\begin{abstract}
\noindent
\noindent
Inspired by pioneering work of Kyûya Masuda in the 1980s, only much more recent PDE studies address global boundedness versus finite-time blow-up in complex time.
The two phenomena are related by passage from real to purely imaginary time.

\smallskip\noindent
As a most simplistic ODE example, we study scalar rational vector fields
\begin{equation*}
\label{*}
\dot{w}=P(w)/Q(w)\,,  \tag{*}
\end{equation*}
for complex polynomials $P,Q$.
We impose mild generic nondegeneracy conditions, including simplicity of poles and hyperbolicity of zeros.
Generically, the real-time dynamics then become gradient-like Morse.
Poles play the role of hyperbolic saddle points.
At poles, however, solutions may blow up in finite time.

\smallskip\noindent
On the Riemann sphere $w\in\widehat{\mathbb{C}}$, we classify the resulting global dynamics up to $C^0$ orbit equivalence, in real time.
This relies on a global description of the connection graph of blow-up orbits, from sources towards saddles/poles, in forward time.
Time reversal identifies the dual graph of blow-down orbits.
We show that the blow-up and blow-down graphs of (*) on $\widehat{\mathbb{C}}$ \emph{realize all finite multi-graphs on $\mathbb{S}^2$}, equivalently.

\smallskip\noindent
The purely polynomial case $Q=1$ realizes all planar trees, alias diagrams of \emph{non-intersecting circle chords}.
The anti-holomorphic cousin $P=1$ realizes all \emph{noncrossing trees} with vertices restricted to circles.
This classification provides combinatorial counts for the number of global phase portraits, which only depend on the degrees of $P$ and $Q$, respectively.

\end{abstract}

\vspace{2cm}
\tableofcontents
\listoffigures


\newpage
\pagenumbering{arabic}
\setcounter{page}{1}

\section{Introduction} \label{Int}

\numberwithin{equation}{section}
\numberwithin{figure}{section}
\numberwithin{table}{section}

\subsection{PDE motivation}\label{PDE}
In the present ODE paper, we study scalar rational vector fields
\begin{equation}\label{ODE}
\dot{w}=f(w)=P(w)/Q(w)\,,
\end{equation}
for coprime complex polynomials $P,Q$.
Our motivation originates from striking attempts by Masuda \cite{Masuda1,Masuda2} in the 1980s to circumnavigate PDE blow-up, in finite real time $t=r\nearrow r^*<\infty$, by a detour via the complex time domain $t=r+\mi s$.
Reaction-diffusion PDEs, such as the nonlinear heat equation
\begin{equation}
\label{PDEw}
w_r=w_{xx}+f(w)\,,
\end{equation}
distinguish the real time direction $t=r$ as physical time, by traditional derivation, and solutions $w=w(t,x)$ are thought of as real-valued.
Indices indicate partial derivatives.

Evaluation of complex analytic power series solutions $t\mapsto w(t,x)$ of \eqref{PDEw} in the \emph{imaginary time direction}, in contrast, leads to solutions of nonconservative variants
\begin{equation}
\label{PDEpsi}
\mi\psi_s=\psi_{xx}+f(\psi)
\end{equation}
of the Schrödinger equation.
An example is the 1-parameter family 
\begin{equation}\label{psiw}
\psi(s,x):=w(r_0-\mi s,x)\,.
\end{equation}
The parameter $r=r_0$ is any fixed real part of complex time $t=r_0-\mi s$, and the Schrödinger solution proceeds along the vertical imaginary ``time'' direction $\mi s$.
In this way, complex time extension provides families of Schrödinger solutions $\psi$, from a single heat solution $w$, and vice versa. 
Moreover, the families of one type are related by a (semi)flow of the other type.
In fact, the two semiflows commute, locally and on analytic solutions.
The simplest case of homogeneous solutions $w=w(t)$ leads back to the ODE case \eqref{ODE}.

Here and below, the terms \emph{``analytic''} and \emph{``analyticity''} refer to local expansions by convergent power series.
\emph{``Holomorphy''} refers to complex differentiability and, therefore, Cauchy-Riemann equations.
For real differentiable functions the two notions coincide, by the Cauchy-Goursat theorem.
\emph{``Entire''} functions are globally analytic, e.g. for all complex time arguments $t\in\mathbb{C}$.

Decomposing PDEs like \eqref{PDEw} into real and imaginary parts $w=u+\mi v$ leads to real reaction-diffusion systems of the special form
\begin{equation}
\label{PDEuv}
\begin{aligned}
   u_t=&u_{xx}+\Re f(u+\mi v)\,,   \\
   v_t=&v_{xx}+\Im f(u+\mi v)\,,  
\end{aligned}
\end{equation}
with equal diffusion coefficients for $u$ and $v$.
See for example \cite{Yanagida}.

When available, however, the complex viewpoint offers enormous advantages.
Examples are the residue theorem, analytic Poincaré linearization, or the fundamental fact that local flows or semiflows in different complex time directions commute.
We refer to \cite{Ilyashenko} for background on complex ODEs.
For a quick summary and further references see \cite{FiedlerShilnikov}.

Spatially homogeneous PDE solutions of \eqref{PDEw} lead back to ODE \eqref{ODE}, as the simplest interesting case.
Explicit solution of \eqref{ODE}, by standard separation of variables, makes the ODE look deceptively innocent and completely classical.
However, the precise details of blow-up $|w|\nearrow\infty$ in finite real time $t$, and complex circumnavigation of blow-up in the Masuda sense, are less classical topics.
See \cite{MuciRiem,DiPol,DiRat,Rousseaud=4} for some advanced literature, in this context, which we will return to in the discussion sections \ref{Degens}, \ref{Riems}.
Our main focus will be on nondegenerate Morse flows of \eqref{ODE} on the Riemann sphere $w\in\widehat{\mathbb{C}}=\mathbb{S}^2$.
Our PDE motivation, moreover, prods us to keep track of the real axis, as a distinguished direction of time.
Our approach is elementary, almost on textbook level, but nontrivial.
We hope to convince the reader that even a supposedly ``trivial'' ODE like \eqref{ODE} does hold some interest, today.

\subsection{Complex blow-up}\label{Blow-up}

Let $w(r,x)$ denote the maximal forward solution of \eqref{PDEw} for prescribed initial values $w_0(x)$ at time $r=0$, i.e. $w(0,x):=w_0(x)$.
\emph{Blow-up} usually means that the solution $w(r,x)$ becomes unbounded at some finite real positive time $r=r^*>0$, and in some suitable state space $X$ of $w_0$ and $w(t,\cdot)$.
For blow-up in parabolic real PDEs see the most diligent monograph \cite{Quittner}.
The term \emph{blow-down} is used to describe the same phenomenon in reverse, finite negative time.
More generally, we extend this terminology to include solutions of rational ODEs \eqref{ODE} which reach singularities $Q=0$ in finite time.

To be specific, let us consider the PDEs \eqref{PDEw} and \eqref{PDEpsi} on the real interval $0<x<\ell$ of length $\ell$, and under Neumann boundary conditions $u_x=\psi_x=0$, at $x=0,\ell$.
Forty years ago, pioneering work \cite{Masuda1,Masuda2} by Kyûya Masuda focused on orbits $w$ of the purely quadratic Fujita variant $w_r=w_{xx}+w^2$, in any space dimension.
Under Neumann boundary conditions, elliptic maximum principles imply that the spatially homogeneous solution $W=0$ is the only real equilibrium.
Real solutions $w(r,x)$ starting at initial profiles $w(0,x)=w_0(x)>0$, for $r=0$, then blow up in finite real time $0\leq r\nearrow r^*<\infty$.
For almost homogeneous real initial conditions $ w_0$\,, Masuda was then able to circumnavigate blow-up, at real time $r=r^*(w_0)$, via a sectorial detour venturing into complex time $t=r+\mi s$.
Notably, Masuda also cautioned that the detours via positive and negative imaginary parts $s$ agree in their real-time overlap after blow-up, if \emph{and only if} $w_0$ is spatially homogeneous.

For some subsequent literature in this PDE context, see for example \cite{COS,Yanagida,Stukediss,Stukearxiv,Jaquetteqp,JaquetteHet,JaquetteStuke,JaquetteMasuda}.
Classical studies observed the absence of any entire, i.e. globally holomorphic, solutions $t\mapsto w(t)$ for important classes of entire ODEs; see for example \cite{Rellich, Wittich1, Wittich2}. 
For blow-up in one-dimensional unstable manifolds see \cite{Ushiki1, Ushiki2}.
For further discussion of these references see \cite{FiedlerClaudia}.

One PDE feature of the real quadratic heat equation \eqref{PDEw} with $f(w)=w^2-1$ is the abundance of heteroclinic orbits, in real time $r$.
\emph{Heteroclinic} solutions $\Gamma(r)$ in real time $r\in\mathbb{R}$ connect time-independent equilibria $w=W_\pm$\,. In symbols
\begin{equation}
\label{leadsto}
\begin{aligned}
   \Gamma:\ W_- &\leadsto W_+\,; \quad \textrm{i.e.}\\
   \Gamma(r) &\rightarrow W_\pm\,, \quad\textrm{for}\  r\rightarrow\pm\infty\,,
\end{aligned}
\end{equation}
in suitable Banach spaces $X$ of solutions.
\emph{Connecting orbits}, \emph{traveling fronts} (or backs), or \emph{solitons} are other names for $\Gamma$.
Unless specified otherwise, explicitly, we subsume the \emph{homoclinic} case $W_+=W_-$ of non-constant $\Gamma(r)$ under the heteroclinic label.

It turns out that any globally bounded, real-valued solution $w(r,x), \ r\in\mathbb{R}$ of \eqref{PDEw} is in fact real heteroclinic, $w=\Gamma:W_-\leadsto W_+$\,, between distinct equilibria $W_\pm$.
This is due to a decreasing Lyapunov functional on real solutions; see \cite{Matano, Zelenyak, Lappicy}.
For further developments, see also \cite{MatanoLap,brfi88,brfi89,Galaktionov,firoSFB,firoFusco,FiedlerFila} and the many references there.

In \cite{FiedlerFila} we have embarked on a preliminary analysis of the complex heat-Schrödinger correspondence \eqref{psiw} between PDEs \eqref{PDEw} and \eqref{PDEpsi}, for quadratic Riccati nonlinearities like $f(w)=w^2-1$.
For the spatially homogeneous ODE case see section \ref{Ex2}.
The interval length $\ell$ serves as a bifurcation parameter.
Our results hold for most lengths $\ell$, with mostly just discrete sets of exceptions.
We establish that the extensions of heteroclinic orbits $w(r,\cdot)=\Gamma(r){:}\ W_-\leadsto W_+$\,, from real time $r$ to complex time $t=r+\mi s$, cannot be complex entire.
Failure to extend holomorphically to all $t\in\mathbb{C}$ entails blow-up and blow-down of the corresponding Schrödinger solutions \eqref{psiw}, for suitably fixed real $r_0$ and $0\leq \pm s\nearrow s^*<\infty$.

Assumptions actually restrict the above results to the case where the target is the unique stable equilibrium $W_+=-1$, and hence spatially homogeneous.
The source equilibrium $W_-$ is required to be unstable hyperbolic, and of unstable dimension (alias Morse index) $i(W_-)$ not exceeding 22.
Only in the homogeneous case $W_-=+1$ were we able to drop that Morse limitation, so far.
On the other hand, Schrödinger blow-up always occurs for $\Gamma$ which emanate from $W_-$ via a fast unstable manifold of dimension $d<1+i(W_-)/\sqrt{2}$.

These results are based on local Poincaré linearization \cite{ArnoldODE,Ilyashenko} in the complex analytic, finite-dimensional, (fast) unstable manifold of the heteroclinic source $W_-$\,.
The prerequisite spectral non-resonances impose the above restrictions, so far.
By Sturm-Liouville theory, the spectra of linearizations consist of simple real eigenvalues.
Let us fix the real time parameter $r_0$ sufficiently negative, so that $\Gamma(r_0)$ falls into the domain where Poincaré linearization rules.
For the Schrödinger solutions \eqref{psiw} associated to real heteroclinic orbits $w(r,\cdot)=\Gamma(r):W_-\leadsto W_+$\,, this implies global quasi-periodicity for all $s\in\mathbb{R}$, rather than blow-up for finite $s$.
For entire $\Gamma$, the resulting quasiperiodic Fourier coefficients would lead to contradictions.
Again, see \cite{FiedlerFila} for complete details, and a discussion of some further literature on PDEs in complex time.

For systems $\dot{w}=f(w),\ w\in X=\mathbb{C}^N$, of ODEs with complex entire nonlinearities $f$, our results are less restrictive.
Real-time ODE heteroclinic orbits $\Gamma:e_-\leadsto e_+$ between hyperbolic equilibria $f(e_\pm)=0$ then experience blow-up in imaginary time, provided the unstable spectrum of the linearization $f'(e_-)$ and the stable spectrum of $f'(e_+)$ are each nonresonant, separately.
See \cite{FiedlerClaudia}.

For some further references concerning the notoriously difficult Navier-Stokes and Euler PDEs, as well as complex $x$ and vorticity ODEs, see also \cite{LiSinai, Fasondini24}. 
Notably, \cite{LiSinai} establish blow-up of the 3d Navier-Stokes equation in real time, for certain \emph{complex} initial conditions.

\subsection{Morse flows on $\mathbb{S}^2$}\label{S2Morse}

Our main result will establish the equivalence of a generic class of (regularized) real-time flows of rational ODEs \eqref{ODE} on the Riemann sphere $\widehat{\mathbb{C}}=\mathbb{S}^2$ with the class of Morse flows on $\mathbb{S}^2$, much in the spirit of \cite{Peixoto1, Peixoto2}.
See section \ref{ResRat}, and in particular theorems \ref{thmRat2Portrait}, \ref{thmPortrait2Rat}, and corollary \ref{Rat=Morse} there.
We therefore introduce and explore this Morse class first, with a focus on real dimension 2.
See also \cite{Palis,PalisSmale,Sotomayor,PalisdeMelo, Hale} for the broader background of Morse-Smale systems,  structural stability, and bifurcation.

\subsubsection{The Morse properties}\label{DefMorse}

We abandon the complex viewpoint, for the moment, and consider real $C^1$ vector fields
\begin{equation}
\label{ODEFMorse}
\dot{w}=F(w)
\end{equation}
on the real $2$-sphere $w\in \mathbb{S}^2$.
In other words, $F$ is a $C^1$-section in the real tangent bundle $T\mathbb{S}^2$.
The resulting real-time flow is global.

An equilibrium $e$ is called \emph{hyperbolic} if the linearization $F'(e)$ does not possess purely imaginary (or zero) eigenvalues.
In two real dimensions, this means that the real part of both eigenvalues (counted with multiplicity) is strictly positive, strictly negative, or else both signs occur.
We then call $e$ a \emph{source}, \emph{sink}, or \emph{saddle}, respectively, and denote the unstable dimension, alias the \emph{Morse index}, by $i(e)=2,0$, or $1$.\,
A \emph{saddle-saddle connection} is a homoclinic or heteroclinic orbit $e_1\leadsto e_2$ between one or two saddle equilibria $e_1,e_2$\,.

\begin{defi}\label{defMorse}
We call $F$, or the global flow of ODE \eqref{ODEFMorse} on $\mathbb{S}^2$, \emph{Morse} if the following three conditions all hold:
\begin{enumerate}[(i)]
   \item equilibria are hyperbolic sources, saddles, or sinks;
   \item there are no nonstationary periodic orbits;
   \item there are no saddle-saddle connections.
\end{enumerate}
\end{defi}
In higher dimensions, condition \emph{(iii)} is replaced by the requirement that stable and unstable manifolds intersect transversely.

\subsubsection{The gradient-like property}\label{Gradlike}

We call ODE \eqref{ODEFMorse} \emph{gradient-like} if there exists a continuous \emph{Lyapunov function} $V:\mathbb{S}^2\rightarrow\mathbb{R}$ which decreases strictly along trajectories $w(t)$, except at time-independent equilibria $w(t)=e,\ F(e)=0$.
Standard terminology requires Morse vector fields to be gradient-like, up front.
On $\mathbb{S}^2$, the following ``converse'' proposition reconciles our definition \ref{defMorse} with standard terminology.

\begin{prop}\label{propgrad}
Morse vector fields on $\mathbb{S}^2$ are gradient-like.
\end{prop}

\begin{proof}
By property \emph{(i)} of definition \ref{defMorse}, the number of equilibria on $\mathbb{S}^2$ is finite.
Property \emph{(ii)} excludes periodic orbits.
The Poincaré-Bendixson theorem on the $2$-sphere therefore asserts that any $\boldsymbol{\alpha}$- and $\boldsymbol{\omega}$-limit sets can only consist of equilibria and heteroclinic orbits $e_1\leadsto e_2$ among them.
By absence of saddle-saddle connections \emph{(iii)}, 
\begin{equation}
\label{idrop}
e_1\leadsto e_2\ \Rightarrow\ i(e_1)>i(e_2)\,.
\end{equation}
By chain recurrence, any $\boldsymbol{\alpha}$- and $\boldsymbol{\omega}$-limit sets are therefore single equilibria. 
Moreover, the nonwandering set consists of equilibria, only.
This absence of nontrivial recurrence implies the gradient-like property; see \cite{Conley}.
\end{proof}

The following theorem, a straightforward corollary to Peixoto \cite{Peixoto1, Peixoto2}, shows that Morse flows are ``typical'', topologically, among gradient-like vector fields \eqref{ODEFMorse} on $\mathbb{S}^2$.
Moreover, they are ``robust'', in the following sense.
Two $C^1$ flows are called $C^0$ \emph{orbit equivalent}, if there exists a homeomorphism $\frak{H}$, in our case
\begin{equation}
\label{C0equiv}
\frak{H}{:}\ \mathbb{S}^2\rightarrow\mathbb{S}^2, 
\end{equation}
which maps orbits of the first flow onto orbits of the second flow, as sets.
\emph{Structural stability} requires, that arbitrary $C^1$-small perturbations of the vector field $F$ lead to $C^0$ orbit equivalent flows.
See also \cite{PalisdeMelo,Pilyugin}.
In the setting of structural stability, $\frak{H}$ is a perturbation of identity and therefore preserves the orientation of $\mathbb{S}^2$ and the time direction along orbits.
Depending on (nonlocal) context, general orbit equivalences $\frak{H}$ may, or may not, reverse orientation or time direction.

\begin{thm}\label{thmPeixoto}
In the $C^1$-topology of gradient-like vector fields $F$, Morse flows are generic and structurally stable.
\end{thm}

We recall that a property is \emph{generic}, if it holds for a countable intersection of subsets, each of which is open and dense.
By Baire's theorem, countable intersections of generic sets remain generic and hence dense, e.g., in complete metric spaces.

\subsubsection{Portraits and connection graphs}\label{Portraits}

To describe Morse flows \eqref{ODEFMorse} on $\mathbb{S}^2$, we  introduce two \emph{portraits} $\mathcal{C}^\pm$.
The portraits are finite undirected \emph{multi-graphs} on the $2$-sphere $\mathbb{S}^2=\widehat{\mathbb{C}}$.
Multiple loops to a single vertex and multiple edges between the same two vertices are admitted, in multi-graphs.

Let $E^+$ denote the set of source equilibria $i(e^+)=2$, and $E^-$ the sinks $i(e^-)=0$. 
Note $E^\pm\neq\emptyset$.
The possibly empty set $E'$ collects the saddles $i(e')=1$.
The sources $E^+$ provide the \emph{vertices} of the \emph{unstable portrait} $\mathcal{C}^+$. 
\emph{Edges} are the stable manifolds of saddles $e'\in E'$, including the saddles themselves.
Note how \emph{faces}, i.e. the connected components of $\mathbb{S}^2\setminus\mathcal{C}^+$, are the domains of attraction of the sinks.
In particular, forward time contracts each face to its unique defining sink.

Vertices of the \emph{stable portrait} $\mathcal{C}^-$ are the sinks $E^-$, and 
edges are the unstable manifolds of saddles $e'\in E'$, including the saddles.
Faces are now the domains of attraction of the original sources in $E^+$, in backward real time.
Time reversal $F\mapsto -F$ just swaps the roles of $\mathcal{C}^+$ and $\mathcal{C}^-$.

The pair $(\mathcal{C}^+, E^-)$ is a repellor-attractor pair in the sense of Conley index theory \cite{Conley}, just as the opposite pair $(E^+,\mathcal{C}^-)$ is.
Because the stable and unstable manifolds intersect transversely, at saddles $e'\in E'=\mathcal{C}^+\cap\mathcal{C}^-$, the stable and the unstable portraits are dual to each other on $\mathbb{S}^2$.

We summarize some of these properties.
Let the multi-graph $\mathcal{C}^+$ with $d^+$ vertices, $d'$ edges, and $d^-$ faces denote the unstable portrait of a Morse flow \eqref{ODEFMorse} on $\mathbb{S}^2$. 
Here $d^+$ counts the sources $e^+\in E^+\neq\emptyset$, alias the vertices of $\mathcal{C}^+$.
Similarly,  $d^-$ counts the sinks $e^-\in E^-\neq\emptyset$ alias faces, and $d'$ counts the saddles $e'\in E'$ alias edges.
Note $d^\pm\geq 1$.

\begin{prop}\label{propPortraits}
Assume $F$ in \eqref{ODEFMorse} is Morse.
Then the unstable portrait $\mathcal{C}^+$ is nonempty, finite, and connected.
By time reversal, analogous statements hold for the dual stable portrait $\mathcal{C}^-$, swapping exponents $\pm$ and sources with sinks.
Moreover,
\begin{equation}
\label{dd'Morse}
d:=d^++d^-=d'+2\geq 2\,.
\end{equation}
\end{prop}

\begin{proof}
Claim \eqref{dd'Morse} holds because the local Brouwer degrees $\textrm{sign}\,\det F'(e)$ of hyperbolic  equilibria $F(e)=0$ add up to the Euler characteristic 2 of $\mathbb{S}^2$, by the Lefschetz fixed point theorem.
With the remarks above, it only remains to prove connectivity of $\mathcal{C}^+$.

For multi-graphs on $\mathbb{S}^2$ with $m$ connected components, Euler's formula for the alternating counts of vertices, edges, and faces reads
\begin{equation}
\label{Eulerconn}
d^+-d'+d^-=m+1.
\end{equation}
Indeed, edge or loop contractions affect neither side.
Contraction of each connected component to a point, i.e. to $d^+=m,\ d'=0,\ d^-=1$, verifies \eqref{Eulerconn}.
Comparison of \eqref{Eulerconn} with \eqref{dd'Morse} then proves connectivity, $m=1$.
\end{proof} 

The only portraits without edges, $E'=\emptyset$, arise from quadratic polynomials like $\dot{w}=w^2-1$ on the Riemann sphere $\widehat{\mathbb{C}}=\mathbb{S}^2$.
See section \ref{Ex2} and figure \ref{fig4} for details.
Note $\mathcal{C}^\pm=\{\pm 1\}$ are then singletons; the portraits are self-dual.
For further examples see sections \ref{Ex3}, \ref{Exd}.

The directed and graded \emph{connection graph} $\mathcal{C}$ is slightly more detailed than the two portraits.
It consists of \emph{all} equilibria, as vertices.
Edges are the stable separatrices $e^+\leadsto e'$\,, from sources $e^+$  ($i=2$) to saddles $e'$ ($i=1$), and the unstable separatrices $e'\leadsto e^-$\,, from saddles to sinks $e^-$ ($i=0$).
In particular $\mathcal{C}$ is of rank -1: any edge decreases the vertex grading $i$ by 1; see \eqref{idrop}.
For many other applications of this concept, see for example \cite{Conley, brfi88, brfi89, Mischaikow, firoSFB, firoFusco, Yorke-a, Yorke-b} and the many references there.
A peculiarity of the present concept is planarity of $\mathcal{C}$, by construction.

Heteroclinic orbits from sources to sinks can be omitted here.
For Morse flows, indeed, the $\lambda$-Lemma \cite{PalisdeMelo} implies transitivity of heteroclinicity. 
In our case, conversely, the following \emph{cascading property} holds.
For any heteroclinic orbit $e^+\leadsto e^-$ from a source $e^+\in E^+$ to a sink $e^-\in E^-$ there exists an intermediate saddle $e'\in E'$ such that $e^+\leadsto e'\leadsto e^-$.
This assumes any saddles exist, i.e.
\begin{equation}
\label{d>2}
d'=d-2>0;
\end{equation}
see also \eqref{dd'}.
For counts $d^++d^-=d\geq 3$ of all sources and sinks, therefore, the saddles $E'=\mathcal{C}^+\cap\mathcal{C}^-$ and the resulting separatrix pieces of 
$\mathcal{C}^+\cup\mathcal{C}^-$ define the connection graph $\mathcal{C}$.
The cascading property does not hold for general gradient-like Morse systems. 
Note its violation for the antipodal flow $\dot w=-w$ on $\widehat{\mathbb{C}}=\mathbb{S}^2$; see section \ref{Ex2}.
Cascading has first been encountered in the context of scalar reaction-advection-diffusion PDEs; see \cite{Cascading} and section \ref{Sturm}.

Homeomorphic embeddings of connection graphs $\mathcal{C}$ imply $C^0$ \emph{orbit equivalence of the flows} \eqref{ODEFMorse}.
Indeed, consider any connected component $C$ of $\mathbb{S}^2\setminus\mathcal{C}$, complementary to the connection graph.
Because edges lower the Morse grading by 1, the boundary of $C$ contains only one source $e_1\in E$  and one sink $e_2\in E$. 
Moreover, all orbits in $C$ are heteroclinic from $e_1$ to $e_2$.
This allows us to extend any orbit equivalence on the boundaries, i.e. on $\mathcal{C}$, to the corresponding connected components $C$, e.g. via \hbox{small circles} centered at the sources and transverse to the flow; see Poincaré linearization \eqref{linequi}.

We summarize our discussion with a corollary which justifies the name ``portrait'' for the multi-graphs $\mathcal{C}^\pm$.

\begin{cor}\label{corPortraits}
For Morse flows on $\mathbb{S}^2$, any of the portraits $\mathcal{C}^\pm\subset\mathbb{S}^2$ implies its dual $\mathcal{C}^\mp$, by time reversal.
For $d=d'+2\geq3$, as assumed in \eqref{d>2}, we also obtain the transitive connection graph $\mathcal{C}=\mathcal{C}^+\cup\mathcal{C}^-$, and hence the global phase portrait of the gradient-like Morse flow \eqref{ODEFMorse}, up to $C^0$ orbit equivalence.
\end{cor}

\subsection{Complex ODE setting}\label{ODEset}

Real-time flows of polynomial or rational vector fields $f=P/Q$ on the Riemann sphere $w\in\widehat{\mathbb{C}}=\mathbb{S}^2$ have been addressed in complete generality.
See for example \cite{MuciRiem,DiPol,DiRat}.
The classification is quite involved.
We restrict attention to the explicit and generic notion of Morse flows, as in section \ref{Morse}.
The notion of nondegeneracy, in definition \ref{defnondeg}, is our analogue for complex flows.
This leads to the much easier classification by gradient-like Morse flows, in sections \ref{ResRat} and \ref{ResPol}, including some combinatorial counting for the latter.
The reduced scope also allows for comparatively elementary proofs, in sections \ref{Gen} -- \ref{PfPol}:
much of our phase plane analysis proceeds at almost textbook level.
We first collect some folklore on relations among Poincaré linearization, residues, periods, iso-periodic families of real-time periodic orbits, and heteroclinic cycles for scalar rational ODEs; see also \cite{Rousseaud=4,FiedlerShilnikov}.
For an expert presentation of much background on complex analysis of ODEs, we again recommend the beautiful book by Ilyashenko and Yakovenko \cite{Ilyashenko}.

For given initial condition $\mathbf{w}(0)=\mathbf{w}_0$ with ODE solution $\mathbf{w}(t)$, let $\Phi^t(\mathbf{w}_0):=\mathbf{w}(t)$ define the local solution flow of the (not necessarily scalar or polynomial) ODE $\dot{ \mathbf{w}}(t)=\mathbf{f}(\mathbf{w}(t))$, with locally holomorphic vector field $\mathbf{f}$.
Then the flow property 
\begin{equation}
\label{flow}
\Phi^{t_2}\circ\Phi^{t_1}= \Phi^{t_1+t_2}, \qquad \Phi^0=\mathrm{id},
\end{equation}
holds, for any argument $\mathbf{w}_0$ and, locally, for all $t_1,t_2\in\mathbb{C}$ such that the  closed complex parallelogram spanned  by $t_1,t_2$ in $\mathbb{C}$ is contained in the domain of existence of the local semiflow $\Phi^t(\mathbf{w}_0)$.
Moreover, and for example, the local flow $\Phi^t$ in real time $t=t_1=r$ commutes with the local flow in imaginary (Schrödinger) time $t=t_2=\mi s$.
This follows from the ODE and Cauchy's theorem.
Locally, in the scalar case and for fixed $t$, the ODE implies that the flow maps $\mathbf{w}_0\mapsto \Phi^t(\mathbf{w}_0)$ are biholomorphically conformal.
Indeed, the holomorphic local inverse of $\Phi^t$ is $\Phi^{-t}$.

\subsubsection{Basic notation and terminology}\label{Basic}

In our present ODE study we explore the real-time dynamics of ODE \eqref{ODE} for rational vector fields $f(w)=P(w)/Q(w)$, with coprime polynomials $P,Q$.
Like \cite{MuciRiem,DiRat}, in other words, we study $\dot{w}=f(w)$ on the Riemann sphere $w\in\widehat{\mathbb{C}} $ for holomorphic maps $f:\widehat{\mathbb{C}}\rightarrow\widehat{\mathbb{C}}$.
For some aspects of transcendental $f$ with essential singularities, see \cite{LebEss,MuciEss}.
For general background on complex analysis and Riemann surfaces, we refer to the textbooks \cite{Forster, Jost, Lamotke} and the standard references there

The group of biholomorphic automorphisms of the Riemann sphere $\widehat{\mathbb{C}}$, alias the \emph{Möbius group}, consists of the fractional linear automorphisms of $\widehat{\mathbb{C}}$, alias the \emph{Möbius transformations}
\begin{equation}
\label{Mobius}
w\mapsto z:=\frac{\frak{a}w+\frak{b}}{\frak{c}w+\frak{d}}\quad\mathrm{with}\quad
\begin{pmatrix}
    \frak{a}  &  \frak{b}  \\
    \frak{c}  &  \frak{d}
\end{pmatrix}
\in \mathrm{PGL}(2,\mathbb{C})\cong \mathrm{SL}(2,\mathbb{C})/\{\pm \mathrm{Id}\}\,.
\end{equation} 
The group is generated by \emph{shifts} $z=w+\frak{b}$, \emph{complex scalings} by nonzero $a=\rho\exp(\mi\theta)$, and the involutary \emph{inversion} $z=1/w$.

Möbius transformations leave the class of rational vector fields invariant, but may affect the degrees $d,d'$ of the polynomials $P,Q$.
In fact, inversion $z=1/w$ of \eqref{ODE} leads to
\begin{equation}
\label{ODEz}
\dot{z}=-z^2f(1/z)=-az^{2-d+d'}+\ldots\,.
\end{equation}
This reveals a hidden pole or equilibrium $w=e_*=\infty\in\widehat{\mathbb{C}}$, unless relation 
\begin{equation}
\label{dd'}
d'=d-2,
\end{equation} 
holds for the polynomial degrees.
Specifically we observe a zero or pole of multiplicity 
\begin{equation}
\label{dd'no}
|d'-(d-2)|\neq 0
\end{equation} 
at $w=\infty$, alias $z=0$, depending on the strictly positive or negative sign of $d'-(d-2)$.
See section \ref{Ex} for specific examples.
Only in case \eqref{dd'} \emph{all} zeros and poles on the Riemann sphere are finite, and visibly accounted for in \eqref{ODE}, by the bounded zero sets $E$ and $E'$ of $P$ and $Q$, respectively.

We assume that all zeros and poles on $\widehat{\mathbb{C}}$ are simple.
Möbius transformations then allow us to reduce any rational vector fields $f=P/Q$ to the case \eqref{dd'}, without loss of generality.
Moreover, all zeros $e_j$ of $P$ and $e_{-j'}$ of $Q$ remain simple and distinct, $1\leq j\leq d,\ 1\leq j'\leq d'$.
In other words, the disjoint sets $E:=\{e_j\,,\ j>0\}$ and $E':=\{e_{-j'}\,,\ j'>0\}$ enumerate all complex equilibria and poles in  $\widehat{\mathbb{C}}$.
Explicitly,
\begin{equation}
\label{PQ}
\begin{aligned}
P(w)    &:=a(w-e_1)\cdot\ldots\cdot(w-e_d)   \\
Q(w)    &:=(w-e_{-1})\cdot\ldots\cdot(w-e_{-d'})   \\
\end{aligned}
\end{equation}
with some constant scaling factor $0\neq a\in\mathbb{C}$.

\begin{figure}[t]
\centering \includegraphics[width=\textwidth]{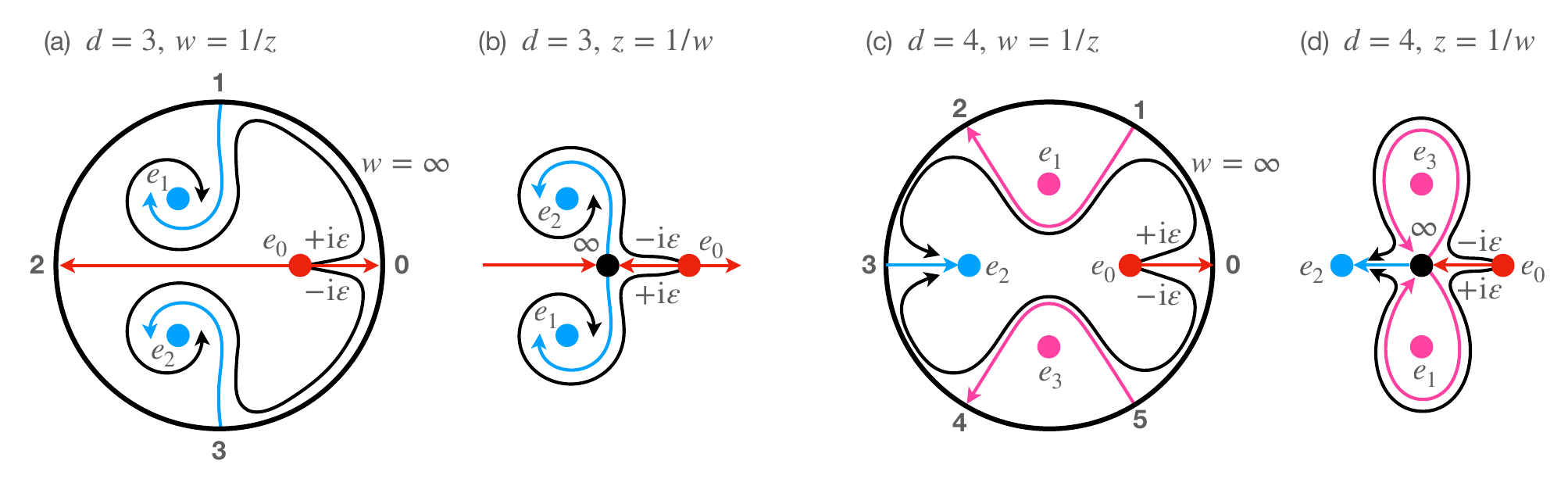}
\caption[Cyclotomic ODEs]{\emph{
Schematic phase portraits, in real time, of complex-valued ODEs \eqref{ODE}, \eqref{ODEPregw} -- \eqref{alpha} for cyclotomic polynomial vector fields $\dot{w}=f(w):=w^d-1$.
For $d=3$ see $w$ in (a), and $z=1/w$ in (b).
Similarly, $w$ in (c) and $z=1/w$ in (d) refer to $d=4$.
The invariant \emph{Poincaré circle} $\rho=|z|=0,\ \alpha\in\mathbb{S}^1$ of the Poincaré compactification at $w=\infty$ is marked black in (a), (c).
Interior equilibria $e_j\in\mathbb{D}$ are stable sinks (blue), Lyapunov centers (purple) surrounded by periodic orbits, or totally unstable sources (red).
Unstable blow-down separatrices (blue) emanate from $w=\infty$ at odd-labeled vertices $\mathbf{k}=\mathbf{1},\mathbf{3},\mathbf{5}$.  
Stable blow-up separatrices (red) run towards the even-labeled saddles $\mathbf{k}=\mathbf{0},\mathbf{2},\mathbf{4}$.
In (c), two blue-red pairs of interior saddle separatrices concur (purple), in each of the two interior saddle-saddle connections $\mathbf{1}\leadsto \mathbf{2}$ and $\mathbf{5}\leadsto \mathbf{4}$.
In (d), the two purple separatrices become homoclinic to $w=\infty$.\\
The black trajectories marked $\pm\mi\eps$ are complex perturbations of the red real blow-up separatrices $e_0\leadsto \mathbf{0}$, say with initial conditions $w(0)=2\pm\mi\eps$.
For $d=3$ in (a), they closely follow the heteroclinic chains $e_0\leadsto \mathbf{0}\leadsto \mathbf{1}\leadsto e_1$ and $e_0\leadsto \mathbf{0}\leadsto \mathbf{3}\leadsto e_2$, respectively.
The intermediate boundary connections $\mathbf{0}\leadsto \mathbf{1}$ and $\mathbf{0}\leadsto \mathbf{3}$ within the invariant boundary circle $\mathbb{S}^1$ concatenate initial red blow-up to terminal blue blow-down.
In the polar view (b), centered at $w=\infty,\ z=0$, the boundary parts contract into $z=0$.
Note the \emph{markedly distinct limits of the two perturbations}, for $\eps\searrow 0$, given by the two distinct remaining blow-up-down concatenations $e_0\leadsto \infty\leadsto e_1$ and $e_0\leadsto \infty\leadsto e_2$\,.\\
In case $d=4$ (d), the initial red blow-up and terminal blue blow-down limits of both perturbations $\pm\mi\eps$ coincide.
Each limit $e_0\leadsto \infty\leadsto\infty\leadsto e_2$ contains an additional interior blow-down-up separatrix $\infty\leadsto\infty$ (purple) which is homoclinic to $z=0$.
The two counter-rotating purple homoclinic separatrices, however, are markedly distinct: their lobes surround the distinct centers $e_1$ and $e_3$, respectively, in opposite direction.
In (c), this is manifest by the two purple interior saddle-saddle connections $\mathbf{1}\leadsto \mathbf{2}$ and $\mathbf{5}\leadsto \mathbf{4}$.
Together with the heteroclinic boundary connections within the black invariant boundary circle $\mathbb{S}^1$, which run between the same saddles but in opposite direction, we obtain two heteroclinic cycles. 
Their interiors are properly foliated by families of nested, synchronously iso-periodic orbits of minimal periods $\mp\pi/2$ around the Lyapunov centers $e_1$ and $e_3$\,; see section \ref{ODEper}.
}}
\label{fig1}
\end{figure}

As a consequence of degree relation \eqref{dd'} we may rescale ODE \eqref{ODE}, replacing $a'w$ by $w$ and $e_{\pm j}$ by $\tfrac{1}{a'}e_{\pm j}$.
Whenever convenient, we may therefore normalize $a$ to
\begin{equation}
\label{a=1}
a=1
\end{equation} 
or any other nonzero complex value.
The purely polynomial case $Q=1$ studied in \cite{FiedlerShilnikov} and in section \ref{ResPol} below, in contrast, introduces a single pole of multiplicity $d'=d-2$, located at $w=\infty$, for $d>2$.
Examples of polynomial degrees $d=3,4$, and $5$ are illustrated in figures \ref{fig1}, and \ref{fig6}.
The ``anti-polynomial'' case $P=1$, analogously, leads to a single zero of multiplicity $d:=d'+2$, located at $w=\infty$, for any $d'\geq 1$.
See sections \ref{ResAntipol}, \ref{PfAntipol}, and figure \ref{fig9}.
Only the Masuda case $d=2$ of quadratic $f=P$ possesses neither pole nor equilibrium at infinity.
See section \ref{Ex2} and figure \ref{fig4}.

We may always regularize the ODE \eqref{ODE} by noting 
\begin{equation}
\label{ODEPQbar}
f=P/Q=|Q|^{-2}P\overline{Q}\,.
\end{equation}
Multiplication of \eqref{ODE} with an Euler multiplier $(1+|w|^2)^{-(d-2)}|Q|^2\geq0$, and rescaling real time to absorb the Euler multiplier on the left, we obtain the \emph{regularized ODE}
\begin{equation}\label{ODEreg}
\dot{w}=f_\mathrm{reg}(w):=(1+|w|^2)^{-(d-2)}P(w)\overline{Q(w)}\,.
\end{equation}
Although we lose complex meromorphy, the regularization preserves real-time orbits, except at poles $Q(e_{-j'})=0$.
Poles of $f$ become equilibria of the regularized ODE \eqref{ODEreg}.

The resulting regularized compact flow is smooth and global on the whole Riemann sphere $\widehat{\mathbb{C}}=\mathbb{S}^2$.
Indeed, in coordinates $z=1/w$ we obtain 
\begin{equation}\label{ODEregz}
\dot{z}=-z^2 f_\mathrm{reg}(1/z)=(1+|z|^2)^{-(d-2)}\widetilde{P}(z)\overline{\widetilde{Q}(z)}\,.
\end{equation}
Here we use the polynomial transformations 
\begin{equation}
\label{PQtilde}
\begin{aligned}
\widetilde{P}(z)&:=-z^{d\phantom{'}}P(1/z)\,,\\
\widetilde{Q}(z)&:=\phantom{-}z^{d'}Q(1/z)\,,
\end{aligned}
\end{equation}
in the notation \eqref{PQ} with $d'=d-2$ of \eqref{dd'}.
In particular $\widetilde{P}(0)=-a,\ \widetilde{Q}(0)=1$, and $\dot{z}=-a+\ldots$ near $z=0,\ w=\infty$.

Since $F:=f_\mathrm{reg}$ is a smooth vector field \eqref{ODEFMorse} on $\mathbb{S}^2=\widehat{\mathbb{C}}$, the general Morse terminology of section \ref{S2Morse} applies to the regularization \eqref{ODEreg} of rational vector fields \eqref{ODE}.
Solutions $w(t)\in\widehat{\mathbb{C}}\setminus E'$ of ODE \eqref{ODE} which reach a pole $e_{-j'}\in E'$ in finite real \emph{blow-up time} $0<t\nearrow t^*<\infty$, previously called stable separatrices, for example, now become \emph{blow-up solutions}.
Their unstable separatrix cousins, which reach poles in negative time $0>t\searrow t_*>-\infty$ indicate \emph{blow-down}.
For further detail, we investigate the linearization at equilibria next.

\emph{Linearization} $\dot{w}=f'(e_j)w$ at a simple equilibrium $e_j\in E$ is given by
\begin{equation}
\label{linequi}
f'(e_j)=P_{j}(e_j)/Q(e_j)\neq 0,\quad\textrm{where}\quad P_{j}(w):= \prod_{0<k\neq j} (w-e_k)\,.
\end{equation}
Since $\dim_\mathbb{C} w = 1$, the nonzero linearization $f'(e_j)$ is trivially nonresonant.
By Poincaré linearization, the flow is then locally biholomorphically equivalent to its linearization; see \cite{Ilyashenko}.
Depending on the real part $\Re f'(e_j)$, several cases arise; see figure \ref{fig1}.
The equilibrium $e_j$ is \emph{hyperbolic} if $\Re f'(e_j)\neq0$.
For $\Re f'(e_j)<0$, we obtain a \emph{sink} (blue).
\emph{Sources} (red) have $\Re f'(e_j)>0$.
For purely imaginary $\Re f'(e_j)=0$ with nonzero imaginary part, we obtain a \emph{Lyapunov center} (purple).
In other words, a neighborhood of $e_j\in\mathbb{C}$ is foliated by periodic orbits $\gamma$ nested around $e_j$\,, all of the same minimal period $|T|>0$; see \eqref{Hopfperiod} below, and $e_1,e_3$ in figure \ref{fig1}(c),(d).

Next, we address the local behavior of ODE \eqref{ODE}, alias \eqref{ODEPQbar}, at simple poles $e_{-j'}$\,.
We linearize the regularization \eqref{ODEreg} at $e_{-j'}$ and obtain
\begin{equation}
\label{linpolew}
\dot{w}=f_\mathrm{reg}'(e_{-j'})\overline{w}\,.
\end{equation}
Here the anti-holomorphic complex conjugate $\overline{w}$ appears, because $Q(e_{-j'})=0$ forces us to linearize the factor $\overline{Q}$ in \eqref{ODEreg}:
\begin{equation}
\label{linpole}
\begin{aligned}
   f_\mathrm{reg}'(e_{-j'})&\phantom{:}=P(e_{-j'})\overline{Q}_{-j'}(e_{-j'})\neq 0,\quad\textrm{where}\\
    \overline{Q}_{-j'}(w)&:= (1+|w|^2)^{-(d-2)}\prod_{0<k'\neq j'} (\overline{w}-\overline{e_{-k'}})\,. \end{aligned}
\end{equation}
In components $w=u+\mi v$ and for any complex number $f_\mathrm{reg}'(e_{-j'})=\alpha+\mi\beta\neq 0$, the anti-holomorphic linearization \eqref{linpolew} takes the symmetric form
\begin{equation}
\label{linpoleuv}
\begin{pmatrix}
      \dot{u}    \\
      \dot{v}  
\end{pmatrix}=
\begin{pmatrix}
 \alpha     &  \beta  \\
 \beta     &  -\alpha
\end{pmatrix}
\begin{pmatrix}
      u    \\
      v  
\end{pmatrix}
\end{equation}
with negative Jacobian determinant $\det=-(\alpha^2+\beta^2)<0$.
In particular, any simple pole $e_{-j'}$ is a hyperbolic saddle point of the regularized flow \eqref{ODEreg}, in real time $t$.
Any saddle $e_{-j'}$ comes with unique, real analytic stable and unstable manifolds which intersect orthogonally at $e_{-j'}$\,.
Heteroclinic \emph{separatrices}, as in section \ref{S2Morse}, are the two half-branches of each, when we omit the saddle itself.

This shows that stable (red) and unstable (blue) separatrices are precisely the blow-up and blow-down orbits, respectively, which reach pole singularities $e_{-j'}$ of the vector field $f$ in finite real time.
Indeed, implicit solutions $w=w(t)$ of the initial value problem $w(0)=w_0$ for the ODE \eqref{ODE} are given by explicit separation of variables as
\begin{equation}
\label{sov}
t \,=\, \int_{w_0}^{w(t)}\omega \,=\, \int_{w_0}^{w(t)}Q(w)/P(w)\,dw \,=\, c \,+ \, \sum_{j=1}^d\,\eta_j \log(w(t)-e_j)\,,
\end{equation}
with some integration constant $c\in\mathbb{C}$.
Here we use the abbreviations
\begin{equation}\label{omega,eta}
\omega:=dw/f(w)\quad \textrm{and} \quad \eta_j:=1/f'(e_j)
\end{equation}
for the 1-form integrand $dw/f$ and the residues $\eta=1/f'$  at the equilibria $e_j$\,, alias the coefficients of the partial fraction decomposition of $f=Q/P$.
Note the constraint
\begin{equation}
\label{sumeta0}
\sum_{j=1}^d \eta_j \ =\  \frac{1}{2\pi\mi} \int_{|w|=\rho} \omega \ =\  0\,,
\end{equation}
due to the residue theorem for large circles $\textrm{max}_j\,|e_j|<\rho\nearrow\infty$; see \eqref{ODEz} for $d'\leq d-2$.

The integration constant $c$ in \eqref{sov} depends on the path $\Gamma$ of $w$ for the integration from $w_0$ to $w(t)$, in the punctured Riemann sphere
\begin{equation}\label{C-E}
w\in\widehat{\mathbb{C}}\setminus E\,.
\end{equation}
More precisely, the residue theorem (alias the multi-valued complex logarithm) implies uniqueness of the integral, up to an integration constant
\begin{equation}
\label{calP}
c\in\mathcal{P}\,:=\,2\pi\mi\,\langle \eta_1,\ldots,\eta_d\rangle_\mathbb{Z}\,.
\end{equation}
We call $\mathcal{P}$ the \emph{period module} over $\mathbb{Z}$ for the rational form $\omega=dw/f=Q/P\, dw$ of \eqref{ODE}. 

Period modules aside, we conclude that integration up to the poles $e_{-j'}$\,, i.e. up to the zeros of the numerator $Q(w)$ of the differential form $\omega$, poses no singularity obstacle.
In particular, any stable or unstable separatrix indeed reaches its singularity $e_{-j'}$ in finite original blow-up time $t^*>0$ or blow-down time $t_*<0$, respectively.
Locally at a pole, say at $e'=0\in E'$, we can expand the integral \eqref{sov} for $w=w(t)$ as 
\begin{equation}
\label{we'}
t= \tfrac{1}{2}\,Q'(e')/P(e')\,w^2+...\,.
\end{equation}
Locally, we may view this as the classical Riemann surface covering map $w\mapsto t$ of the square root.
The two local preimages $w$ of $t>0$ then define the real analytic unstable separatrices (blue) with blow-down from $t\searrow t^*=0$.
Analogously, $t<0$ identifies the red stable separatrices which blow up at $t\nearrow t^*=0$.

Summarizing, we can now compare the original "flow" of ODE \eqref{ODE} and the flow of the regularization \eqref{ODEreg} on the compact Riemann sphere $\widehat{\mathbb{C}}$, from a global point of view.
Real-time orbits coincide, except at saddles $e_{-j'}$\,.
At saddles, only stable separatrices (red) proceed towards the singularities $e_{-j'}$ of $f$, taking finite original blow-up time in \eqref{ODE}.
Afterwards, they may take off along either of the two unstable separatrix branches (blue), non-uniquely and in finite blow-down time.
This contrasts sharply with the quadratic Masuda case $d=2,\ d'=0$, or $d-d'=2$ more generally, where the vector field $f=P/Q$ of the original complex ODE \eqref{ODE} is neither vanishing nor singular at $w=\infty$.
All other real-time orbits remain regular, in forward time.

\subsubsection{Periodic orbits and heteroclinic cycles}\label{ODEper}

Let $\gamma\subset\widehat{\mathbb{C}}\setminus (E\cup E')$ denote any nonstationary real-time periodic orbit of \eqref{ODE} with minimal period $|T|>0$. 
Let $j\in J$ enumerate the equilibria ``inside'' the closed Jordan curve $\gamma$.
Then explicit integration \eqref{sov}, \eqref{omega,eta} along $\gamma$ and the residue theorem imply
\begin{equation}
\label{T}
0\neq T = \int_\gamma \omega = 2\pi\mi \sum_{j\in J} \eta_j\,.
\end{equation}
The sign of $T\in\mathbb{R}$ indicates the clockwise or anti-clockwise time direction of the periodic solution $\gamma$, with respect to the chosen interior.
Swapping interior and exterior swaps $\emptyset\neq J\subsetneq\{1,\ldots,d\}$ with its nonempty complement $J^c$ and reverses the sign of $T$, in accordance with \eqref{sumeta0}.
Moreover \eqref{T} implies
\begin{equation}
\label{sumetaJ0}
\Re \sum_{j\in J} \eta_j \ =\  \Re \sum_{j\in J^c} \eta_j \ = 0\,.
\end{equation}
Let $d_J$ and $d_J'$ count the equilibria and poles inside the periodic orbit $\gamma$ of ODE \eqref{ODE}, i.e. the zeros and poles of $f=P/Q$, respectively.
Then, equivalently, $d_J$ and $-d_J'$ sum the Brouwer degrees $+1$ and $-1$ of the corresponding zeros of $f_\mathrm{reg}=(1+|w|^2)^{-(d-1)}P\overline{Q}$ in the regularization \eqref{ODEreg}, inside $\gamma$.
Since the vector field along the periodic orbit possesses winding number $+1$, in either case and orientation, we conclude
\begin{equation}
\label{Brouwer}
d_J'=d_J-1,
\end{equation}
and likewise for the exterior $J^c$.
In particular, $d_J\geq 1$ and $J, J^c$ are both nonempty -- as follows from standard Poincaré-Bendixson theory \cite{Hartman}, or directly from \eqref{T}.
Summing over the interior and exterior recovers degree relations \eqref{dd'Morse}, \eqref{dd'}.
For the purely polynomial case $d_J'=0,\ d_J=1$, we obtain a unique equilibrium $J=\{j\}$ in the bounded interior component, i.e. a Lyapunov center $e_j$\,; see also Lemma 4.1 in \cite{FiedlerShilnikov}.

Real-time periodic orbits $\gamma$ of ODE \eqref{ODE} occur within \emph{local foliations by nested iso-periodic} (alias \emph{isochronous}) \emph{orbits}.
To see this, we just invoke that flows in real time $t_1=r$ and in imaginary time $t_2=\mi s$ commute; see \eqref{flow}.
Since the conformal map $t=r+\mi s \mapsto \Phi^t(w_0)$ preserves angles, the local flow $\Phi^{\mi s}$ applied to the real-time periodic orbit $\gamma$ provides the iso-periodic foliation, as claimed.
In particular, the local iso-periodic foliation excludes \emph{limit cycles}, i.e. periodic orbits $\gamma$ which are $\boldsymbol{\alpha}$- or $\boldsymbol{\omega}$-limit sets of initial conditions $w_0\not\in\gamma$.

\emph{Maximal iso-periodic families} of $\gamma\subset\widehat{\mathbb{C}}\setminus (E\cup E')$ are connected components of local iso-periodic foliations.
They extend, as open subsets of $\widehat{\mathbb{C}}$, until they hit boundary points in $E\cup E'$.
Sources and sinks $e_j\in E$ are excluded, because they attract neighborhoods in backward or forward time; see \eqref{linequi}.
Lyapunov centers in $e_j\in E$ only arise as centers of punctured iso-periodic disks, with signed periods
\begin{equation}
\label{Hopfperiod}
T=2\pi\mi \eta_j = 2\pi\mi/f'(e_j)\in\mathcal{P}\,;
\end{equation}
compare \eqref{linequi}, \eqref{sov}, \eqref{calP}, and \eqref{T}.
See also lemma 4.1 in \cite{FiedlerShilnikov} for further details on periodic orbits, in the holomorphic case.

Following \cite{Sotomayor}, in spirit if not by the letter, the only remaining option for the boundary of a maximal iso-periodic family are certain poles $e_{-j'}\in E'$ and their separatrices.
See \cite{DiRat,Rousseaud=4} for very detailed accounts.
We may encounter \emph{heteroclinic cycles of length} $\ell$ such that
\begin{equation}
\label{hetcycle}
e_{-j'_1}\,\leadsto\, \ldots \,\leadsto\, e_{-j'_\ell} \,\leadsto\, e_{-j'_{\ell+1}} =\, e_{-j'_1} \,.
\end{equation}
See \eqref{leadsto} for notation.
A saddle-saddle connection $e_{-j'_1} \,\leadsto\, e_{-j'_2}$\,, for example, means that one blue unstable separatrix of $e_{-j'_1}$ coincides with one red stable separatrix of $e_{-j'_2}$\,, globally.
The special case of a homoclinic orbit $e_{-j'} \,\leadsto\, e_{-j'}$ constitutes a ``heteroclinic cycle'' of length $\ell=1$.
See figure \ref{fig1} for examples.

Heteroclinic cycles decompose into one or several closed Jordan curves $\Gamma$ which are disjoint, except at their saddle equilibria $e_{-j'}$.
Let $j\in J$ connect the equilibria in the interior component of any $\Gamma$.
Since each $\Gamma$ comes with an associated finite real time $T$, in the original ODE \eqref{ODE}, the residue constraint \eqref{sumetaJ0} remains valid for the equilibria $e_j$\,.

\subsubsection{Saddle-saddle connections}\label{Saddlehet}

Motivated by the heteroclinic cycles above, we now address saddle-saddle connections
\begin{equation}
\label{hetsingle}
\Gamma:\ e_{-j'_1}\,\leadsto\,e_{-j'_2} \neq e_{-j'_1}
\end{equation}
which are not homoclinic.
According to \eqref{sov}, the finite real time $T=T(\Gamma)>0$ which the solution $\Gamma$ of ODE \eqref{ODE} spends from $e_{-j'_1}$ to $e_{-j'_2}$ is given explicitly by
\begin{equation}
\label{sovpath}
0<T(\Gamma) \,=\, \int_\Gamma \omega \,=\, c \,+ \, \int_{\Gamma_0} \omega\,.
\end{equation}
Here $\Gamma_0\subset\widehat{\mathbb{C}}\setminus E$ is any path from $e_{-j'_1}$ to $e_{-j'_2}$, which need not coincide with the saddle heteroclinic $\Gamma$, and $c\in\mathcal{P}$ is a period which depends on $\Gamma_0$, for given $\Gamma$.
In particular,
\begin{equation}
\label{hetTP}
   \int_{\Gamma_0} \omega\ \in\ T+\mathcal{P}\,. 
\end{equation}
In the homoclinic case $e_{-j'_2} = e_{-j'_1}$\,, this simplifies to
\begin{equation}
\label{homTP}
0<T\ \in\ \mathcal{P}\,.
\end{equation}

\subsubsection{Nondegeneracy, Morse structure, and connection graphs}\label{Morse}

For rational ODEs \eqref{ODE} we henceforth assume the following nondegeneracy condition.

\begin{defi}\label{defnondeg}
We call the rational vector field $f=P/Q$ \emph{nondegenerate} if the following conditions hold:
\begin{enumerate}[(i)]
  \item all zeros $e_j$ and $e_{-j'}$ of the polynomials $P$ and $Q$ are simple, their zero sets $E$ and $E'$
  	 are disjoint, and their degrees $d$ and $d'$ satisfy $d'=d-2\geq 0$; see \eqref{dd'}, \eqref{PQ};
  \item none of the equilibria $e_j\in E$ lies on a straight line $\Gamma_0$\,, of real dimension 1, through any pair of poles 
  	$e_{-j'}\in E'$;
  \item the residues $\eta_j=1/f'(e_j)$ satisfy
	\begin{equation}
	\label{sumetaJno}
	\Re \sum_{j\in J} \eta_j \ \neq 0\,,
	\end{equation}  
	for any nontrivial subset $\emptyset\neq J\subsetneq\{1,\ldots,d\}$; compare \eqref{omega,eta}, 
	\eqref{sumeta0}, and \eqref{sumetaJ0};
   \item integration of $\omega=dw/f$ along any straight line $\Gamma_0$ between poles satisfies
     \begin{equation}
     \label{hetTPno}
	 \int_{\Gamma_0} \omega\ \not\in\ \mathbb{R}+\mathcal{P}\,.
     \end{equation}	
\end{enumerate}
\end{defi}

In theorem \ref{thmRatGen} and section \ref{Gen} we will show that nondegeneracy is a generic property in the very restrictive classes of rational vector fields $f=P/Q$ with degrees $d'=d-2$.
As an easy consequence of nondegeneracy, the smooth regularized flow on the Riemann sphere is Morse and gradient-like; see definition \ref{defMorse} and proposition \ref{propgrad}.

\begin{cor}\label{corMorse}
Assume $f=P/Q$ is nondegenerate. Then the regularized flow \eqref{ODEreg} is Morse, in real time, and hence gradient-like.
\end{cor}

\begin{proof}
We have to verify requirements \emph{(i)--(iii)} of Morse definition \ref{defMorse}, based on properties \emph{(i)--(iv)} of nondegeneracy assumption \ref{defnondeg}.
Regularization \eqref{ODEreg} works by nondegeneracy property \emph{(i)}, and provides the nonlinearity $F$ of ODE \eqref{ODEFMorse}.
Property \emph{(iii)} excludes Lyapunov centers $e_j$\,; see \eqref{sumetaJno} for $J=\{j\}$.
This leaves us with linearly non\-degenerate sources, sinks and, by \eqref{linpolew}--\eqref{linpoleuv}, saddles -- verifying Morse requirement \ref{defMorse}\emph{(i)}.
To prove \ref{defMorse}\emph{(ii)}, i.e. absence of periodic orbits, just compare nondegeneracy \eqref{sumetaJno} with \eqref{sumetaJ0}.
In view of \eqref{homTP}, the same argument shows absence of homoclinic orbits.
To prove \ref{defMorse}\emph{(iii)} on $\mathbb{S}^2$, it remains to exclude saddle-saddle connections.
This follows from nondegeneracy properties \emph{(ii)} and \emph{(iv)}; compare \eqref{hetTP} and nondegeneracy \eqref{hetTPno}.
This proves Morse requirement \ref{defMorse}\emph{(iii)}, and the corollary.
\end{proof}

Therefore our main results aim for a classification of connection graphs $\mathcal{C}$ up to homeomorphic embeddings or, equivalently, of portraits $\mathcal{C}^\pm$ or, again equivalently, for a classification of regularized nondegenerate rational flows \eqref{ODEreg} up to $C^0$ orbit equivalence.

\subsubsection{Glossary and color coding}\label{Color}

We summarize this introduction with a brief glossary of terms and color codings concerning the real-time dynamics of rational ODEs \eqref{ODE}, \eqref{ODEreg} on the $2$-sphere $\mathbb{S}^2=\widehat{\mathbb{C}}$.
See figures \ref{fig1}, \ref{fig6}, \ref{fig7} for an illustration of all terms.
The \emph{equilibria} $e_j\in E,\ j=1,\ldots,d$ possess nonzero complex linearizations $f'(e_j)=1/\eta_j$ and are called (blue) \emph{sinks}, (purple) \emph{Lyapunov centers}, and (red) \emph{sources} in case $\Re f'(e_j)$ is strictly negative, zero, or strictly positive, respectively.
Heteroclinic orbits $e_j\leadsto e_k$ (black) from sources $e_j$ to sinks $e_k$ foliate open regions and are called \emph{source/sink heteroclinics}.

The $d'=d-2\geq 1$ \emph{simple poles} $e_{-j'}\in E',\ j'=1,\ldots,d'$ are indicated by black circles.
They are hyperbolic saddles.
Their \emph{stable separatrices} (red), i.e. the half branches of their stable manifolds, 
are called  \emph{blow-up orbits}.
Their \emph{unstable separatrix} counterparts (blue) are called \emph{blow-down orbits}.
If two such separatrices happen to coincide, by non-generic coincidence, they form a (purple) saddle-saddle heteroclinic orbit or homoclinic loop which we indiscriminately call a \emph{saddle-saddle connection} or \emph{blow-down-up} orbit.
In all other cases, blue blow-down unstable separatrices are heteroclinic orbits towards sinks, and red blow-up separatrices are heteroclinic orbits emanating from sources.

In the \emph{nondegenerate} case of definition \ref{defnondeg}, purple objects do not occur.
The red graph of the \emph{unstable portrait} $\mathcal{C}^+$ consists of the red unstable sources, as vertices, and the nonoriented stable separatrix pairs, as red edges, marked by their generating saddles (black circles).
In the rational complex context, where the red separatrices coincide with the blow-up orbits, we also call the unstable portrait $\mathcal{C}^+$ the \emph{blow-up portrait}.
The blue \emph{stable portrait}  $\mathcal{C}^-$ consists of blue sink vertices and the nonoriented unstable separatrix pairs, as blue edges, marked by the same saddles (black circles).
Since the blue separatrices coincide with the blow-down orbits, we also call the stable portrait $\mathcal{C}^+$ the \emph{blow-down portrait}.
Regular trajectories passing through $w=\infty$ are not counted among the blow-up or blow-down phenomena because they may be circumnavigated, in complex time, by a local Masuda detour.
Indeed, their ``fake'' blow-up disappears under Möbius transformations.

Blow-up and blow-down portraits are dual to each other, on $\widehat{\mathbb{C}}=\mathbb{S}^2$.
Their union $\mathcal{C}^+\cup \mathcal{C}^-$, decomposed into blow-up and blow-down separatrices, defines the planar \emph{connection graph} $\mathcal{C}\subset\widehat{\mathbb{C}}$, with vertex grading by Morse indices $i$. 
All edges are of rank -1, and can be oriented in the direction of real time, alias strictly decreasing Morse index.

All figures are understood up to homeomorphic embeddings of isomorphic graphs or, equivalently, global $C^0$ orbit equivalence.

\subsection{Outline} \label{Out}

We proceed as follows.
In section \ref{ResRat} we state our main equivalence result, theorems \ref{thmRatGen} and \ref{thmPortrait2Rat}, between Morse flows on $\mathbb{S}^2$ and regularized global flows \eqref{ODEreg} of nondegenerate rational ODEs \eqref{ODE} on the Riemann sphere $\widehat{\mathbb{C}}$.
In preparation for proofs, we include some nondegeneracy and genericity results for rational flows on $\widehat{\mathbb{C}}$.
In \cite{FiedlerShilnikov} theorem 1.4, a somewhat analogous result was announced for polynomials.
Section \ref{ResPol} states that polynomial equivalence result as theorems \ref{Classthm} and \ref{thmTree2Pol}.
Section \ref{Ex} illustrates the results by several examples.
Rational genericity of the Morse property is addressed in section \ref{Gen}.
Section \ref{PfthmRat} proves rational realization theorem \ref{thmPortrait2Rat}.
The polynomial variant, theorem \ref{thmTree2Pol}, is proved in section \ref{PfPol}.
The curious case of anti-holomorphic polynomials turns out to be, simultaneously, gradient-like \emph{and} Hamiltonian on the Riemann sphere.
See section \ref{ResAntipol} and the proof of theorem \ref{nctree=antipol} in section \ref{PfAntipol}.
Section \ref{Dis} concludes with a discussion of some higher-dimensional and philosophical aspects.

\subsection{Acknowledgment}\label{Ack} 
This paper is dedicated to the memory of Professor Masaya Yamaguti, renowned for his inspired and visionary influence on mathematics in Japan, and far beyond. 
Unforgettable is my encounter with him at the International Congress of Industrial and Applied Mathematics (ICIAM) at Hamburg, 1995.
Approaching him as the nobody I was, I immediately felt unreservedly and warmly welcome. 
In his mild-mannered way, with eyes of alert curiosity, he asked me: ``so, who do you know in Japan?'' 
And with almost every name I mentioned -- Japanese mathematicians of wide and international acclaim, some of them now editors of this volume -- he kindly nodded and remarked: ``oh yes, my student'', or ``ahso, a student of my student'', and on and on.
I am much indebted to the editors of this volume, for their initiative and long-lasting friendship and support.
What an admirable and flourishing academic family tree!

Vassili Gelfreich very patiently explained his profound and motivating work on exponential splitting of homoclinic orbits. Anatoly Neishtadt influenced me with lucid conversations on adiabatic elimination. 
Tibor Szabó pointed me at the combinatorics of nc-trees, with instant precision.
Karsten Matthies, Carlos Rocha, Jürgen Scheurle, Hannes Stuke, Nicolas Vassena, and the late colleagues Marek Fila and Claudia Wulff, have provided lasting interest and motivation in real and complex times.
The concluding remarks rely much on gentle and forgiving initiation, by Jia-Yuan Dai, Yuya Tokuta, Biao Xiang, and many others including DeepL, to the unfathomable realms of Asian thought .

\section{Main Results}\label{Res}

\subsection{Rational ODEs}\label{ResRat}

We can now formulate our main results.
Throughout we consider $d\geq2,\ d'=d-2$.
See section \ref{Ex2} for the polynomial quadratic case $d=2$.
For nondegeneracy of rational ODEs \eqref{ODE} and their regularizations \eqref{ODEreg} on the $2$-sphere $\mathbb{S}^2=\widehat{\mathbb{C}}$ see definition \ref{defnondeg}.
By corollary \ref{corMorse}, the nondegeneracy assumption implies the Morse property.
By proposition \ref{propgrad}, the flows are then also gradient-like.
As in sections \ref{Morse} and \ref{Color}, we associate the mutually dual red blow-up and blue blow-down portraits $\mathcal{C}^\pm$ to these flows. 

The Peixoto theorem \ref{thmPeixoto} established genericity of $C^1$ Morse flows on $\mathbb{S}^2$.
Our present setting of (regularized) rational vector fields on the Riemann sphere, however, is much more restrictive.
In section \ref{Gen} we prove the following genericity result.

\begin{thm}\label{thmRatGen}
Nondegeneracy holds for any fixed complex coefficient $a\neq0$ and a generic, in particular dense, subset $\mathcal{E}\subset\mathbb{C}^{2d-2}$ of configurations of $d$ simple zeros and $d'=d-2$ simple poles of $f=P/Q$; see \eqref{PQ}.
\end{thm} 

For nondegenerate rational ODEs, the description of gradient-like Morse flows in section \ref{S2Morse} applies, as follows.
See section \ref{Color} for terminology and color codes.

\begin{thm}\label{thmRat2Portrait}
Any nondegenerate flow \eqref{ODE}, regularized by \eqref{ODEreg}, is gradient-like Morse.
Therefore it is characterized, up to $C^0$ orbit equivalence, by its blow-up  portrait $\mathcal{C}^+\subset\mathbb{S}^2=\widehat{\mathbb{C}}$ or its dual blow-down portrait $\mathcal{C}^-$.
Proposition \ref{propPortraits} applies, with $d=d^++d^-$ counting red sources and blue sinks, and $d'=d-2$ counting poles of \eqref{ODE}, \eqref{PQ}.
Blow-up and blow-down orbits are the red and blue separatrices, respectively.
\end{thm}

\begin{proof}
By corollary \ref{corMorse}, the regularized flow \eqref{ODEreg} on $\mathbb{S}^2=\widehat{\mathbb{C}}$ of the nondegenerate ODE \eqref{ODE} is Morse.
By proposition \ref{propgrad}, the regularized flow is also gradient-like.
In particular section \ref{S2Morse} on Morse flows and their terminology applies, including proposition \ref{propPortraits}.
Orbit characterization by portraits follows from corollary \ref{corPortraits}.
This proves the theorem.
\end{proof}

\begin{thm}\label{thmPortrait2Rat}
Conversely, prescribe any mutually dual pair $\widetilde{\mathcal{C}}^\pm$ of nonempty, finite, connected multi-graphs on the $2$-sphere $\mathbb{S}^2$ with $d'\geq 0$ edges, each, and a total of $d=d^++d^-=d'+2$ vertices.\\
Then there exists a nondegenerate flow \eqref{ODE} with rational nonlinearity $f=P/Q$, such that the red blow-up and blue blow-down portraits $\mathcal{C}^\pm$ on $\mathbb{S}^2=\widehat{\mathbb{C}}$ are isomorphic to the given multi-graphs $\widetilde{\mathcal{C}}^\pm$, with orientation preserving homeomorphic embeddings.\\
In particular $d'=\mathrm{deg} \,Q, \ d=\mathrm{deg} \,P$, and the prefactor $a$ in \eqref{PQ} can be chosen as $a=1$.
The choice $a=-1$ reverses time and swaps $\mathcal{C}^\pm$.
\end{thm}

For a proof of this \emph{rational realization theorem}, see section \ref{PfthmRat}.

Corollary \ref{corPortraits} for Morse flows on $\mathbb{S}^2$ asserts orbit equivalence for homeomorphically embedded portraits $\widetilde{\mathcal{C}}^\pm\subset \mathbb{S}^2$.
We rephrase this result, for regularized rational ODEs on the Riemann sphere.

\begin{cor}\label{Rat=Morse}
The nondegenerate flows of regularized rational ODEs \eqref{ODEreg} on $\mathbb{S}^2=\widehat{\mathbb{C}}$, in totality, are $C^0$ orbit equivalent to the Morse flows on the $2$-sphere $\mathbb{S}^2$.
\end{cor}

\subsection{Polynomial ODEs}\label{ResPol}

The purely polynomial case $\dot{w}=P,\ Q=1$ of rational ODEs \eqref{ODE}, \eqref{PQ} has already been studied in the prequel \cite{FiedlerShilnikov}.
Complex rescalings $w\mapsto a'w+b'$ normalize \eqref{ODE} to 
\begin{equation}
\label{ODEP}
\dot{w}=P(w)=(w-e_1)\cdot\ldots\cdot(w-e_d)=w^d+\ldots+1;
\end{equation} 
compare \eqref{a=1}.
Again we consider $d\geq2$; see section \ref{Ex2} for $d=2$.
As a first result for polynomial nonlinearities $f=P$, we recall that the complex residues $\eta_j=1/P'(e_j)\neq 0$ of \eqref{omega,eta} can be prescribed arbitrarily, provided that
\begin{equation}
\label{sumseta}
\sum_{j=1}^{d}\,\eta_j\,=\,0\,,\quad\mathrm{but\ also}\quad
 \sum_{j\in J}\,\eta_j\,\neq\,0\,,
\end{equation}
for any nonempty subset $\emptyset\neq J \subsetneq\{1,\ldots,d\}$.
See theorem 1.1 in \cite{FiedlerShilnikov}.
For the necessary left constraint, see \eqref{sumeta0} above.
The nonvanishing condition may be partly technical.

However, $Q=1$ in \eqref{ODEP} violates the above degree relation \eqref{dd'}, unless $d=2$.
Passing to $z=1/w$, as in \eqref{dd'no}, drags the pole of $P$ from its hiding at $w=\infty$ to $z=0$:
\begin{equation}\label{ODEPtilde}
\dot{z}=-z^2 P(1/z)=\widetilde{P}(z)/\widetilde{Q}(z)\,.
\end{equation}
In the coordinate $z$, the denominator $\widetilde{Q}(z)=z^{d'}$ with $d':=d-2$ now complies with \eqref{dd'}.
We recover a rational ODE \eqref{ODE} with a new numerator polynomial $\widetilde{P}(z):=-z^dP(1/z)=-z^d+\ldots-1$, as defined in \eqref{PQtilde}, and with  zeros $\widetilde{P}(1/e_j)=0$. 
For polynomial degrees $d\geq3$, the single pole at $e_-=0$ in the denominator $\widetilde{Q}$ now reveals a hidden pole $w=\infty$ of multiplicity $d'=d-2\geq1$, in the original polynomial ODE \eqref{ODEP}; compare \eqref{dd'no}.

Smooth global regularization \eqref{ODEreg}, this time of \eqref{ODEPtilde}, provides
\begin{equation}\label{ODEPregz}
\dot{z}=\widetilde{f}_\mathrm{reg}(z):=(1+|z|^2)^{-(d-2)}\widetilde{P}(z)\overline{z}^{d-2}=-\overline{z}^{d-2}+\ldots \,.
\end{equation}
Applying the involutions \eqref{PQtilde} once again, $\tilde{\tilde{P}}=P$ and $\tilde{\tilde{Q}}=Q=1$ recover the proper regularization \eqref{ODEreg} in original variables:
\begin{equation}\label{ODEPregw}
\dot{w}=f_\mathrm{reg}(w):=(1+|w|^2)^{-(d-2)}P(w)\,.
\end{equation}
See figure \ref{fig1} for examples with $d=3,4$ and poles at $w=\infty$ of multiplicities $d'=1,2$.

We collect results from \cite{FiedlerShilnikov}, in our adapted formulation.
Polar coordinates $z= r \exp(-\mi\alpha)$ desingularize the pole $z=0$, and expand the regularization \eqref{ODEPregz} at the trivial zero of order $d-2$ as
\begin{align}
\label{r}
\dot{ r }\ &= r \big(-\cos((d-1)\alpha)+\ldots \big)\,,   \\
\label{alpha}
\dot{\alpha}\ &=\phantom{ r \big(-}\ \sin((d-1)\alpha)+\ldots\,,
\end{align}
again after rescaling time.
Such ``blow up'' at $z=0$, in the sense of singularity theory, replaces $w=\infty$ by the circle $ r =0,\ \alpha\in\mathbb{S}^1$.
In other words, local polar coordinates for $z=1/w$ essentially compactify the plane $w\in\mathbb{C}=\mathbb{R}^2$ to the closed unit disk $\mathbf{D}$, instead of the Riemann sphere $\widehat{\mathbb{C}}$.
Indeed the boundary circle $\mathbb{S}^1$ of $\mathbf{D}$ corresponds to the circle $\alpha\in\mathbb{S}^1$ at $ r =0$.

The disk $\mathbf{D}$ is called the \emph{Poincaré compactification} of $w\in\mathbb{C}$.
See figure \ref{fig1} (a),(c).
On the invariant boundary circle $ r =0$, a total of $2(d-1)$ hyperbolic saddle equilibria $\alpha_k=\pi k/(d-1)\in\mathbb{S}^1$ appear, for $k\ \textrm{mod}\,2(d-1)$.
In figure \ref{fig1} (a) and (c), we label the saddles at $\alpha_k$ by $\mathbf{k}$\,.
Restricted to the boundary circle, they are alternatingly unstable, at even $\mathbf{k}$, and stable, at odd $\mathbf{k}$.

Their corresponding stable and unstable separatrix counterparts, received from and sent into $|w|=|1/z|<\infty$, are marked red and blue, respectively.
In \eqref{ODEP}, they signify red blow-up and blue blow-down of $w(t)$, in finite real original time.
In the regularization \eqref{ODEPregz}, they also delimit the $2(d-1)$ local hyperbolic sectors of the equilibrium $z=0$, according to the planar classification of degenerate saddles by Poincaré. 
See section VII.9 in \cite{Hartman}, and (b), (d) of figure \ref{fig1}.
In (b), i.e. for $d=3$, red and blue simply mark the $2(d-1)=2+2$ separatrix half branches of the stable and unstable manifolds at the hyperbolic equilibrium $z=0$ of \eqref{ODEPregz}.
Colors match glossary \ref{Color}.

\begin{figure}[t]
\centering \includegraphics[width=0.9\textwidth]{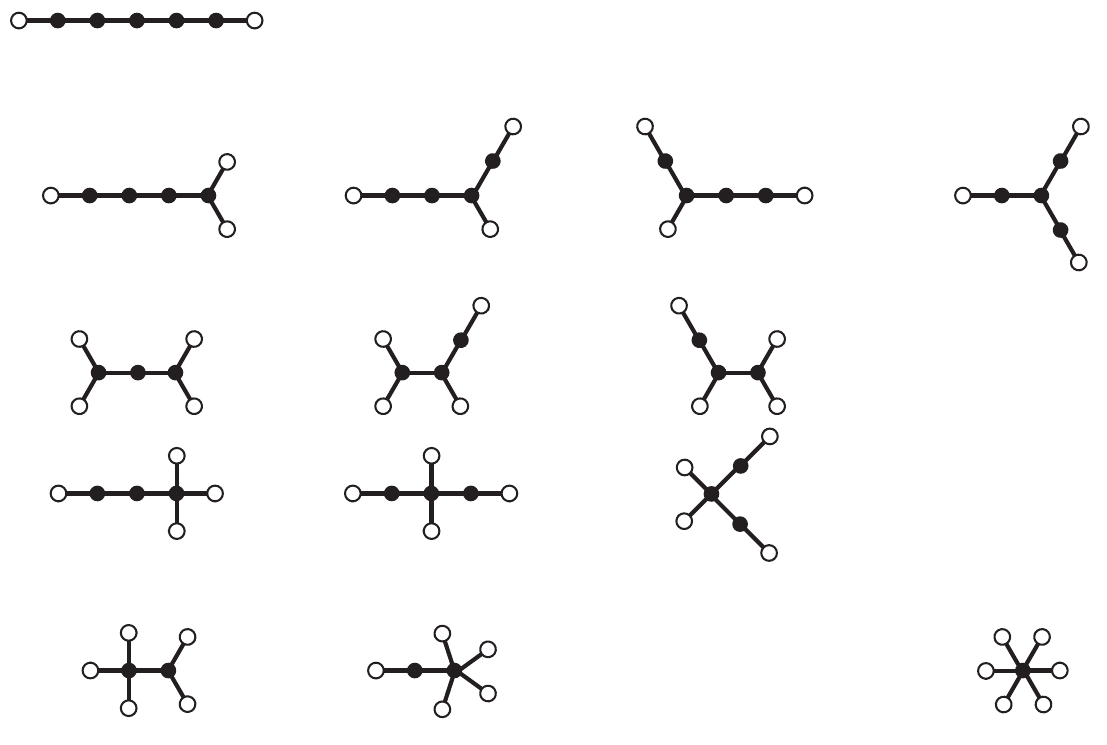}
\caption[The 14 planar trees with 7 vertices]{\emph{
The $14$ planar trees $\mathcal{T}$, alias reduced connection graphs $\mathcal{C}^*$, with $7$ vertices, up to orientation preserving equivalence.
Note the two pairs which are mirror-symmetric to each other, but not mirror-symmetric, individually.
Leaf vertices, of edge degree one, are distinguished by circles. 
Solid dots mark all other vertices.
The trees are sorted by increasing number of leaves.
Each tree possesses two bi-colorations which mark the vertices as red sources and blue sinks, alternatingly.
This also defines edge orientations, from sources to sinks.
Color swaps correspond to time reversal.
Since the number of vertices is odd, time reversible graphs $\mathcal{C}^*$, i.e., with isomorphic bi-colorations, cannot occur in this example.
}}
\label{fig2}
\end{figure}

In \cite{FiedlerShilnikov}, we have classified the regularized real-time flows, up to orientation preserving $C^0$ orbit equivalence; see \eqref{C0equiv}.
Abstractly, we compare orbit equivalence classes to equivalence classes of certain planar trees with $d$ vertices and $d-1$ edges.
A \emph{tree} is a finite, connected, undirected graph without cycles.
\emph{Planar trees} $\mathcal{T}$ are embedded, say, in the open unit disk $\mathbb{D}$.
We consider two planar trees $\mathcal{T}_1,\mathcal{T}_2$ as \emph{equivalent}, if there exists an orientation preserving homeomorphism $\frak{H}$ of $\mathbb{D}$ which acts as a graph isomorphism on the trees.
In other words, $\frak{H}{:}\ \mathcal{T}_1\rightarrow\mathcal{T}_2$ maps vertices to vertices, and edges to edges, preserving their adjacency relations and the left cyclic orderings of corresponding edges around each vertex.
Any direction of edges may be reversed under $\frak{H}$.
We do not require any distinguished vertices or edges to be mapped to each other, e.g. any vertices marked as ``roots'' or otherwise labeled. 
See figure \ref{fig2} for the $14$ planar trees with $7$ vertices.

Our specific tree $\mathcal{T}$ will neither coincide with the full connection graph $\mathcal{C}$, nor with its constituent portraits $\mathcal{C}^\pm$ discussed in section \ref{Morse}.
The full connection graph connects each source vertex to the same, single bottleneck, by its $d-1$ red stable separatrices: the pole at $w=\infty, z=1/w=0$.
The $d-1$ blue unstable separatrices connect the pole to the sink vertices.
Differently from the generic rational case of section \ref{ResRat}, however, the pole at $w=\infty$ is degenerate of multiplicity $d-2$ for $d\geq 4$.
In particular, transitivity fails: the heteroclinic orbits from sources to sinks are not the all-to-all bipartite graph.

Instead, the tree $\mathcal{T}$ will be the \emph{reduced connection graph} $\mathcal{C}^*$ of heteroclinic orbits from sources to sinks.
I.e., vertices of $\mathcal{C}^*$ correspond to sources and sinks in the Poincaré compactified flows \eqref{ODEPregw}.
Single edges now indicate the existence of connected \emph{families of heteroclinic orbits} between two vertices.
Since the union of heteroclinic orbits $e_1\leadsto e_2$\,, from any source to any sink, defines a quadrangle bounded by separatrices, the regularized ODE defines equivalent planar embeddings of the reduced connection graph $\mathcal{C}^*=\mathcal{T}$ in $\mathbb{D}$, independently of our specific choice of representative heteroclinic orbits.
Similarly, small perturbations of nondegenerate polynomial coefficients $P$ lead to equivalent embeddings.
Differently from the full connection graph $\mathcal{C}=\mathcal{C}^-\cup\mathcal{C}^+$ of rank -1 under the grading by Morse index, however, edges of the reduced connection graph $\mathcal{C}^*$ are of rank -2.

A third viewpoint are dual \emph{diagrams of} $d-1$ \emph{non-intersecting chords} of the closed unit disk $\mathbf{D}$, up to rotation.
Here the chords are closed undirected straight lines between their $2(d-1)$ endpoints on the disk boundary $\mathbb{S}^1$. 
Their endpoints are spaced equidistantly at angles $\beta_k:=\pi (k+\tfrac{1}{2})/(d-1),\ 0\leq k<2(d-1)$.
The non-intersecting chords are neither allowed to share endpoints, nor are they allowed to cross each other.
We call two chord diagrams \emph{equivalent}, if they coincide up to rotation of $\mathbf{D}$.
Chord diagrams are also called non-crossing (single-armed) \emph{handshakes} of $2(d-1)$ people $\beta_k$ at a round table, or partitions of $2(d-1)$ elements into non-crossing blocks of size 2.

\begin{thm}\label{Classthm}\emph{(\cite{FiedlerShilnikov}, Theorem 1.3)}
Consider flows \eqref{ODEP} of polynomials $P$ with $d\geq 2$ equilibria $e_j$\,.
For the reciprocal linearizations $\eta_j=1/f'(e_j)$ at $e_j$ we assume the \emph{nondegeneracy condition}
\begin{equation}
\label{sumsetaRe}
 \sum_{j\in J}\,\Re\,\eta_j\,\neq\,0\,,
\end{equation}
on any nonempty subset $\emptyset\neq J \subsetneq\{1,\ldots,d\}$.
Then the following two statements hold true, in terms of the orientation preserving equivalence classes just described.
\begin{enumerate}[(i)]
\item  The real-time phase portraits of the Poincaré compactifications correspond, one-to-one, to certain unlabeled, unrooted, undirected, planar trees $\mathcal{T}$ with $d$ vertices.
\item Equivalently, they correspond, one-to-one, to certain chord diagrams of $d-1$ unlabeled, non-intersecting chords of the unit circle.
\end{enumerate}
\end{thm}

Time reversal introduces a little glitch here:
the undirected trees are invariant under time reversal.
One remedy is to allow time reversing homeomorphisms in \eqref{C0equiv}.
A second remedy distinguishes trees $\mathcal{T}$ by red/blue bi-colorations of their vertices.
Theorems \ref{Classthm} and \ref{thmTree2Pol} hold in either case.
The combinatorial corollary \ref{corCountpol} will require the first option.

The nondegeneracy condition \eqref{sumsetaRe} excludes (purple) Lyapunov centers and periodic orbits.
It also excludes saddle-saddle loops via the interior $\mathbb{D}$ of the Poincaré compactification $\mathbf{D}$.
This establishes condition \emph{(iii)} of Morse definition \ref{defMorse} in the interior $\mathbb{D}$.
Only hyperbolicity of saddles is violated at the single pole $w=\infty$, in the regularization \eqref{ODEPregz}.
In absence of Lyapunov centers and interior saddle-saddle connections, the reduced connection graph $\mathcal{C}^*$ nevertheless  persists, under small perturbations of the polynomial coefficients $P$.
Although we saw transitivity of heteroclinicity fail, this recovers standard results on structural stability of Morse systems \cite{PalisSmale, Sotomayor, PalisdeMelo, Pilyugin}: small perturbations of equilibria only produce $C^0$ orbit equivalent Poincaré compactifications.

Correspondence theorem \ref{Classthm}, however, extends well beyond local perturbations.
See figure 59 in \cite{PalisdeMelo} for an example of planar phase portraits with planar trees $\mathcal{T}$ which are isomorphic, in the sense of abstract graph theory,
but fail to be equivalent, under any planar homeomorphism.
The planar phase portraits thus fail to be $C^0$ orbit equivalent.

Correspondence theorem \ref{Classthm} only asserts that each polynomial phase portrait corresponds to \emph{some} planar tree $\mathcal{T}$, alias circular handshake or non-intersecting chord diagram.
It does not assert that \emph{all} planar trees actually do arise.
In other words: the correspondence $\mathcal{C}^*\mapsto\mathcal{T}$ is injective. 
But how about surjectivity?
The following \emph{polynomial realization theorem} on surjectivity has been announced, but not proved, in theorem 1.4 of \cite{FiedlerShilnikov}.
We complete the proof in section \ref{PfPol} below.

\begin{thm}\label{thmTree2Pol}
Each unlabeled, unrooted, undirected, planar tree $\mathcal{T}$ with $d\geq2$ vertices, alias each circular handshake or each chord diagram of $d-1$ non-intersecting chords, is realized as the reduced connection graph $\mathcal{C}^*$ of \eqref{ODEP}, for suitable polynomials $P$ of degree $d$ and up to orientation preserving homeomorphisms.
\end{thm}

Explicit counts of planar trees have been provided by \cite{oeispol} and in theorem 2 of \cite{countpol1}.
See also $\omega_D^2$ in \cite{countpol2}, equation (10) and table 1. 
By correspondence theorem \ref{Classthm} and realization theorem \ref{thmTree2Pol}, we therefore obtain the following counts of real-time global phase portraits for polynomial ODEs \eqref{ODEP} with $d$ nondegenerate source/sink equilibria.
See figure \ref{fig2} again, for the $A_6=14$ cases with $d=7$ vertices.

\begin{cor}\label{corCountpol} \emph{\cite{countpol1, oeispol, countpol2}}
Up to orientation preserving equivalence, the number of unlabeled, unrooted, undirected, planar trees with $d$ vertices, alias chord diagrams of $d-1$ unlabeled, non-intersecting chords of the unit circle, is given by entry $A_{d-1}$ of sequence \emph{A002995} in the online encyclopedia of integer sequences \emph{\cite{oeispol}}.
The counts $A_{d-1}$ for $2\leq d\leq 16$ are
\begin{equation}
\label{countpol}
1, 1, 2, 3, 6, 14, 34, 95, 280, 854, 2694, 8714, 28640, 95640, 323396.
\end{equation}
An explicit expression for the general counts $A_{d-1}$ in closed form is
\begin{equation}
\label{countpol1}
\begin{aligned}
A_{d-1}\,=\,\tfrac{1}{2(d-1)d}\, \tbinom{2(d-1)}{d-1}\,&+\,\tfrac{1}{4(d-1)}\tbinom{d}{d/2}\,+\,\tfrac{1}{d-1}\,\phi(d-1) \,+\\
                  &+\,\tfrac{1}{2(d-1)}\,\sum_{k=2}^{d-2}\,\tbinom{2k}{k}\, \phi(\tfrac{d-1}{k}) \,.           
\end{aligned}
\end{equation}
Here $\phi$ denotes the Euler totient count of coprime elements, and the sum only runs over proper divisors $k$ of $d-1$.
For odd $d$, the second summand on the right is omitted.
\end{cor}

\subsection{Anti-polynomial ODEs}\label{ResAntipol}

In this section we discuss \emph{anti-holomorphic gradient} ODEs 
\begin{equation}
\label{ODEF}
\dot{\mathbf{w}}=\overline{\mathbf{f}(\mathbf{w})}\,,\quad\mathrm{for}\ \mathbf{f}=-\mathrm{grad}\,F\,.
\end{equation}
Here the potential $F{:}\ \mathbb{C}^N\rightarrow\mathbb{C}$ of $\mathbf{f}$ is assumed holomorphic.
Decomposing $\mathbf{w}=\mathbf{u}+\mi \mathbf{v}\in\mathbb{C}^N$ and $F=G+\mi H\in\mathbb{C}$ into real and imaginary parts, we can rewrite \eqref{ODEF} in components $\mathbf{w}=(\mathbf{u},\mathbf{v}),\ \mathbf{u}=(u_1,\ldots,u_N),\ \mathbf{v}=(v_1,\ldots,v_N)$ as
\begin{equation}
\label{ODEwn}
\begin{aligned}
    \dot u_n&=-G_{u_n}  = -H_{v_n}\,;\\
    \dot v_n&=-G_{v_n}  = +H_{u_n}\,.
\end{aligned}
\end{equation} 
Suppressing the arguments $(\mathbf{u},\mathbf{v})$ of $G,H$, we have used the Cauchy-Riemann equations for $F$ here. 
In particular this proves the following.

\begin{prop}\label{propGradHam}
In real time, the anti-holomorphic gradient ODE \eqref{ODEF} is gradient \emph{and} Hamiltonian, simultaneously.
More precisely, the middle parts of \eqref{ODEwn} show that \eqref{ODEF} is the negative gradient flow associated to $G=\Re F$, under the standard Euclidean scalar product of $(\mathbf{u},\mathbf{v})\in\mathbb{R}^{2N}$.
The right hand sides show that \eqref{ODEF} is Hamiltonian with respect to $H=\Im F$, under the standard symplectic form $\sum du_n\wedge dv_n$ on $\mathbb{R}^{2N}$.
\end{prop}

A scalar example, $N=1$, are rational ODEs \eqref{ODE} with trivial numerators $P=1$ and polynomial denominator $Q$.
After multiplication by an Euler multiplier $|Q|^2$, indeed, this amounts to anti-holomorphic vector fields
\begin{equation}
\label{Qbar}
\dot{w}=\overline{Q(w)}.
\end{equation}
For the primitive polynomial $F'=-Q$, say with $F(0)=0$, proposition \ref{propGradHam} and \eqref{ODEwn} imply
\begin{equation}
\label{ODEwbar}
    \dot{w}=-G'(w)=\,\mi\, H'(w)\,.
\end{equation}
The imaginary unit $\mi$ in \eqref{ODEwbar} defines the standard symplectic form $du\wedge dv$ on $\mathbb{C}=\mathbb{R}^2$.
The commuting perpendicular flow, in imaginary time, rotates the roles of $G$ and $H$.

To be specific, consider $Q(w)=(w-e_1')\cdot\ldots \cdot (w-e_{d'})=w^{d'}+\ldots$ of degree $d'\geq1$, with simple zeros $e_k'\in E'$\,.
In \eqref{Qbar}, these become saddle equilibria.
For generic nondegeneracy, we only assume distinct values of the Hamiltonian $H$ on $E'$.
By \cite{critval}, these values can be prescribed arbitrarily.
For the original ODE \eqref{ODE}, the substitution $z=1/w$ yields
\begin{equation}
\label{zdbar}
\dot{z}=-z^{d'+2}\,(1+\ldots)\,.
\end{equation}
Here we have omitted anti-holomorphic terms of higher orders $\overline{z},\ldots,\overline{z}^{d'}$\,, and an Euler multiplier $|z|^{2d'}$.
Expectedly, this regularization reveals an equilibrium of order $d=d'+2$ at $w=\infty$.
Analogously to  section \ref{ResPol}, this allows a Poincaré compactification of the anti-polynomial ODE \eqref{Qbar} in polar coordinates $z=r\exp(-\mi\alpha)$ as 
\begin{align}
\label{rzbar}
\dot{ r }\ &= -r \big(\cos((d'+1)\alpha)+\ldots \big)\,;   \\
\label{alphazbar}
\dot{\alpha}\ &= \ - \sin((d'+1)\alpha)+\ldots\,;
\end{align}
compare \eqref{r}, \eqref{alpha}.
We now obtain $d-1=d'+1$ sinks at $r=0,\ \alpha_j=j\pi/(d'+1)$, for even $j$.
The sinks alternate with $d'+1$ sources, at odd $j$.

\begin{figure}[t]
\centering \includegraphics[width=\textwidth]{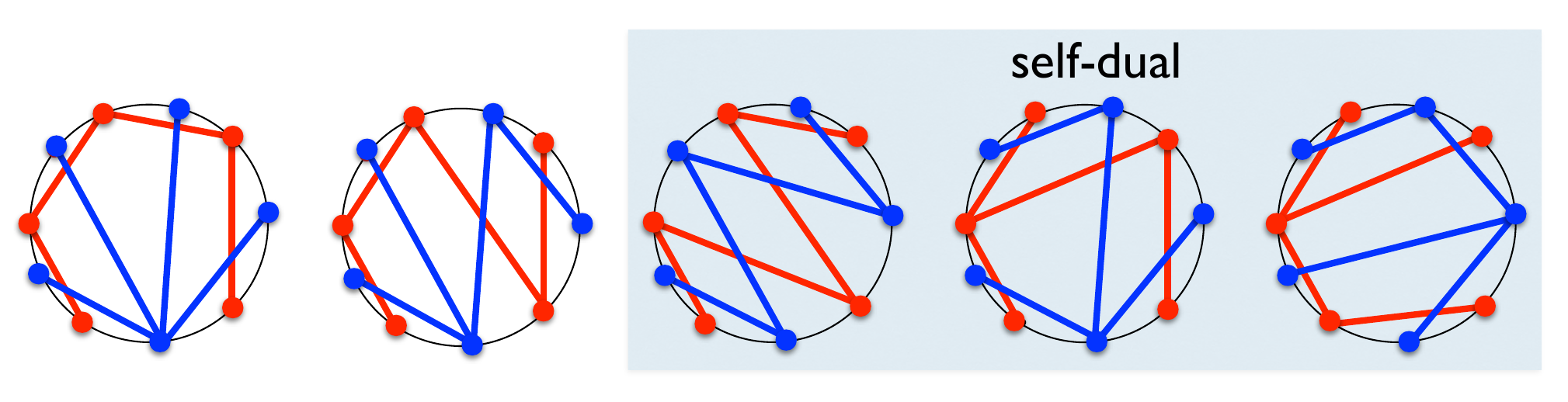}
\caption[The 7 noncrossing trees with 5 vertices]{\emph{
Poincaré compactifications of the anti-polynomial ODE \eqref{Qbar}, for nondegenerate polynomials $Q$ of degree $d'=4$.
Counting up to orientation preserving or orientation reversing $C^0$ equivalence, we obtain $A_{d'}=7$ non-crossing trees $\mathcal{C}^\pm$ with $d'+1=5$ vertices on the circle, each, and with $d'=4$ chords, as edges.
The blue trees $\mathcal{C}^-$ are dual to the red trees $\mathcal{C}^+$.
Edges are stable and unstable separatrix pairs between their saddles and red/blue sources/sinks at the boundary circle.
In the three self-dual cases, on the right, the dual portraits $\mathcal{C}^+$ and $\mathcal{C}^-$ are equivalent, and we obtain three connection graphs $\mathcal{C}=\mathcal{C}^-\cup\mathcal{C}^-$, each reversible in real time.
In each of the two cases on the left, time-reversal leads to a new, non-homeomorphic connection graph of swapped colors, contributing two more connection graphs.
Together, this accounts for the $A_{d'}=7$ connection graphs $\mathcal{C}$, alias orbit equivalence classes of \eqref{Qbar}.
}}
\label{fig3}
\end{figure}

The $d'$ hyperbolic saddles $e_k'\in E'$ are finite, each equipped with red blow-up separatrices and blue blow-down separatrices, as usual.
We describe the Poincaré compactified flow in the closed unit disk $\mathbf{D}$.
Compare glossary \ref{Color} for terminology, and see figure \ref{fig3}.

The flow is a Morse flow, in the sense of definition \ref{defMorse}.
Indeed, all orbits are heteroclinic because the flow \eqref{ODEwbar} is gradient with respect to the real potential $G$.
Each portrait $\mathcal{C}^\pm$ then consist of $d'+1$ source/sink vertices on the boundary circle $\mathbb{S}^1$, respectively, and $d'$ edges in the interior $\mathbb{D}$.
The Hamiltonian structure \eqref{ODEwbar} preserves the real Hamiltonian $H=\Im F$, along real-time solutions, and hence prevents saddle-saddle connections between the saddles $e_k'$\,, which reside at different levels of $H$.
Note duality of $\mathcal{C}^\pm$.

Each portrait $\mathcal{C}^\pm$ is acyclic: absence of interior source and sink vertices prevents faces.
In particular each portrait is a tree $\mathcal{C}^\pm=\mathcal{T}^\pm$, each connected by Euler's formula.
Tangencies of the separatrix edges at the boundary sinks/sources are along their slow stable/unstable manifolds, and hence perpendicular to the boundary circle.

Up to orientation preserving homeomorphisms $\frak{H}$, in other words, each portrait $\mathcal{C}^\pm$ is a (connected) tree of \emph{noncrossing} chords (straight line edges) between $d'+1$ vertices on a circle.
Noncrossing means absence of intersections, except at endpoints.
Following \cite{Noy}, we call such (unlabelled) trees \emph{nc-trees}.
See figure \ref{fig3} for all seven nc-trees with $d'=4$ edges.
In analogy to polynomial theorems \ref{Classthm} and \ref{thmTree2Pol}, we obtain the following result.

\begin{thm}\label{nctree=antipol}
Consider regularized anti-polynomial flows \eqref{Qbar} for univariate polynomials $Q$ of degree $d'\geq1$.
Let $F'=Q$ denote the primitive polynomial $F$ of $Q$.
We assume simplicity of the $d'$ zeros $e_k'\in E'$\, and pairwise distinct values of the Hamiltonian $H=\Im F$ on them.
Then the anti-polynomial flow of \eqref{Qbar} is nondegenerate.
Moreover, the following two statements hold true, up to orientation preserving homeomorphisms.
\begin{enumerate}[(i)]
\item  The real-time portraits $\mathcal{C}^\pm$ of the Poincaré compactifications correspond, one-to-one, to certain dual pairs $\mathcal{T}^\pm$ of nc-trees with $d'+1$ vertices on the circle.
\item Conversely, any dual pair $\mathcal{T}^\pm$ of nc-trees with $d'+1$ vertices arises as portraits $\mathcal{C}^\pm$ from an anti-polynomial flow \eqref{Qbar}, for some nondegenerate univariate $Q$ of degree $d'$, as above.
\end{enumerate}\end{thm}

We have proved claim \emph{(i)} above.
We prove realization claim \emph{(ii)} in section \ref{PfAntipol}.

In analogy to counting corollary \ref{corCountpol}, and due to Noy, these claims allow us to determine the number of anti-polynomial flows, via the following (corrected and adapted) count of nc-trees.
See figure \ref{fig3} again.

\begin{cor}\label{corCountanti} \emph{ \cite{Noy, oeisanti}}
Up to orientation preserving or reversing equivalence, the number of unlabeled nc-trees with $d'+1$ vertices on the circle, and hence $d'$ chords as edges, is given by entry $A_{d'}$ of sequence \emph{A296533} in the online encyclopedia of integer sequences \emph{\cite{oeisanti}}.
The counts for $1\leq d'\leq 13$ are
\begin{equation}
\label{countanti}
A_{d'}' = 1, 1, 3, 7, 28, 108, 507, 2431, 12441, 65169, 351156, 1926372, 10746856.
\end{equation}
Explicit expressions \emph{\cite{Noy}} for the general counts $A_{d'}$ in closed form are
\begin{equation}
\label{countanti1}
\begin{aligned}
   A_{d'}'\,=\,\frac{1}{(2d'+1)(2d'+2)}\binom{3d'}{d'}&+\frac{3}{3d'+1}\binom{(3d'+1)/2}{(d'-1)/2)}\,,\quad\ \mathrm{for}\ d'\ \mathrm{odd}\,;   \\[2mm]
    A_{d'}'\,=\,\frac{1}{(2d'+1)(2d'+2)}\binom{3d'}{d'}&+\quad\frac{1}{2d'+2}\binom{3d'/2}{d'/2}\,,\qquad\quad\mathrm{for}\ d'\ \mathrm{even}\,.  
\end{aligned}
\end{equation}
\end{cor}

The Noy count requires us to admit, both, orientation preserving and orientation reversing homeomorphisms $\frak{H}$ in our definition of $C^0$ orbit equivalence.
This is not in conflict with the finer classification up to orientation preserving homeomorphisms $\frak{H}$, in theorem \ref{nctree=antipol}.
Indeed, the finer classification implies that theorem \ref{nctree=antipol} remains valid, verbatim, for the coarser equivalence under \emph{all} homeomorphisms.
The Noy count requires the latter.
Time reversal is not admitted.
Indeed, the Noy count does not collate dual pairs $\mathcal{C}^\pm$, but counts non-homeomorphic nc-trees $\mathcal{C}^+$ and $\mathcal{C}^-$ separately.
Analogously, non-homeomorphic time-reversed copies of directed connection graphs $\mathcal{C}=\mathcal{C}^+\cup\mathcal{C}^-$, graded by Morse index, are counted separately.
Any of the two nc-trees $\mathcal{C}^\pm$ determines the other,  by time reversal or, equivalently, by duality.
We therefore obtain the Noy counts for graded connection graphs $\mathcal{C}$, up to orientation preserving or reversing equivalence $\frak{H}$, but excluding time reversal.

\section{Three examples}\label{Ex}

\subsection{Quadratic vector fields: $d=2,\ d'=0$}\label{Ex2}
We recall the polynomial ODE case $Q=1$ of quadratic Riccati nonlinearities $\dot{w}=f(w)=P(w)$, which the perturbative PDE results by Masuda were based on.
Up to Möbius transformations \eqref{Mobius}, there are only two cases in complex time: the purely quadratic Masuda case $\dot{w}=w^2$ and the ``general'' quadratic case 
\begin{equation}
\label{ODE2}
\dot{w}=w^2-1\,.
\end{equation}
In the purely quadratic case, all nonstationary orbits in $\widehat{\mathbb{C}}$ are homoclinic loops to the algebraically double equilibrium $w=0$, e.g. by explicit separation of variables.
In either case, the flow maps $w\mapsto\Phi^t(w)$ are globally biholomorphic, i.e. Möbius transformations \eqref{Mobius} of $\widehat{\mathbb{C}}$.
Along real one-dimensional time rays $r\mapsto t=r \exp(\mi \theta)$, the flows $\Phi^t$ form one-parameter subgroups of the Möbius group $\mathrm{PSL}(2,\mathbb{C})$ which fix the equilibria.

\begin{figure}[t]
\centering \includegraphics[width=0.86\textwidth]{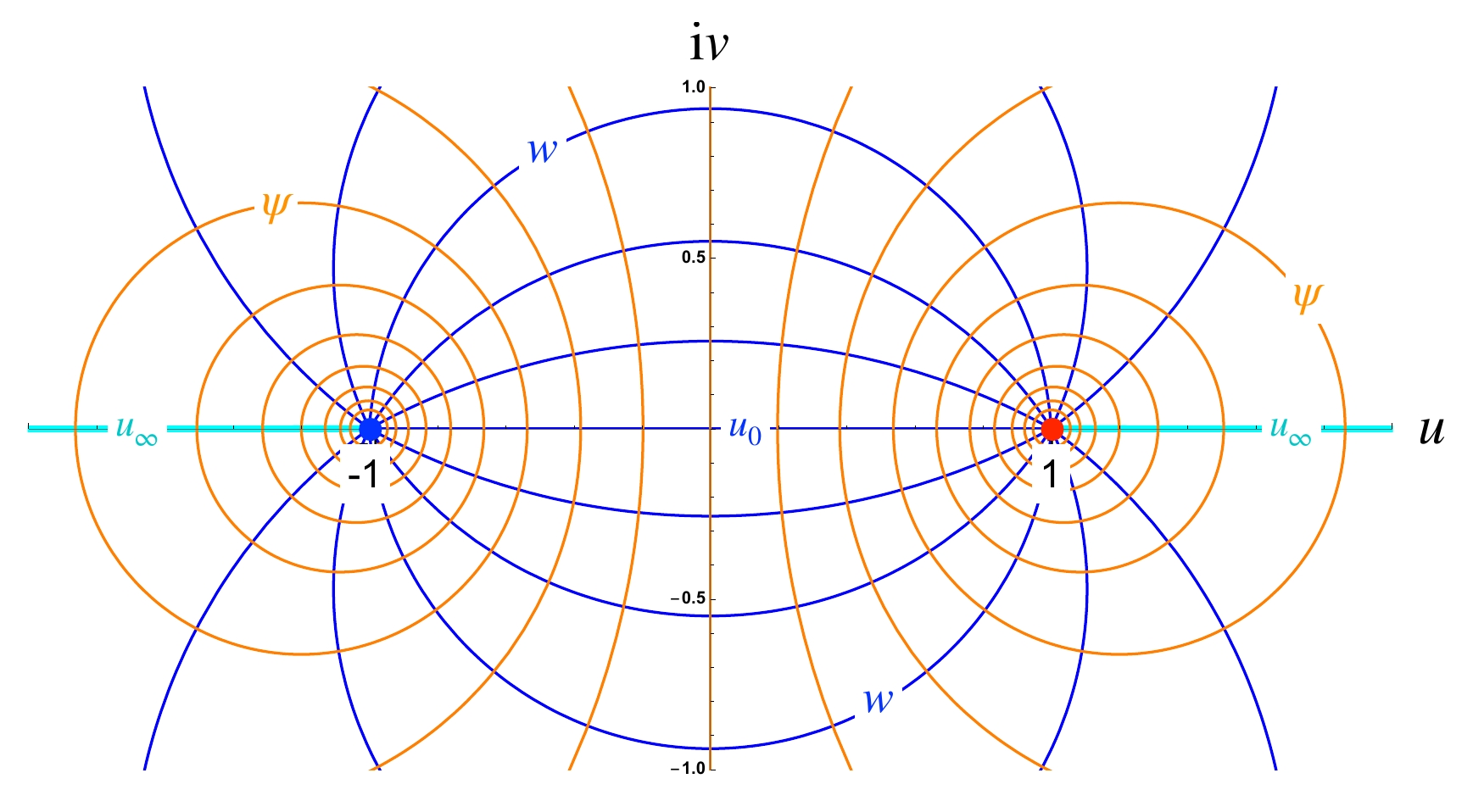}
\caption[The quadratic complex ODE]{\emph{
Phase portrait of the quadratic complex ODE \eqref{ODE2} with source $e_1=+1$ (red) and sink $e_2=-1$(blue).
All other real-time orbits are heteroclinic $+1\leadsto -1$ (blue).
They foliate the cylinder $\widehat{\mathbb{C}}\setminus\{\pm 1\}$ by circular segments with midpoints on the imaginary axis.
The fake ``blow-up'' and ``blow-down'' orbits $u_\infty$ through $w=\infty$ (cyan) may be concatenated to become heteroclinic, too, and may be circumnavigated by a Masuda detour.
We call them fake, because their definition and selection is not invariant under Möbius transformations.
The blow-up and blow-down portraits $\mathcal{C}^\pm=\{\pm1\}$ are single-vertex graphs without edges.
\\
Nonstationary orbits in imaginary time (orange), in contrast, provide an iso-periodic foliation of minimal period $\pm\pi$, nested around $\pm 1$.
Since the flow $t\mapsto\Phi^t(w_0)$ is conformal, for any fixed $w_0\in\widehat{\mathbb{C}}\setminus\{e_1,e_2\}$, the blue and orange circle families are mutually orthogonal.
Again, the iso-periodic orbit on the imaginary axis experiences fake  ``blow-up'' and ``blow-down'' in finite imaginary time.
\\
Still, as in the PDE context of \cite{FiedlerFila}, such blow-up and blow-down in imaginary time remains related to real heteroclinicity $u_0$ in real time.
}}
\label{fig4}
\end{figure}

In the general case \eqref{ODE2}, we illustrate the flows in real time, $\theta=0$, and in imaginary time $\theta=\pi/2$. See figure \ref{fig4}, and also the discussion in \cite{FiedlerShilnikov} for further details.
In real time, we obtain a unique source $e_1=+1$ (red) and sink $e_2=-1$ (blue), but no poles. 
All other real-time orbits (here blue) are heteroclinic $+1\leadsto -1$.
On the real axis, we encounter the fake ``blow-up'' and ``blow-down'' orbits through $w=\infty$ (cyan) which can be concatenated to become heteroclinic, too.
``Blow-up'' can therefore be circumnavigated by a Masuda detour.

In imaginary time, the two equilibria become Lyapunov centers.
All nonstationary orbits (orange) are iso-periodic of minimal period $\pm\pi$.
They foliate the cylinder $\widehat{\mathbb{C}}\setminus\{e_1\,,\ e_2\}$.
Because the flow $t\mapsto\Phi^t(w_0)$ is conformal, for any fixed $w_0\in\widehat{\mathbb{C}}\setminus\{e_1,e_2\}$, the heteroclinic and iso-periodic foliations are mutually orthogonal, just as the real and imaginary time directions are.

The Möbius involution $z=(w+1)/(w-1)$ globally linearizes \eqref{ODE2}, mapping  
the equilibria $w=\pm 1$ to $z=\infty$ and $z=0$, respectively.
Indeed $\dot{z}=-2z$ becomes linear.
In $z$, the heteroclinic orbits of \eqref{ODE2} become inward radial and, perpendicularly, the iso-periodic orbits become circles of constant radius $|z|$.
Since Möbius transformations preserve families of circles (including lines) in $\mathbb{C}$, this observation also shows the circular nature of all orbits in figure \ref{fig4}, for both real and imaginary time.

\subsection{A cubic rational example: $d=3,\ d'=1$}\label{Ex3}

\begin{figure}[t]
\centering \includegraphics[width=0.5\textwidth]{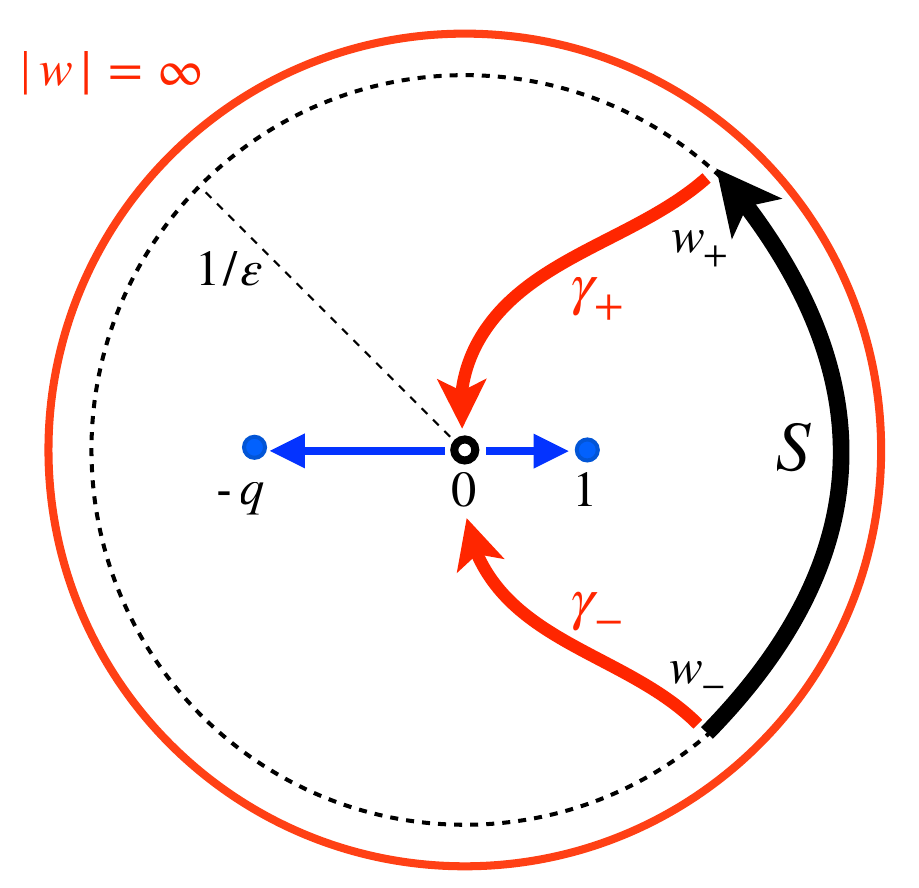}
\caption[A cubic rational example]{\emph{
Phase portrait of the real case $\dot{w}=-(w-1)(w+q)/w$ of ODE \eqref{cubic}, i.e. for specific parameters $a=-1,\ \theta=0,\ q>0$.
Note the three equilibria at $w=1,\,-q$ (sinks, blue), and $w=\infty$ (source, red).
Since real $f$ commute with complex conjugation, the flow is symmetric to the invariant real axis. 
In particular, the unstable blow-down separatrices (blue) of the unique pole at $w=0$ (black circle) are real heteroclinic to the sinks.
In backward time $t\searrow-\infty$, the two stable separatrices (red) of $w=0$ form a symmetric heteroclinic pair $\infty\leadsto 0$, sharing the unique source $w=\infty$.
They cross any large circle of radius $1/\eps$ (dashed black) transversely at $w_\pm=\eps^{-1}\exp(\mi\psi_\pm)$.
This defines a circle segment $S$ (solid black) between symmetric opening angles $\psi_-=-\psi_+=-\pi/(1+q)$; see \eqref{limpsi0}.
}}
\label{fig5}
\end{figure}

To prepare the proof of rational realization theorem \ref{thmPortrait2Rat}, we study the \emph{cubic example}
\begin{equation}
\label{cubic}
\dot{w}=h(w):=a\,(w-e_1)(w+q e_1)/w, \quad\mathrm{for}\quad \Re a<0<q\ \mathrm{and}\ e_1:=\exp(\mi\theta)\,.
\end{equation}
Reversing the sign of $a\in\mathbb{C}$ reverses the flow and swaps colors.
Note the single pole at $w=0$.
For the real case of $a=-1$ and rotation angle $\theta=0$, we obtain the regularized flow illustrated in figure \ref{fig5}.
The only finite equilibria $w=e_1$ and $w=e_2:=-qe_1$ are sinks, with linearizations
\begin{equation}
\label{cubiclin}
\begin{aligned}
h'(e_1)&=a(1+q)\,,\\
h'(e_2)&= a(1+1/q)\,.
\end{aligned}
\end{equation}
We call the example cubic, because substitution of $z=1/w$ reveals a source at $z=0$,
\begin{equation}
\label{cubicz}
\dot{z}=-a\,(1-e_1z)(1+q e_1z)z\,,
\end{equation}
with linearization $-a$; see also \eqref{sumeta0}. 
A Möbius transformation, which relocates the source away from $w=\infty$, exhibits the rational cubic example as nondegenerate; see definition \ref{defnondeg}.
By corollaries \ref{corMorse} and \ref{corPortraits}, the associated regularized flow is gradient-like Morse. 
In particular, the phase portraits $\mathcal{C}^\pm$ are mutually dual, by time reversal. 

The two red stable separatrices $\gamma_\pm$ of $w=0$ have to converge to the unique source $w=\infty$, in backward time.
The unique edge of the red blow-up portrait $\mathcal{C}^+$ is that heteroclinic pair $\gamma_\pm{:}\  \infty\leadsto 0$, which forms a loop to the unique source vertex $w=\infty$.
The two blue unstable separatrices, which emanate from $w=0$ in forward time, define the single blow-down edge of $\mathcal{C^-}$.
By duality of $\mathcal{C}^\pm$, they converge to the two sinks $e_1$ and $e_2$\,, respectively.
Their convergence is spiraling, for $\Im a\neq 0$.
Topologically, $\mathcal{C^+}$ is a circle which crosses the dual interval $\mathcal{C^-}=[-q,1]$.

We aim for a more analytic description.
For general rotation angles $\theta$, choose $\eps>0$ small enough such that the vector field $f$ points inward at large circles $|w|=1/\eps$.
For $\eps\searrow 0$, indeed, we obtain scalar products
\begin{equation}
\label{cubictrv}
\langle w,h(w)\rangle \,=\,\Re( \, \overline{w}\,h(w))\,=\,\eps^{-2}\,\Re a+\mathcal{O}(\eps^{-1})\,<\,0\,.
\end{equation}
By transversality, the two red separatrices $\gamma_\pm$ intersect the circle at unique points $w_\pm=\eps^{-1}\exp(\mi\psi_\pm)$\,.
Let $S$ denote the left winding circle segment from $\psi_-$ to $\psi_+$\,.

\begin{lem}\label{lempsi}
For $q>0$ and small $\eps> 0$, the opening angles $\psi_\pm=\psi_\pm(\eps,q,\theta))$ satisfy
\begin{align}
\label{psitheta}
\psi_\pm(\eps,q,\theta)&=\psi_\pm(\eps,q,0)+\theta\,;\\
\label{limpsi0}
\psi_+(\eps,q,\theta)-\psi_-(\eps,q,\theta)&=2\pi/(1+q)\,.  
\end{align}
By \eqref{psitheta}, the right hand side of \eqref{limpsi0} is independent of $\theta$ and, notably, of $\eps$.\\
In other words, prescribe any crossing angles $\psi_+\neq \psi_- \mod 2\pi$.
Then $\psi_\pm=\psi_\pm(\eps,q,\theta))$ can be realized by the cubic example \eqref{cubic}, for suitable choices of $q>0,\ \theta\in\mathbb{R}$, and any small enough $\eps>0$.
\end{lem}

\begin{proof} 
We first prove claim \eqref{psitheta}, by equivariance of \eqref{cubic} under rotation by $\theta\in\mathbb{R}$.
Let $w^\theta$ denote any solution of \eqref{cubic}, for any fixed angle $\theta$.
Then $w^0:=\exp(-\mi\theta)w^\theta$ solves the same equation, but for $\theta=0$\,.
This proves claim \eqref{psitheta}.

To prove the asymptotic claim \eqref{limpsi0}, we integrate $\omega=dw/h$ along a closed contour $\Gamma$.
The contour consists of three parts: the circle segment $S$ from $w_-$ to $w_+$\,, the red separatrix part $\gamma_+$ inside the circle of radius $1/\eps$ and, in reverse time direction, the remaining separatrix part $\gamma_-$\,.
The separatrix parts $\gamma_\pm$ are orbits of \eqref{cubic} in real time.
Their contribution to the imaginary part of the integral of $\omega$ along $\Gamma$ therefore vanishes, by separation of variables as in \eqref{sov}.
The residue theorem then implies
\begin{equation}
\label{Sres}
\Im \int_S \omega = \Im \int_\Gamma \omega =\Im\Big(  2\pi\mi /h'(e_1) \Big) =  \Re(1/a)\, 2\pi/(1+q)\,.  
\end{equation}
Explicit evaluation of the same integral for $|w|=1/\eps$ yields
\begin{equation}
\label{Sint}
\Im \int_S \omega = \Im \int_{\psi_-}^{\psi_+} \mi w\, d\psi /h(w)=\Re(1/a)\,(\psi_+-\psi_-)
\end{equation}
for all small $\eps>0$, by Cauchy's theorem in annuli of $\eps$.
Combining \eqref{Sres} and \eqref{Sint} proves claim \eqref{limpsi0}, and the lemma. 
\end{proof}

\subsection{A polynomial induction}\label{Exd}

\begin{figure}[t]
\centering \includegraphics[width=0.9\textwidth]{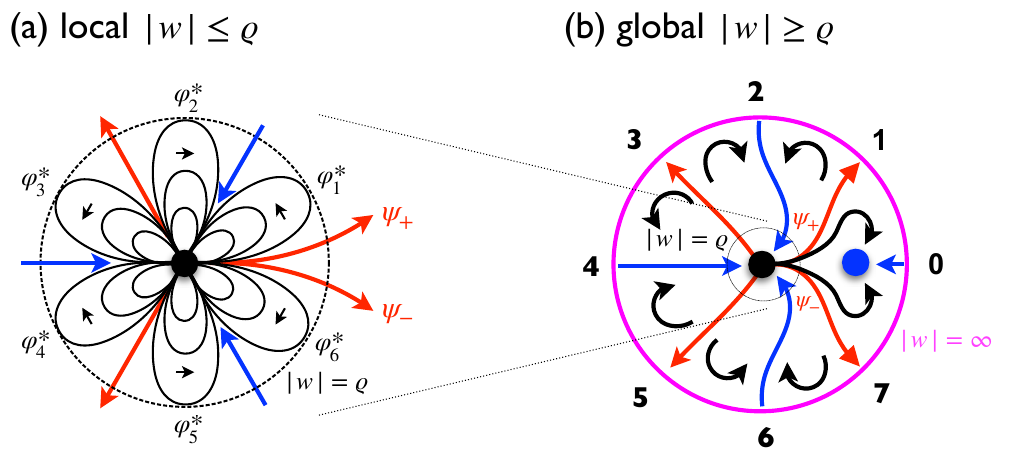}
\caption[A  degenerate polynomial example]{\emph{
Phase portrait of ODE \eqref{wdtheta} for multiplicity $d=4$ and $\theta=0$.\\
(a) Local view. Note the $2(d-1)$ elliptic sectors of homoclinic loops.
Elliptic sectors alternate with $2(d-1)$ separating sectors, which contain the separatrices.\\
(b) Global view. The real blue blow-down separatrix $\mathbf{0}$ terminates at the nontrivial sink $w=\exp(-\mi\theta)=1$.
All other separatrices $\mathbf{k}\neq0$ are heteroclinic between $w=\infty$ and $w=0$.
Note the two red blow-up separatrices $\mathbf{1}$ and $\mathbf{2d-1}$, which straddle the real axis.
All other orbits between the straddling separatrices and $\mathbf{0}$ are heteroclinic, from $w=0$ to $w=1$.
All remaining non-separatrices are loops to $w=0$.
}}
\label{fig6}
\end{figure}

To prepare the proof of polynomial realization theorem \ref{thmTree2Pol}, we study the polynomial examples
\begin{equation}
\label{wdtheta}
\dot{w}=P(w):=w^d(1-\exp(\mi\theta)w)
\end{equation}
of degrees $d+1\geq 3$, for certain fixed rotations $\theta=\theta_\ell$\,.
For $\theta=0$ see figure \ref{fig6}.
We provide a summary at the end of this section.

In the terminology of \cite{Hartman}, chapter VII, the trivial equilibrium $w=0$ of multiplicity $d$ possesses $2(d-1)$ elliptic and no elliptic sectors.
The elliptic sectors are separated, in particular, by $2(d-1)$ separatrices which are stable and unstable, alternatingly.
The regularized pole $z=1/w=0$ at $w=\infty$ of multiplicity $d-1$, in contrast, possesses $2d$ hyperbolic and no elliptic sectors.
Based on separatrices, we adapt this concept to our specific purposes.

We first investigate equilibria $w=0$ of multiplicity $d\geq2$, locally.
Consider
\begin{equation}
\label{wd}
\dot{w}=w^d\,.
\end{equation}
In polar coordinates $w= r \exp(\mi\psi)$, this reads
\begin{align}
\label{rw}
\dot{ r }\ &= r \cos\big((d-1)\psi\big)\,;   \\
\label{alphaw}
\dot{\psi}\ &=\,\phantom{ r }\sin\big((d-1)\psi\big)\,;
\end{align}
compare \eqref{rzbar}, \eqref{alphazbar}.
Equilibria are located on the ``circle'' $ r =0$, at angles
\begin{equation}
\label{psij}
\psi=\psi_j:=j \pi/(d-1)\,,
\end{equation}
for $0\leq j< 2(d-1)$. 
For even $j$, we obtain source equilibria with strong unstable manifolds along the circle $ r =0$.
The equilibria at odd $j$ are sinks with strong stable manifolds along the circle $ r =0$.
Convergence outside the strong manifolds occurs tangent to the perpendicular eigendirections of constant $\psi$.
Globally, \eqref{rw}, \eqref{alphaw} is a skew product over $\psi$. 
In particular, the $d-1$ red blow-up orbits, and the $d-1$ blue blow-down orbits are strictly radial, at the invariant and alternating source/sink angles $\psi=\psi_j$\,.

Let $ r =|w|\leq\rho>0$ denote any closed disk centered at $w=0$.
\emph{Elliptic sectors}, in the sense of \cite{Hartman}, consist of loops to $w=0$ contained in the disk.
The sectors are bounded by orbits through tangents of the vector field to the boundary circle $ r =\rho$.
See figure \ref{fig6}(a).
By \eqref{rw}, the circle tangencies occur at the intermediate angles
\begin{equation}
\label{phij*}
\psi=\phi_j^*:=(j-\tfrac{1}{2})\pi/(d-1)=\tfrac{1}{2}(\psi_{j-1}+\psi_j)\,,
\end{equation}
for $1\leq  j\leq 2(d-1)$. 
Orbits through the $2(d-1)$ tangencies remain in the $\rho$-disk, converging to the neighboring sink/source at $r=0$, heteroclinically.
Their interiors, likewise, are foliated by heteroclinic orbits of the same type.
For $w$, of course, these are nested homoclinic loops to $w=0$.
The $2(d-1)$ alternating entry and exit intervals of \eqref{wd} on the circle $r=\rho$, between the interior tangencies at $\phi_j^*$\,, define the \emph{separating sectors}.
They are marked by blue and red arrows in figure \ref{fig6}(a), respectively.
All trajectories in entry intervals converge to $w=0$, tangentially to the sink angle $\psi=\psi_j$ of \eqref{psij} in that interval, i.e. for the associated odd $j$.
For exit intervals, i.e.~for even $j$ and in reversed time, we obtain the corresponding source angles $\psi_j$\,.

Due to the homoclinic loops, by the way, the flow \eqref{wd} is not gradient-like.
In fact, the equilibrium $w=0$ is \emph{not} an isolated invariant set, in the sense of Conley index theory; see \cite{Conley, Mischaikow}.
Indeed all tangencies occur from inside the disk $|w|\leq\rho$.
The exterior disk $ r =|w|\geq\rho$ of the same circle on the Riemann sphere $w\in\widehat{\mathbb{C}}$, on the other hand, provides an isolating neighborhood of the pole $w=\infty$, alias of the equilibrium $z=1/w=0$, of multiplicity $d-2$.
Indeed, the same interior tangencies now occur from outside the exterior disk $|w|\geq\rho$.
The $2(d-1) $ separatrices, which emanate from the pole $w=\infty$ by Poincaré compactification, enter/leave the $\rho$-disk between the tangencies $\phi_j^*$ in corresponding order -- one per interval.
In the language of \cite{Hartman}, they define the $2(d-1)$ \emph{hyperbolic sectors} of $z=1/w=0$.

We now return to global aspects of the full example \eqref{wdtheta}. 
Let $ r =|w|\leq\rho>0$ denote a \emph{small} disk centered at $w=0$.
Inside the small disk, \eqref{wdtheta} is a small perturbation, of order $d+1$, of the case $\dot{w}=w^d$ just discussed; see \eqref{wd}.
By the implicit function theorem, the elliptic sectors defined by the tangencies at angles $\phi_j^*$ persist at angles
\begin{equation}
\label{phijtheta}
\psi=\phi_j(\rho,\theta)=\phi_j^*+\mathcal{O}(\rho)\,,
\end{equation}
for radii $|w|=\rho\searrow 0$.
For $\phi_j^*$ see \eqref{phij*}.
Trajectories in the $2(d-1)$ separating entry and exit intervals still converge to $w=0$, in the appropriate time direction, tangentially to the unperturbed asymptotic source and sink angles $\psi=\psi_j$ of \eqref{psij} in those intervals.

Analogously to section \ref{ResPol}, we observe a unique pole at $w=\infty$, albeit of increased multiplicity $d-1$.
Specifically, inversion $z=1/w$ leads to
\begin{equation}
\label{zdtheta}
\dot{z}=-z^2P(1/z)=z^{-(d-1)}(\exp(\mi\theta)-z)\,.
\end{equation}

\textbf{Case 1:} $\theta=0$.\\
Poincaré compactification \eqref{r}, \eqref{alpha} in polar coordinates $z= r \exp(-\mi\alpha)$ now reads
\begin{align}
\label{rz}
\dot{ r }\ &= r \big(\cos(d\alpha)+\ldots \big)\,;   \\
\label{alphaz}
\dot{\alpha}\ &= \ - \sin(d\alpha)+\ldots\,;
\end{align}
Compare \eqref{r}, \eqref{alpha}.
The $d$ blue unstable blow-down separatrices $\mathbf{k}$ emanate from $z=0,\  r =0$ at
\begin{equation}
\label{alphak}
\alpha_k:=k \pi/d\,,
\end{equation}
for even $0\leq k< 2d$.
The $d$ red stable separatrices $\mathbf{k}$ terminate at $z=0,\  r =0$, for the intermediate odd $k$.
For $d=4$ see figure \ref{fig6}(b).

Since the polynomial $P$ commutes with complex conjugation, the resulting phase portrait is symmetric to the real axis.
By invariance of the real axis, the blue unstable blow-down separatrix $\mathbf{0}$ emanating from the saddle $r=\alpha_0=0$ remains real and terminates at the unique nontrivial equilibrium $e_{d+1}=1$, a sink.
Also observe the real heteroclinic orbit $w: 0\leadsto e_{d+1}$.

In particular, the $2d$ hyperbolic sectors of $z=0$, prevent any separatrix loops to $z=0,\ w=\infty$, locally.
See section \ref{ODEper}.
Indeed, the mandatory residue constraint \eqref{sumetaJ0} cannot be satisfied.
For the same reason, periodic orbits and heteroclinic cycles are excluded; see \eqref{hetcycle}.

Let us now follow any of the $d-1$ remaining blue unstable blow-down separatrices $w(t)\in\mathbf{k}$, for even $k\neq 0$.
They have to converge to $w=0$ in forward time $t\nearrow+\infty$.
Hence they have to pass via the entry segments between the tangencies $w=\rho\exp(\phi_j(\rho,\theta))$; see \eqref{phijtheta}.
On the other hand, each of the $d-1$  entry segments of the small circle $|w|=\rho$ has to contain at least one of those $d-1$ blue blow-down separatrix from $w=\infty$.
Indeed, the elliptic endpoints $w=\rho\exp(\mi\phi_j^*)$ possess different sources on the $\psi$-circle of $w=0$, as their $\boldsymbol{\alpha}$-limit sets.
At least one blue separatrix has to separate their domains of attraction, in reverse time.
Therefore the $d-1$ remaining blue blow-down separatrices match the $d-1$ entry segments at $|w|=\rho$, bijectively.

Since $e_{d+1}=1$ is a sink, all  $d$ red stable blow-up separatrices $\mathbf{k}$, for odd $\mathbf{k}$, have to emerge from the only source $w=0$ via one of the $d-1$ exit segments of the small circle $|w|=\rho$: at least one per segment.
Since red and blue separatrix trajectories cannot cross (except at poles), this identifies the segment from $\phi_{2(d-1)}$ to $\phi_1$ as the only segment which can accommodate two red blow-up separatrices.
In fact those separatrices have to be $\mathbf{k}=\mathbf{1}$, passing at an angle $w=:\rho\exp(\psi_+)$, 
and $\mathbf{k=2d-1}$, passing at an angle $w=:\rho\exp(\psi_-)$.
In general, $\psi_\pm=\psi_\pm(\rho,\theta)$.
For $\theta=0$, both trajectories have to emerge tangentially to the real axis, i.e. to the slow unstable manifold of the source at $r=\psi=0$.
Therefore
\begin{equation}
\label{psik}
\psi_\pm(\rho,0)= \mathcal{O}(\rho)\,.
\end{equation}

Recall the $d$-homogeneous case \eqref{wd}.
There, all $2(d-1)$ blow-up and blow-down separatrices were strictly radial, at constant angles $\psi=\psi_j$\,, alternatingly out- and inward; see \eqref{psij}.
We now see how the insertion of the simple equilibrium $e_{d+1}=1$ interrupts the red real blow-up separatrix at $\psi_0=0$.
That separatrix gets replaced by the straddling blow-up pair $\mathbf{1},\ \mathbf{2d-1}$ which bounds a region of heteroclinic orbits $0\leadsto e_{d+1}$ and a blue blow-down separatrix $\infty\leadsto e_{d+1}$\,.
See figure \ref{fig6}(b).

\textbf{Case 2:} $\theta=\theta_\ell=\ell\pi/(d-1),\ \ell$ even.\\
As in section \ref{Ex3}, let $w^\theta$ denote any solution of \eqref{wdtheta}, for any fixed angle $\theta=\theta_\ell$\,.
Then $w:=\exp(\mi\theta)w^\theta$ solves 
\begin{equation}
\label{wtheta}
\dot{w}=\mathrm{e}^{\mi\theta}P(\mathrm{e}^{-\mi\theta}w)=(-1)^\ell\, w^d(1-w)\,.
\end{equation}
I.e., any solution $w$ of \eqref{wdtheta} for $\theta=0$ induces a solution $w^\theta$ of \eqref{wdtheta}, provided that
\begin{equation}
\label{thetal}
\theta=\theta_\ell:=\ell\pi/(d-1)
\end{equation}
and $0\leq \ell< 2(d-1)$ is even.
For $w^{\theta_\ell}$, in particular, figure \ref{fig6} becomes rotated clockwise by $\theta_\ell$\,.
For $\theta=\theta_\ell$\,, this shows the asymptotics
\begin{align}
\label{alphakl}
\alpha_k(\theta_\ell)&= \,\quad\alpha_k\,\quad -\theta_\ell=\alpha_{k-\ell}-\theta_\ell/d\,, \\
\label{psipml}
\psi_\pm(\rho,\theta_\ell)&= \psi_\pm(\rho,0)-\theta_\ell=-\theta_\ell\,+\,\mathcal{O}(\rho)\,,\\
\label{phijl}
\phi_j(\rho,\theta_\ell)&=\quad \phi_{j}^*\, +\,\mathcal{O}(\rho)= \tfrac{1}{2}(\theta_{j-1}+\theta_j)+\mathcal{O}(\rho)\,.
\end{align}
See \eqref{alphak}, \eqref{psik}.
For circle tangencies at $\phi_j(\rho,\theta_\ell)$, recall \eqref{phij*}, \eqref{phijtheta}, \eqref{thetal}.

Now, the two red blow-up separatrices $\mathbf{2d\pm 1-\ell}\  \mathrm{mod}\ 2d$, team up to replace the straddling blow-up separatrices $\mathbf{1}$ and $\mathbf{2d-1}$.
They emanate from $w=0$ along the same tangent $\psi=-\theta_\ell=\psi_{-\ell}$\,. 
I.e., they straddle heteroclinic orbits $0\leadsto e_{d+1}=\exp(-\mi\theta_\ell)$ with the same tangent at $w=0$.
By \eqref{alphakl}, the asymptotic angles $\alpha_k$ of the $2d$ blow-up orbits all shift by the same angle, counter-clockwise.
After only $d-1$ rotations by $\theta_2$\,, we recover the original phase portrait of $\theta=0$, albeit with labels $\mathbf{k}$ shifted by $2$.

\textbf{Case 3:} $\theta=\theta_\ell=\ell\pi/(d-1),\ \ell$ odd.\\
For $\theta=\theta_\ell$ given by \eqref{thetal}, as in the previous case but with odd $0<\ell<2(d-1)$, the ODE \eqref{wtheta} for $w$ becomes
\begin{equation}
\label{wtheta-}
\dot{w}=- w^d(1-w)\,.
\end{equation}
This is the time reversal of the original ODE \eqref{wdtheta}, for $\theta=0$.
We therefore obtain the analogous results as in case 2, where $\ell$ was even, except that the time direction is reversed.
This swaps the colors red and blue, as well as blow-up and blow-down separatrices.
Also, the blue sink $e_{d+1}=\exp(-\mi\theta_\ell)$ becomes a source, with reversed heteroclinic direction $e_{d+1}\leadsto 0$.

\textbf{Summary.}\\
Let us compare the two parity variants of ODE \eqref{wdtheta}, for $\theta=\theta_\ell$\,, with the $d$-homogeneous case of \eqref{wd}, as we already did for $\theta=0$ in case 1.
In cases 2 and 3, we see how the insertion of the simple equilibrium $e_{d+1}=\exp(-\mi\theta_\ell)$ then interrupts \emph{any} radial separatrix of \eqref{wd} at tangent $\psi=-\theta_\ell=\psi_{-\ell}$\,.
That separatrix gets replaced by the straddling separatrix pair with the same tangent at $w=0$.
The straddling pair bounds a region which consists of heteroclinic orbits between $0$ and $e_{d+1}$ and a single separatrix of the opposite color between $w=\infty$ and $e_{d+1}$\,.
This appends the new equilibrium $e_{d+1}$ to $w=0$, a sink or source depending on the parity of $\ell$, along any original tangent direction.

\section{Generic Morse structure}\label{Gen}

In this section we prove genericity theorem \ref{thmRatGen}: the nondegeneracy conditions \emph{(i)--(iv) } of definition \ref{defnondeg} in section \ref{Morse} have to be satisfied by generic configurations $\mathcal{E}\subset\mathbb{C}^{2d-2}$ of $d$ zeros $e_j\in E$ and $d'=d-2$ poles $e_{-j'}\in E'$ of $f=P/Q$.
For the regularized flow \eqref{ODEreg}, corollary \ref{corMorse} then implies the Morse properties of definition \ref{defMorse} and proposition \ref{propgrad} from sections \ref{DefMorse}, \ref{Gradlike}.

\emph{Proof of theorem} \ref{thmRatGen}.\\
\emph{Properties (i).}\quad
We only have to choose $E\cup E'\subset\mathbb{C}$ to consist of $d+d'=2d-2$ distinct elements $e_j\,,\ e_{-j'}$\,.
Each degeneracy $e_j=e_{-j'}$ is linear of real codimension 2.
Avoiding them defines a subset $\mathcal{E}_{(i)}\subset\mathbb{C}^{2d-2}$ which is open and dense, hence generic.

\emph{Properties (ii).}\quad
Pairs of poles $e_{-j'}$ define at most $\tfrac{1}{2}d'(d'-1)$ straight lines $\Gamma_0$ of real dimension 1 in $\mathbb{C}$.
Avoiding these lines by $E$ defines a further subset of configurations $\mathcal{E}_{(ii)}\subset\mathcal{E}_{(i)}$ which is open and dense, hence generic.

\emph{Properties (iii).}\quad
Let $\mathcal{E}_{(iii)}\subset\mathcal{E}_{(ii)}$ denote the set of configurations which satisfy the $2^d-2$ real nondegeneracy conditions \eqref{sumetaJno}.
By continuity, $\mathcal{E}_{(iii)}$ is open.

To show density of $\mathcal{E}_{(iii)}$, fix any subset $\emptyset\neq J\subsetneq\{1,\ldots,d\}$ and  any $j_*\not\in J$.
Let 
\begin{equation}
\label{gej*}
g(e_{j_*}):=\sum_{j\in J} \eta_j(e_{j_*})
\end{equation}
denote the nonempty residue sum in \eqref{sumetaJno}.
The notation highlights the dependency of the residues $\eta_j=1/f'(e_j)$ on $e_{j_*}$\,, with $j_*\neq j$.
To achieve nondegeneracy (ii), we first perturb $e_{j_*}$ such that $\Re g(e_{j_*})\neq 0$ for
\begin{equation}
\label{1/eta}
1/\eta_j(e_{j_*})\,=\,f'(e_j)\,=\, a(e_j-e_{j_*})\,\frac{1}{Q(e_j)}\,\prod_{k\neq j,j_*}(e_j-e_k)\,.
\end{equation}
See \eqref{linequi}, \eqref{omega,eta}.
Since $j\in J\not\ni j_*$\,, the function $e_{j_*}\mapsto g(e_{j_*})$ is rational in $e_{j_*}$ with simple poles at $e_j\,,\ j\in J$.
The open mapping theorem for the nonconstant function $g$ therefore implies that $\Re g(e_{j_*})\neq 0$ can always be achieved by arbitrarily small perturbations of $e_{j_*}$\,.
By finite intersection over all $2^d-2$ nontrivial subsets $J$, Baire's theorem shows that $\mathcal{E}_{(iii)}$ is open and dense.

\emph{Properties (iv).}\quad
Fix any period $c_*=2\pi\mi\,\sum m_j\eta_j \in\mathcal{P}\,:=\,2\pi\mi\langle \eta_1,\ldots,\eta_d\rangle_\mathbb{Z}$\,; see \eqref{calP}.
Without loss, fix any two distinct poles $e_{-1}\,, e_{-2}\in E'$ with the straight line $\Gamma_0$ between them.
By continuity, the set $\mathcal{E}_*\subset\mathcal{E}_{(iii)}$ of configurations which satisfy nondegeneracy condition \eqref{hetTPno} for that particular period $c_*$\,, rather than the full countable but possibly dense set $\mathcal{P}$, is open.

To show density of $\mathcal{E}_*$\,, i.e.~of $\int_{\Gamma_0}\omega\not\in c_*+\mathbb{R}$, we will push $T=T(e_1)$ in \eqref{sovpath}  off the real axis, by perturbations of $e_1$\,.
In detail:
\begin{equation}
\label{hetTPno1}
T(e_1):= c_*+\int_{\Gamma_0}\omega \,=c_*+ \sum_{j=1}^d \eta_j \int_{e_{-1}}^{e_{-2}} dw/(w-e_j)\,=\,b_1(e_1)g(e_1)+b_0(e_1)\,.
\end{equation}
With \eqref{1/eta} and the abbreviations $\lambda_j:=\log(e_{-2}-e_j)-\log(e_{-1}-e_j)$, we have separated the sum over $j>1$ as $b_0(e_1)$ here, to collect
\begin{equation}
\label{ge1}
\begin{aligned}
b_0(e_1)\ &:=\,\ \sum_{j>1} c_j/(e_j-e_1)\,;\\
b_1(e_1)\ &:=\, \phantom{\Big(}a^{-1}\prod_{k'} (e_1-e_{-k'})\, \Big/ \ \prod_{k>1}(e_1-e_k)\,;\\
c_j\ \ \ &:=\, a^{-1}\Big(\prod_{k'} (e_j-e_{-k'})\, \Big/  \prod_{1<k\neq j}(e_j-e_k)\Big)\cdot(2\pi\mi\, m_j+\lambda_j)\quad \mathrm{for}\ j>1\,;\\
g(e_1)\ &:=\,\ 2\pi\mi\, m_1+\lambda_1(e_1)\,.
\end{aligned}
\end{equation}
Recall how $m_1\,, m_j$ account for $c_*\in\mathcal{P}$.
Also see \eqref{omega,eta} and \eqref{1/eta}, for $j_*:=1<j$.

With respect to $e_1$\,, expressions $c_j$ are constant, $b_0(e_1)$ and $b_1(e_1)$ are rational on $\mathbb{C}$, and $g(e_1)$ is holomorphic on the Riemann surface $\mathcal{R}$ of the logarithm $\lambda_1(e_1)$ over $\mathbb{C}\setminus\{e_{-1},e_{-2}\}$.
Since $T(e_1)$ is nonconstant meromorphic on $\mathcal{R}$, the open mapping theorem provides
arbitrarily small perturbations of $e_1$ which push $T(e_1)-c_*$ off the real axis.
This shows density of $\mathcal{E}_*\subset\mathcal{E}_{(iii)}$\,.

To cover all periods $c_*\in\mathcal{P}$ in \eqref{hetTPno}, we define $\mathcal{E}_{(iv)}\subset\mathcal{E}_{(iii)}$ as the countable intersection of all associated open dense configurations $\mathcal{E}_*$\,.
Then $\mathcal{E}_{(iv)}$ is generic, by Baire's theorem, and satisfies all nondegeneracy conditions \emph{(i)--(iv)} of definition \eqref{hetTPno}.
This proves genericity theorem \ref{thmRatGen}. \hfill $\bowtie$

\section{Proof of rational theorem \ref{thmPortrait2Rat}}\label{PfthmRat}

\subsection{Proof, up to outer and inner expansions}\label{OutthmRat}
To prove realization theorem \ref{thmPortrait2Rat} for rational vector fields, we proceed by induction on the degree $d$ of the numerator polynomial $P$, from $d\geq 2$ to $d+1$.
The quadratic case $d=2,\ d'=0$ is trivial; see section \ref{Ex2}.

\begin{figure}[t]
\centering \includegraphics[width=0.9\textwidth]{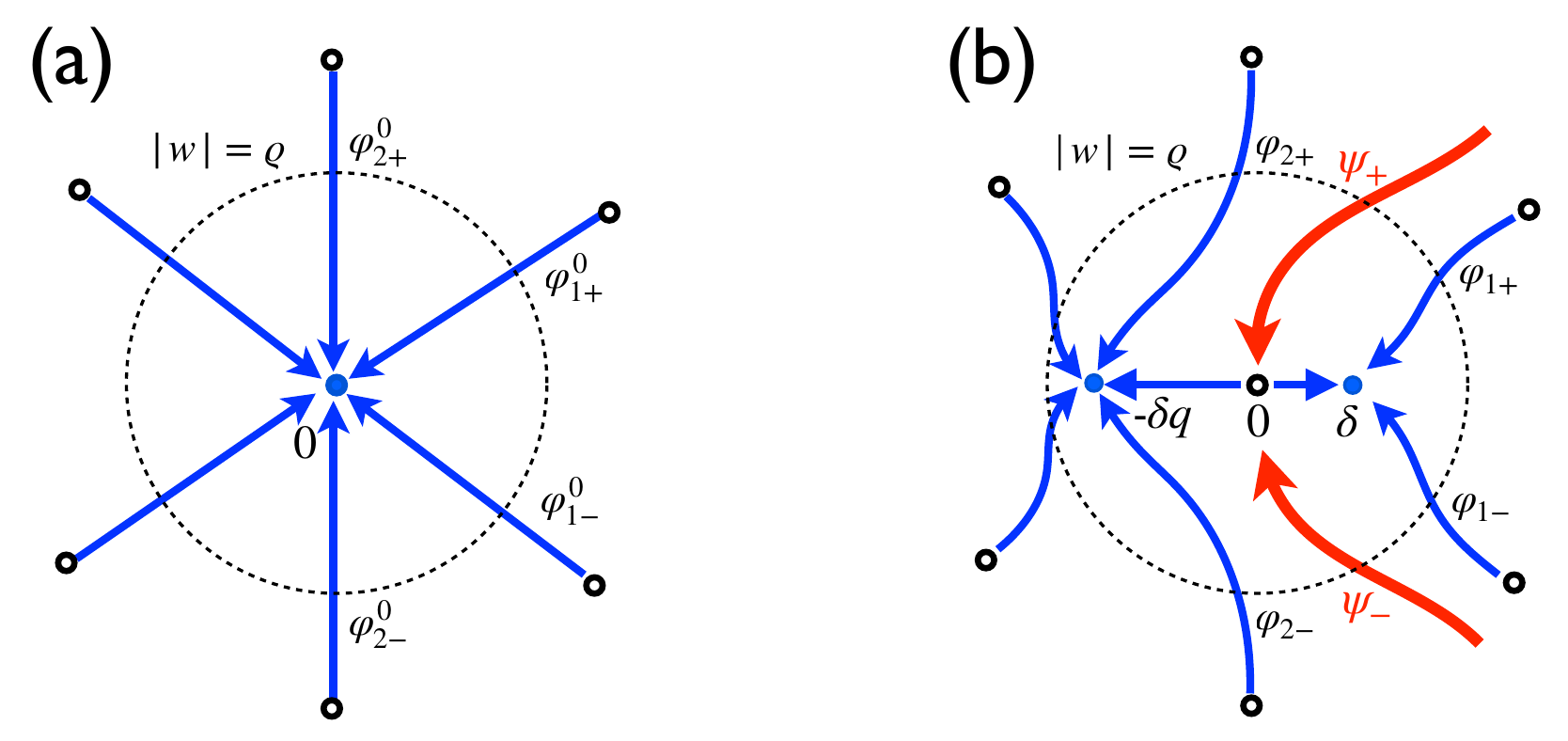}
\caption[Induction over rational degree]{\emph{
Illustration of induction \eqref{fd}, \eqref{fd+1} over degree $d$, by edge contraction for $\theta=0$.\\
(a) Local phase portrait \eqref{fd} after edge contraction, $\dot{w}=f_d=awg(w),\ g(0)=1$. 
Note the sink $w=0$ (blue) with incoming blue unstable separatrices.
These originate from certain nonzero poles of $g$ elsewhere (black circles).
Just like any other trajectories, they transversely enter a fixed small circle $|w|=\rho>0$ (black dashed) around the sink $w=0$, at angles $\phi_{j'}^0$\,.\\
(b) Local phase portrait prior to edge contraction, $\dot{w}=f_{d+1}=h^\delta(w) g(w)$.
The version $h^\delta$ of the rational factor $h$ is downscaled by $\delta=\eps\rho>0$; compare the cubic example of section \ref{Ex3} and \eqref{hdelta}--\eqref{fd+1}.
The two red stable separatrices, which emanate from the simple pole $w=0$ of $f_{d+1}$ and $h^\delta$\,, in backward time, leave the circle $|w|=\rho$ at angles $\psi_\pm$.
By lemma \ref{lempsi}, the angles $\psi_\pm$ can be chosen to separate the incoming angles $\phi_{1\pm}$ and $\phi_{2\pm}$ from each other, as indicated.
Therefore the incoming separatrices at $\phi_{j'}^\delta$ can be ordained to converge to the two sinks $w=\delta e_1=\delta$ and to $w=\delta e_2=-\delta q$ (blue), respectively, for $j'$ in $J_1'$ and $J_2'$\,.
See the outer and inner expansions \eqref{ephi}--\eqref{elimpsi0}.
}}
\label{fig7}
\end{figure}

Consider any given, mutually dual pair $\widetilde{\mathcal{C}}_{d+1}^\pm$ of nonempty, finite, connected multi-graphs on the $2$-sphere $\mathbb{S}^2$ with $d'+1\geq 1$ edges, each, and a total of $d+1=(d'+1)+2$ vertices.
After time reversal, and an orientation preserving homeomorphism if necessary, we may assume that  $\widetilde{\mathcal{C}}_{d+1}^-$ possesses a blue edge through the pole $w=0$ which is not a loop, but joins two blue vertices $e_1>0>e_2$\,.
See figure \ref{fig7}(b).
Contracting the blue edge collapses its endpoints $e_1,e_2$ to a distinguished single blue sink vertex, say at $w=0$. 
We obtain a contracted multi-graph  $\widetilde{\mathcal{C}}_d^-$ with only $d'=d-2$ edges.
In the red dual $\widetilde{\mathcal{C}}_d^+$ , the corresponding red edge of $\widetilde{\mathcal{C}}_{d+1}^+$ through $w=0$ simply gets removed, and its red source vertices remain untouched.
See figure \ref{fig7}(a).

By the induction hypothesis, the dual pair $\widetilde{\mathcal{C}}_d^\pm$ is realized by a nondegenerate flow \eqref{ODE} with rational nonlinearity $f=f_d=P/Q$, such that the red unstable and blue stable portraits $\mathcal{C}_d^\pm$ on $\mathbb{S}^2=\widehat{\mathbb{C}}$ are isomorphic to the given multi-graphs $\widetilde{\mathcal{C}}_d^\pm$, with orientation preserving homeomorphic embeddings.
In symbols, $\mathcal{C}_d^\pm\equiv\widetilde{\mathcal{C}}_d^\pm$\,.
The degrees of the polynomials $P,Q$ with disjoint simple zero sets $E,E'$ are $d$ and $d'=d-2$, respectively.
After a Möbius transformation of $w$, we may assume the collapsed blue sink vertex of $\widetilde{\mathcal{C}}_d^-$ to occur at $w=0$.
This allows us to factorize $f_d(w)=a\,wg(w)$, with $a:=f_d'(0),\ \Re a<0$, and rational $g$.

To realize the full original graphs $\widetilde{\mathcal{C}}_{d+1}^\pm$ by homeomorphic portraits $\mathcal{C}_{d+1}^\pm$ of $\dot{w}=f_{d+1}(w)$, we substitute the factor $w$ of $f_d$ by a downscaled variant $h^\delta(w)$ of $h$ from the cubic example \ref{Ex3}; see \eqref{cubic}.
The scaling $h^\delta$ of $h$ is real linear,
\begin{equation}
\label{hdelta}
h^\delta(w):= \delta h(w/\delta) = a\,(w-\delta e_1) (w+q\delta e_1)/w\,.
\end{equation}
In other words, $\dot{w}=h(w)$ in the cubic example \eqref{cubic} becomes $\dot{w}^\delta=h^\delta(w^\delta)$, for downscaled $w^\delta:=\delta w$.
We consider small $\eps>0$ such that the unscaled $h$ satisfies lemma \ref{lempsi}.
Define $h^0(w):=aw$ and $f^\delta:=h^\delta g$.
Then $f_d$ becomes $f^0$ and we choose $f_{d+1}:=f^\delta$\,, for some
\begin{equation}
\label{deltaepsrho}
\delta:=\eps\rho>0\,.
\end{equation}
This two-scale expansion with respect to small $\rho$ and $\eps$ will downscale the large circle $|w|=1/\eps$ to $|w^\delta|=\rho$.
The two induction-related ODEs now read
\begin{align}
\label{fd}
    \dot{w}=\,\ f_d(w)\,\ &:=aw\,g(w)=h^0(w)\,g(w)\,,\quad \mathrm{versus}  \\
\label{fd+1}
    \dot{w}=f_{d+1}(w)&:= h^\delta(w)\,g(w)\,. 
\end{align}
It remains to fix $q>0$ and $\theta\in\mathbb{R}$, in the definition \eqref{cubic} of unscaled $h$, such that \eqref{fd+1} realizes the prescribed pre-contraction portrait $\mathcal{C}_{d+1}^-\equiv\widetilde{\mathcal{C}}_{d+1}^-$\,, up to an orientation preserving homeomorphism.

Consider circles $|w|=\rho$ of small radii $\rho>0$.
In particular, the $\rho$-circle should isolate the sink $w=0$ from all other equilibria $e_j$ and poles $e_{-j'}$ of $f_d$\,.
For $\rho\searrow 0$, we observe
\begin{equation}
\label{rhod}
f_d(w)/(aw)=g(w)=1+\mathcal{O}(\rho)\,.
\end{equation}
Analogously to \eqref{cubictrv}, therefore, $f_d$ points strictly inward on the circle $|w|=\rho$,
\begin{equation}
\label{fdtrv}
\langle w,f_d(w)\rangle \,=\,\Re(\, \overline{w}\,f_d(w))\,=\,\Re a\,\rho^2+\mathcal{O}(\rho^3)\,<\,0\,,
\end{equation}
for small $\rho$.
In particular, let $w_{j'}=\rho\exp(\mi\phi_{j'}^0),\ j'\in J'$, define the angles at which any blue incoming separatrices $e_{-j'}\leadsto 0$ of \eqref{fd} from outside pole vertices $e_{-j'}\in E_d'$, alias the blue edges of $\mathcal{C}_d^-$ towards the sink $w=0$, cross the $\rho$-circle.
By definition of the contracted graph $\mathcal{C}_d^-\equiv \widetilde{\mathcal{C}}_d^-$ of $f_d$\,, the crossing set $J'$ decomposes into two disjoint parts, $J'=J_1'\cup J_2'$\,: 
the indices $j'\in J_1'$ of blue separatrices in $\widetilde{\mathbb{C}}_{d+1}^-$ which terminate at $e_1$\,, and the  $j'\in J_2'$ towards $e_2=-qe_1$\,.
In counter-clockwise orientation of the crossing angles $\phi_{j'}^0$\,, let $j'=1-$ and $j'=1+$ denote the first and last elements of $J_1'$\,, respectively. 
For $j'\in J_2'$\,, we denote the corresponding elements by $j'=2\pm$, but in clockwise orientation.
The associated angles $\phi_{j'}^0$ are $\phi^0_{1\pm}$ and $\phi^0_{2\pm}$\,.
Here and below we assume that $J_1'\,, J_2'$ contain at least two elements, each.
See figure \ref{fig7}.
The  easy adaptations to one or no elements are left to the reader.
In section \ref{OuterRat} we will invoke the implicit function theorem to establish the following \emph{outer expansion} of the crossing angles
\begin{equation}
\label{ephi}
\phi_{j'}(\eps,\rho,q,\theta)=\phi_{j'}^0(\rho)+\mathcal{O}(\eps)=\phi_{j'}^*+\Theta(\rho)+\mathcal{O}(\eps)+\mathcal{O}(\rho)\,,
\end{equation}
for $f_{d+1}=f^\delta$; see lemma \ref{lemOuterRat}.
Here $\eps,\rho,q,\theta$ and $\delta=\eps\rho$ explicitly account for the dependence on the factor $h^\delta$ of $f^\delta$, and its details.
The expansion is locally uniform in $q>0$.
The distinct constants $\phi_{j'}^*$ depend on $j'$, only.
Note that the uniform rotation by $\Theta=\Theta(\rho)$ is independent of $j'$ and $\delta,q,\theta$.
Only the higher orders $\mathcal{O}(\rho)$ restore that dependence.

For the crossing angles $\psi_\pm$ of the red stable separatrices towards the pole $w=0$ of $f_{d+1}=f^\delta$\,, section \ref{InnerRat} will provide analogous \emph{inner expansions}
\begin{align}
\label{epsitheta}
\psi_\pm(\eps,\rho,q,\theta))&=\psi^0_\pm(\eps,q,0)+\theta+\mathcal{O}(\rho)\,,\\
\label{elimpsi0}
\psi_+(\eps,\rho,q,\theta))-\psi_-(\eps,\rho,q,\theta))&=2\pi/(1+q)+\mathcal{O}(\rho)\,; 
\end{align}
see lemmata \ref{lempsi} and \ref{lemInnerRat}.

For $f_{d+1}=f^\delta$, we now combine the inner and outer expansions to make the separatrices for $|w|\leq\rho$ correspond to the preordained configurations $\widetilde{\mathcal{C}}_{d+1}^\pm$ of figure \ref{fig7}(b).
Let $c$ denote the minimal length of the clockwise intervals $I_-:=(\phi_{2-},\phi_{1-})$ and $I_+:=(\phi_{1+},\phi_{2+})$ which separate $J_1'$ and $J_2'$\,.
We choose $\eps,\rho$ small enough such that any remainders satisfy
\begin{equation}
\label{Orhoetc}
\mathcal{O}(\eps)<c/9\,,\qquad \mathcal{O}(\rho)<c/9\,.
\end{equation}
We then take the difference of the outer and inner expansions \eqref{ephi} and \eqref{epsitheta}, \eqref{elimpsi0} as follows.
Given $\rho$, we define $q,\theta$ such that $\pm\pi/(q+1)+\theta$ become the mid-points of the separating intervals $I_\pm=I_\pm^*+\Theta(\rho)$, respectively, each of length at least $c$.
The perturbation estimates \eqref{Orhoetc} then guarantee, at parameters $\eps,\rho,q,\theta$, that the entry angles $\psi_\pm$ of the two red separatrices towards the inserted pole $w=0$ of $f^\delta=f_{d+1}$, i.e. of the new red edge in $\mathcal{C}_{d+1}^+$, still separate the entry angles $\phi_{j'}$ of the blue separatrices, as prescribed by $\widetilde{\mathcal{C}}_{d+1}^+$ and the decomposition $J'=J_1'\cup J_2'$ of their index set.
Indeed, $4/9<1/2$.
This implies convergence of blue separatrices to the sinks $\delta e_1$ and $\delta e_2=-\delta qe_1$ of $f_{d+1}$ in ${\mathcal{C}}_{d+1}^-$\,, respectively, as prescribed by $\widetilde{\mathcal{C}}_{d+1}^-$\,.
The remaining blue separatrix edges of the new pole $w=0$ cannot leave the $\rho$-circle, by persistence of transversality $\langle w, f_{d+1}(w) \rangle<0$; see \eqref{fdeltatrv}.
Therefore they have to connect to the blue sinks $\delta e_1$ and $\delta e_2$ in ${\mathcal{C}}_{d+1}^-$\,, again as prescribed by $\widetilde{\mathcal{C}}_{d+1}^-$\,.
Up to the outer and inner expansions \eqref{ephi} and \eqref{epsitheta}, \eqref{elimpsi0}, which we prove next, this establishes the preordained separatrix configuration of figure \ref{fig7}(b), and proves theorem \ref{thmPortrait2Rat}.

\subsection{Outer expansion}\label{OuterRat}

In this section we prove the outer expansion \eqref{ephi} for crossing angles of separatrices which enter the circle $|w|=\rho$ from outside, i.e. from $|w|\geq\rho$.
Throughout we choose $\delta=\eps\rho$\,; see \eqref{deltaepsrho}.
We consider small radii $\rho>0$.
In particular all poles $e'\in E_d'$ and all nonzero equilibria $e\in E_d\setminus\{0\}$ of $f_d=aw\,g(w)$ are strictly outside the $\rho$-circle.
The scaling parameter $\eps$ will also be small, as in lemma \ref{lempsi}.
By construction, and for $\eps<1,1/q$, exterior zeros and poles of $f^0=f_d$ and $f_{d+1}=f^\delta=h^\delta\,g$ coincide, because they can all be attributed to $g$.
The only remaining interior zeros and poles of $f_{d+1}$ are
\begin{equation}
\label{EdE'd}
0\in E_d\cap E_{d+1}',\quad\mathrm{and}\quad e_1^\delta\,,\ e_2^\delta\in E_{d+1}\,.
\end{equation}
Here $e_1^\delta=\delta\exp(\mi\theta)$ and $e_2^\delta=-qe_1^\delta$\,.
See section \ref{Ex3} and \eqref{hdelta}--\eqref{fd+1}.

Omitting indices, let $e'\in E'$ denote any exterior pole of $f_d=aw\,g(w)$, i.e. of $g$, and hence also of $f_{d+1}=h^\delta\,g$\,.
Assume $e'$ possesses a blue unstable separatrix $w: e'\leadsto 0\in E_d$ for $f_d$\,. 
Under the flow of $f_{d+1}=f^\delta$ for $\delta=\eps\rho$, we denote the angle $\phi$ where the separatrix crosses the circle $w(t)=\rho\exp(\mi\phi)$ by
\begin{equation}
\label{phirho}
\phi=\phi(\eps,\rho,q,\theta)\,.
\end{equation}
To define any unique crossing angles of separatrices, we first note that inward crossing of all separatrices of \eqref{fdtrv} for $f_d=f_0$ persists by transversality.
Indeed
\begin{equation}
\label{fdeltatrv}
\langle w,f^\delta(w)\rangle \,=\,\Re\big(  \overline{w}\,h^\delta(w)g(w)\big) \,=\,\delta\, \Re\big(  \overline{w}\,h(w/\delta)\big)(1+\mathcal{O}(\rho))\, <\,0
\end{equation}
on the $\rho$-circle $|w|=\rho=\delta/\eps$\,, provided that $\eps$ and $\rho$ are chosen small enough; see \eqref{cubictrv}, \eqref{hdelta}, \eqref{deltaepsrho}, \eqref{rhod}.

Since the separatrix is traversed in real time $t$, solution \eqref{sov} by separation of variables provides an explicit equation for the crossing angle $\phi$.
In fact $\phi$ solves the real equation
\begin{equation}
\label{Phi}
\Phi(\phi,\eps,\rho,q,\theta) := \Im \int_{e'}^{\rho\exp(\mi\phi)} dw/f^{\eps\rho}(w) \,=\, 0\,,
\end{equation}
provided the integration path $\Gamma$ remains homotopic to the separatrix itself.

\begin{lem}\label{lemOuterRat}
For small $\rho\searrow 0$, let $\phi^0=\phi^0(\rho)$ denote the crossing angle of a separatrix from the pole $e'$ to the sink $w=0$ for the nonlinearity $f_d=f^0=aw\,g$\,.

Then $\phi=\phi^0(\rho)$ is a solution of \eqref{Phi} for $\eps=0$ and any $q,\theta$, with an expansion
\begin{align}
\label{phi0rho}
\phi^0(\rho)&=\phi^*+\Theta(\rho)+\mathcal{O}(\rho)\,,\\
\label{Thetarho}
\Theta(\rho)&=\tfrac{\Im a}{\Re a}\, \log\rho +\mathcal{O}(\rho)\,.
\end{align}
While the constant term $\phi^*$ in expansion \eqref{Thetarho} depends on the particular separatrix, the leading logarithmic spiral dependence $\Theta(\rho)$ of $\phi^0$ on $\rho$ does not.

Moreover, $\Phi$ is $C^1$ in a neighborhood and the implicit function theorem provides the unique angle $\phi=\phi(\eps,\rho,q,\theta)$, locally, at which the blue separatrix of $e'$ for $f_{d+1}=f^\delta=h^\delta g,\ \delta=\eps\rho$ enters the circle $|w|=\rho$.
For $\eps\searrow 0$, in particular, this proves
\begin{equation}
\label{phieps}
\phi(\eps,\rho,q,\theta)=\phi^0(\rho)+\mathcal{O}(\eps)\,.
\end{equation}
Combining \eqref{phi0rho} and \eqref{phieps} proves the outer expansion \eqref{ephi}.
\end{lem}

\begin{proof}
To show that the reference solution $\phi=\phi^0$ solves \eqref{Phi} for $\eps=0$, just note that $\delta=\eps\rho=0$ and $h^0=aw$.
Local differentiability for $q>0,\ \theta,\psi\in\mathbb{R}$ and small $|\eps|\geq0,\ \rho>0$ follows because $|w|\geq\rho$, \eqref{hdelta}, and $\delta:=\eps\rho$ imply a denominator
\begin{equation}
\label{eps=0}
f^\delta(w)=h^{\delta}(w)g(w)= a\,(w-\eps\rho e_1) (w+\eps\rho q e_1)\,g(w)/w
\end{equation}
which remains regular for $\rho>0$ and small $|\eps|$.
For the integration path $\Gamma$ from $e'$ to $w=\rho\exp(\mi\phi)$ with $\phi$ near $\phi^0$ we may choose to follow the separatrix of $f_d$\,, to $\phi=\phi^0$, and an arc on the $\rho$-circle afterwards.

To show logarithmic spiraling \eqref{Thetarho} of the reference solution, we holomorphically linearize ODE \eqref{ODE} for $f=f_d$\,, locally at $w=0$.
Since $f_d'(0)=a$, near-identity linearization by $w\mapsto \tilde{w}=w+\ldots$ yields $\tilde{w}(t)=\tilde{w}(t^*) \exp(a(t-t^*))$.
In polar coordinates for $\tilde{w}$, this implies claim \eqref{Thetarho} \emph{without} higher order terms.
Reverting to polar coordinates for $w$ proves the claim, albeit \emph{with} higher order terms.

To expand $\phi$ for small $\eps>0$, we apply the implicit function theorem to \eqref{Phi}.
At the reference solution $\phi=\phi^0(\rho),\ w=\rho\exp(\mi\phi^0),\ \eps=0$, the partial derivative $\Phi_\phi$ satisfies
\begin{equation}
\label{Phiphi}
\begin{aligned}
   \Phi_\phi=\Im\Big(  \, \mi w/f^\delta(w)\Big) &=\Re\Big( \,  w^2\,/\,\big(a\,(w-\eps\rho e_1) (w+\eps\rho q e_1)\,g(w)\big)\Big)=\\
      &=\Re\Big( \,\,w^2/\,(aw^2\,g(w))\Big) =\Re (1/a)+\mathcal{O}(\rho)\neq 0.
\end{aligned}
\end{equation}
We evaluate $\Phi_\eps$ next:
\begin{equation}
\label{Phieps}
\begin{aligned}
   \Phi_\eps&=\Im\int_{e'}^{w}  \partial_\eps \big(1/h^\delta(\hat{w})\big)\,d\hat{w}/g(\hat{w})=\\
		&=\Im \int_{e'}^{w}  \partial_\eps \big(\hat{w}/(a\,(\hat{w}-\eps\rho e_1) (\hat{w}+\eps\rho q e_1))\,d\hat{w}/g(\hat{w})=\\
		&=(1-q)\,\Im\Big( \, e_1\int_{e'}^{\rho\exp(\mi\phi^0)} \rho \, d\hat{w}\,/\big(a\hat{w}^2 g(\hat{w})\big)\Big) \,.
\end{aligned}
\end{equation}
We have to determine the asymptotics of the last integral $I=I(\rho)$, for $\rho\searrow 0$.
Except for a finite constant contribution, say from $e'$ to some fixed $\rho_0>0$, we may invoke \eqref{rhod} to replace $g$ by $1$ and evaluate the integral:
\begin{equation}
\label{Irho}
\begin{aligned}
I(\rho)&=\mathcal{O}(1)+ \int^{\rho\exp(\mi\phi^0)} \rho\, d\hat{w}\,/\big(a\hat{w}^2)\big)=\\
	&= \mathcal{O}(1) - \exp(-\mi\phi^0)/a = \mathcal{O}(1)\,.
\end{aligned}
\end{equation}
This proves boundedness of $\Phi_\eps$\,,
\begin{equation}
\label{Phiepsbd}
\Phi_\eps=\mathcal{O}(1)\,,
\end{equation}
uniformly for $\rho\searrow 0$ and along any separatrix $e'\leadsto 0$ of $f_d$\,.
Combining estimates \eqref{Phiphi} and \eqref{Phiepsbd}, the implicit function theorem applies to \eqref{Phi} and provides the local crossing angles $\phi=\phi(\eps,\rho,q,\theta)$ for the perturbed blue separatrices of $f_{d+1}=f^\delta$.
This proves the lemma. 
\end{proof}

\subsection{Inner expansion}\label{InnerRat}

Complementing the previous section, we now prove the inner expansions \eqref{epsitheta}, \eqref{elimpsi0} for the crossing angles $\psi_\pm$ of the two red stable separatrices towards the new pole $w=0$ which leave the disk $|w|\leq\rho$ through the circle $|w|=\rho$ for good, transversely and in backward time.
See transversality property \eqref{fdeltatrv}.
The nonlinearity remains $f_{d+1}=f^\delta=h^\delta g$ for small $\eps,\rho\searrow 0$ and $\delta=\eps\rho$, as before; see \eqref{deltaepsrho}.
Substituting the pole $w=0$ for $e'$, and $\psi=\psi_\pm$ for $\phi$ in \eqref{Phi}, we now have to solve
\begin{equation}
\label{Psi}
\begin{aligned}
    \Psi(\psi,\eps,\rho,q,\theta)\,&:=\,\Im \int_0^{\eps^{-1}\exp(\mi\psi)} \Big(\hat{w}/(a\,(\hat{w}-e_1) (\hat{w}+q e_1)\Big)\,d\hat{w}/g(\eps\rho\hat{w})\,=\\ 
   &\phantom{:}=\, \Im \int_0^{\rho\exp(\mi\psi)} dw/f^\delta(w) \,= \,0 \,.
\end{aligned}
\end{equation}
Here we have substituted $\hat{w}=w/\delta$ in \eqref{eps=0}, without affecting angles.
As in section \ref{OuterRat}, we integrate along paths $\Gamma$ which remain homotopic to the reference separatrices themselves.
Dependence on $\theta\in\mathbb{R}$ is hidden in $e_1=\exp(\mi\theta)$.

\begin{lem}\label{lemInnerRat}
For small $\eps>0$ in the cubic ODE example $\dot{w}=h(w)$ of \eqref{cubic}, let $\psi_\pm^0=\psi_\pm^0(\eps,q,\theta)$ denote the crossing angles, at the large circle $|w|=1/\eps$, of the two red stable separatrices which terminate at the pole $w=0$\,; see section \ref{Ex3}.

Then $\psi=\psi_\pm^0$ are solutions of \eqref{Psi} for $\rho=0$ and any $\eps,q,\theta$.
By lemma \ref{lempsi}, they satisfy equivariance property \eqref{psitheta} in $\theta$, and $\eps$-expansion \eqref{limpsi0}.

Moreover, $\Psi$ is $C^1$ in a neighborhood of these solutions, and the implicit function theorem provides the unique angles $\psi_\pm=\psi_\pm(\eps,\rho,q,\theta)$, locally, at which the red separatrices of $f_{d+1}=f^\delta=h^\delta g,\ \delta=\eps\rho$ towards the pole $w=0$ enter the circle $|w|=\rho>0$.
For $\rho \searrow 0$ and uniformly for small $\eps>0$, in particular, this proves
\begin{equation}
\label{psieps}
\psi_\pm(\eps,\rho,q,\theta)=\psi_\pm^0(\eps,q,\theta)+\mathcal{O}(\rho)\,.
\end{equation}
Combining \eqref{psieps} with \eqref{psitheta}, \eqref{limpsi0} of lemma \ref{lempsi}, where $\psi_\pm=\psi_\pm^0(\eps,q,\theta)$, proves the inner expansions \eqref{epsitheta}, \eqref{elimpsi0}.
\end{lem}

\begin{proof}
To show that the reference solutions $\psi=\psi_\pm^0=\psi_\pm^0(\eps,q,\theta)$ solves \eqref{Psi} for $\rho=0$, first note $g(\eps\rho\hat{w})=1$ for $\rho=0$\,; see \eqref{rhod}.
By separation of variables \eqref{sov}, therefore, the first integral in \eqref{Psi} represents \emph{real} time along the two separatrices $\gamma_\pm$ of \eqref{cubic}, and the \emph{imaginary} parts must vanish.

To expand $\psi$ for small $\rho>0$, we apply the implicit function theorem to \eqref{Psi}.
Local differentiability, for $q>0$, real $\theta,\psi$, and small $|\rho|\geq 0,\ \eps>0$, follows from the substitution $\hat{w}=w/\eps$ in the first integral:
\begin{equation}
\label{Psi2}
\Psi(\psi,\eps,\rho,q,\theta)\,=\,\Im \int_0^{\exp(\mi\psi)} w/\big(a\,(w-\eps e_1) (w+\eps q e_1)\big)\,dw/g(\rho w)
\end{equation}
in the unit disk $|w|\leq 1$.

At the reference solution $\psi=\psi_\pm^0(\eps,q,\theta)$ for $\rho=0$, we determine the partial derivatives $\Psi_\psi$ and $\Psi_\rho$ as follows.
First, \eqref{Psi2} at $w:=\exp(\mi\psi),\ \rho=0,\ g(\rho w)=1$ implies
\begin{equation}
\label{Psipsi}
   \Psi_\psi\,=\,\Im\Big(  \, \mi w^2\,/\,\big(a\,(w-\eps e_1) (w+ \eps q e_1)\big)\Big) 
       \,=\,\Re (1/a)+\mathcal{O}(\eps)\neq 0\,.
\end{equation}
Next we evaluate $\Psi_\rho$ at $\rho=0,\ g(\rho w)=1$.
We claim that
\begin{equation}
\label{Psirho}
   \Psi_\rho\,
                   =\,-\,\Im\Big( \, \tfrac{g'(0)}{a}\int_0^{\exp(\mi\psi)} w^2/\big((w-\eps e_1) (w+\eps q e_1)\big)\,dw\Big) \,=\,\mathcal{O}(1)
\end{equation}
remains bounded, uniformly for $\theta\in\mathbb{R}$ and small $\eps>0$, and locally uniformly for $q>0$.
To evaluate the integral and show boundedness, we use explicit complex integration of the resulting partial fractions:
\begin{equation}
\label{Psirhobd}
\begin{aligned}
    &\int_0^{\exp(\mi\psi)}  w^2/\big((w-\eps e_1) (w+\eps q e_1)\big)\,dw\,=   \\
    =\,&\int_0^{\exp(\mi\psi)}  \Big( 1+\tfrac{\eps'}{1+q}\big(1 /(w-\eps')-q^2/(w+q \eps')\big)\Big)\,dw\,=  \\ 
    =\,&\quad \exp(\mi\psi) + \tfrac{\eps'}{1+q}\big(\log(w-\eps') -q^2\log (w+q \eps') \big) \bigg\vert_0^{\exp(\mi\psi)}
\end{aligned}
\end{equation}
Here we have substituted $\eps':=\eps e_1$.
The complex logarithms have to be evaluated up to some fixed integer multiple of $2\pi\mi$, depending on a Riemann cut.
This proves boundedness \eqref{Psirho} of $\Psi_\rho$\,.

Combining estimates \eqref{Psipsi} and \eqref{Psirho}, the implicit function theorem provides the local crossing angles $\psi=\psi_\pm(\eps,\rho,q,\theta)$ for the perturbed red stable separatrices of $f_{d+1}=f^\delta$, viz the solutions $\psi$
of $\Psi=0$\,; see \eqref{Psi}. Moreover $\psi_\pm(\eps,\rho,q,\theta)=\psi_\pm^0(\eps,q,\theta)+\mathcal{O}(\rho)\,$, as claimed in \eqref{psieps}.
This proves the lemma. 
\end{proof}

With the proofs of lemmata \ref{lemOuterRat} and \ref{lemInnerRat}, the outer and inner expansions \eqref{ephi} and \eqref{epsitheta}, \eqref{elimpsi0} are now verified.
By section \ref{OutthmRat}, this completes the proof of realization theorem \ref{thmPortrait2Rat}.

\section{Proof of polynomial theorem \ref{thmTree2Pol}: completion}\label{PfPol}

\subsection{Previously: On polynomial realization}\label{thmPolComm}

Beyond correspondence theorem \ref{Classthm}, theorem \ref{thmTree2Pol} claims realization of \emph{all} planar trees $\mathcal{T}$, as the reduced source/sink connection graphs $\mathcal{C}^*\equiv\mathcal{T}$ of Poincaré compactifications for polynomial ODEs \eqref{ODEP}.
The techniques developed in section \ref{PfthmRat} suggest induction over the degree $d$ of the polynomial $P$, i.e.~over the number $d$ of vertices in the planar trees $\mathcal{T}_d$\,.
The cases $d=2,3$ of figures \ref{fig4}, \ref{fig1}(a) are trivial, because the corresponding trees are unique.
Indeed, corollary \ref{corCountpol} then asserts the counts $A_1=A_2=1$ of $A_{d-1}$\,.

Consider any prescribed planar tree $\mathcal{T}_{d+1}\,,\ d\geq3$, which we have to realize as the reduced connection graph $\mathcal{C}_{d+1}^*$ of a polynomial nonlinearity $f=P_{d+1}$\,.
Like any tree of at least two vertices, $\mathcal{T}_{d+1}$ possesses at least two leaves, i.e. vertices with a single edge attached; see figure \ref{fig2}.
Plucking any leaf from  $\mathcal{T}_{d+1}$ therefore leaves us with a planar tree $\mathcal{T}_d$\,.
By induction hypothesis, let $P_d$ denote any polynomial of degree $d$ which realizes the prescribed planar tree, $\mathcal{C}_d^*\equiv\mathcal{T}_d$\,. 
Scaling \eqref{a=1} makes $P_d$ univariate.
We then have to construct a polynomial $P_{d+1}$ of degree $d+1$, which makes $\mathcal{T}_d$ regrow that plucked leaf, in any of the $2(d-1)$ positions which $\mathcal{T}_d$ has to offer.

Based on $P_{d+1}:=P_d\cdot(1-\eps\exp(\mi\theta)w)$, we have sketched such a construction in section 4.6 of \cite{FiedlerShilnikov}.
Rotating the added large equilibrium $w=e_{d+1}^\theta=\eps^{-1}\exp(-\mi\theta)$ by -$\theta$ produces a sequence $e_{d+1}^\theta$ of Lyapunov centers.
At each center, the stability of $e_{d+1}$ flips between source and sink.
We have sketched how the attachment of the leaf vertex $e_{d+1}$ to the persistent tree $\mathcal{T}_d$ of $P_d$ changes, at each Lyapunov center, to realize all $2(d-1)$ positional possibilities.
See section \ref{Degens} for further discussion.
Although this remains a valid approach, we proceed more directly in the present paper.

The cubic example \eqref{cubic} of section \ref{Ex3} was at the heart of our inductive proof for rational realization theorem \ref{thmPortrait2Rat}, in section \ref{PfthmRat}.
Analogously, the polynomial example \eqref{wdtheta} of section \ref{Exd} is central to the induction which proves polynomial realization theorem \ref{thmTree2Pol}.

\subsection{Proof by outer and inner expansions}\label{OutthmPol}
By induction hypothesis, let $P_d(w)=(w-e_1)\cdot \ldots \cdot (w-e_d)$ realize the tree $\mathcal{C}_d^*\equiv\mathcal{T}_d$\,.
Similarly to section \ref{OutthmRat}, we downscale $P_d$ by some $\delta=\eps\rho>0$ as
\begin{equation}
\label{Pddelta}
P_d^\delta(w):= \delta^d\, P_d(w/\delta) = (w-\delta e_1)\cdot \ldots \cdot (w-\delta e_d)\,.
\end{equation}
Note $P_d^0(w)=w^d$\,.
For $\delta>0$ we have downscaled $P_d$ such that $w(t)$ solves $\dot{w}=P_d(w)$, if and only if $w^\delta(t):=\delta w(\delta^{d-1}t)$ solves $\dot{w}^\delta = P_d^\delta(w^\delta)$.
Our choice for $P_{d+1}$ is 
\begin{equation}
\label{Pd+1}
P_{d+1}(w)=P^\delta(w):=P_d^\delta(w)\cdot (1-\exp(\mi\theta)w)\,,
\end{equation} 
for suitable $\theta$.
The three induction-related ODEs now read
\begin{align}
\label{Pddeltaw}
    \dot{w}=P_d^\delta(w)&\,; \\
\label{Pdeltaw}
    \dot{w}=P^\delta(w)&:= P_d^\delta(w)\cdot (1-\exp(\mi\theta)w)\,,\quad\mathrm{for}\ \delta>0\,; \\
\label{P0w}
    \dot{w}=P^0(w)&:=\ \ \,  w^d \ \ \cdot (1-\exp(\mi\theta)w)\,,\quad\mathrm{for}\ \delta=0\,.
\end{align}
As in \eqref{thetal}, we fix $\theta=\theta_\ell:=\ell\pi/(d-1)$, for some integer $0\leq\ell<2(d-1)$.
Let $|w|=\rho>0$ denote a small circle and, as in the two-scale expansion \eqref{deltaepsrho}, define $\delta:=\eps\rho$\,, for small $\eps>0$.
See figure \ref{fig8} for an illustration of the induction step, in case $d=4$. 

\subsubsection{Inner expansions}\label{InnerPol}
For \emph{inner expansions} we first denote the variable $w$ in \eqref{Pdeltaw} as $w^\delta$, i.e. we study $\dot{w}^\delta=P^\delta(w^\delta)$. 
Backtracking the original transformation $w\mapsto w^\delta(t):=\delta w(\delta^{d-1}t)$, from $w^\delta$ to $w$, we then obtain the equivalent ODE
\begin{equation}
\label{Pdwdelta}
\dot{w}=P_d(w)\cdot (1-\delta\exp(\mi\theta)w)\,.
\end{equation}
The inner expansion is concerned with $|\delta w|=|w^\delta|\leq\rho$. 
Therefore \eqref{Pdwdelta} becomes
\begin{equation}
\label{Pdwrho}
\dot{w} =P_d(w)\cdot (1+\mathcal{O}(\rho))\,,
\end{equation}
which is a small perturbation of $\dot{w}=P_d(w)$\,. 
All equilibria $e_j$ of $P_d$ are preserved in \eqref{Pdwdelta}, as sources and sinks.
Their source/sink heteroclinic orbits are therefore defined by finite time segments of \eqref{Pdwdelta} between the local domains of attraction in negative and positive time, respectively.
In particular, the structurally stable connection tree $\mathcal{T}_d$ of $P_d$ is preserved for \eqref{Pdeltaw}, downsized by $\delta=\eps\rho$ to fit inside the $\rho$-disk, and up to an orientation preserving homeomorphism.
Indeed, bounded $w=\mathcal{O}(1)$ on the tree $\mathcal{T}_d$ are downsized to $|w^\delta|=\delta |w| = \eps\rho\, \mathcal{O}(1)\leq \rho/2 $, for sufficiently small $\eps>0$.
Therefore, we simply denote the equivalent perturbed and downsized connection tree $\mathcal{T}_d$ as $\delta\mathcal{T}_d$ with vertices $\delta e_j$\,.

\subsubsection{Outer expansions}\label{OuterPol}
Consider \emph{outer expansions} next, for $|w|\geq \rho>0$.
Then $|\delta/w|\leq \eps$ in \eqref{Pddelta} implies
\begin{equation}
\label{Pdelta=wd}
P_d^\delta(w)=w^d\cdot (1+\mathcal{O}(\eps))\,.
\end{equation}
We explore tangencies to the $\rho$-circle $|w|=\rho$, for $P^\delta$ in \eqref{Pdeltaw} first, i.e.
\begin{equation}
\label{tang}
0=\big\langle w,P^\delta(w) \big\rangle = 
\big\langle w,\, P_d^\delta(w)(1+\mathcal{O}(\rho)) \big\rangle\ =
\big\langle w,\, w^d(1+\mathcal{O}(\eps)+\mathcal{O}(\rho)) \big\rangle\,.
\end{equation}
See \eqref{Pd+1}.
We turn to polar coordinates $w= \rho \exp(\mi\psi)$, as in \eqref{rw}, \eqref{alphaw}, \eqref{phijtheta}.
Analogously to expansion \eqref{phijl} for $P^0$ in \eqref{P0w}, the angles $\psi=\phi_j$ of tangencies to the $\rho$-circle for $\theta=\theta_\ell$ perturb to
\begin{equation}
\label{phijleps}
\psi=\phi_j(\eps,\rho,\theta_\ell) = \tfrac{1}{2}(\theta_{j-1}+\theta_j)+\mathcal{O}(\eps)+\mathcal{O}(\rho)\,.
\end{equation}
This is the outer expansion of tangencies and their exit/entry intervals.

\begin{figure}[t!]
\centering \includegraphics[width=0.8\textwidth]{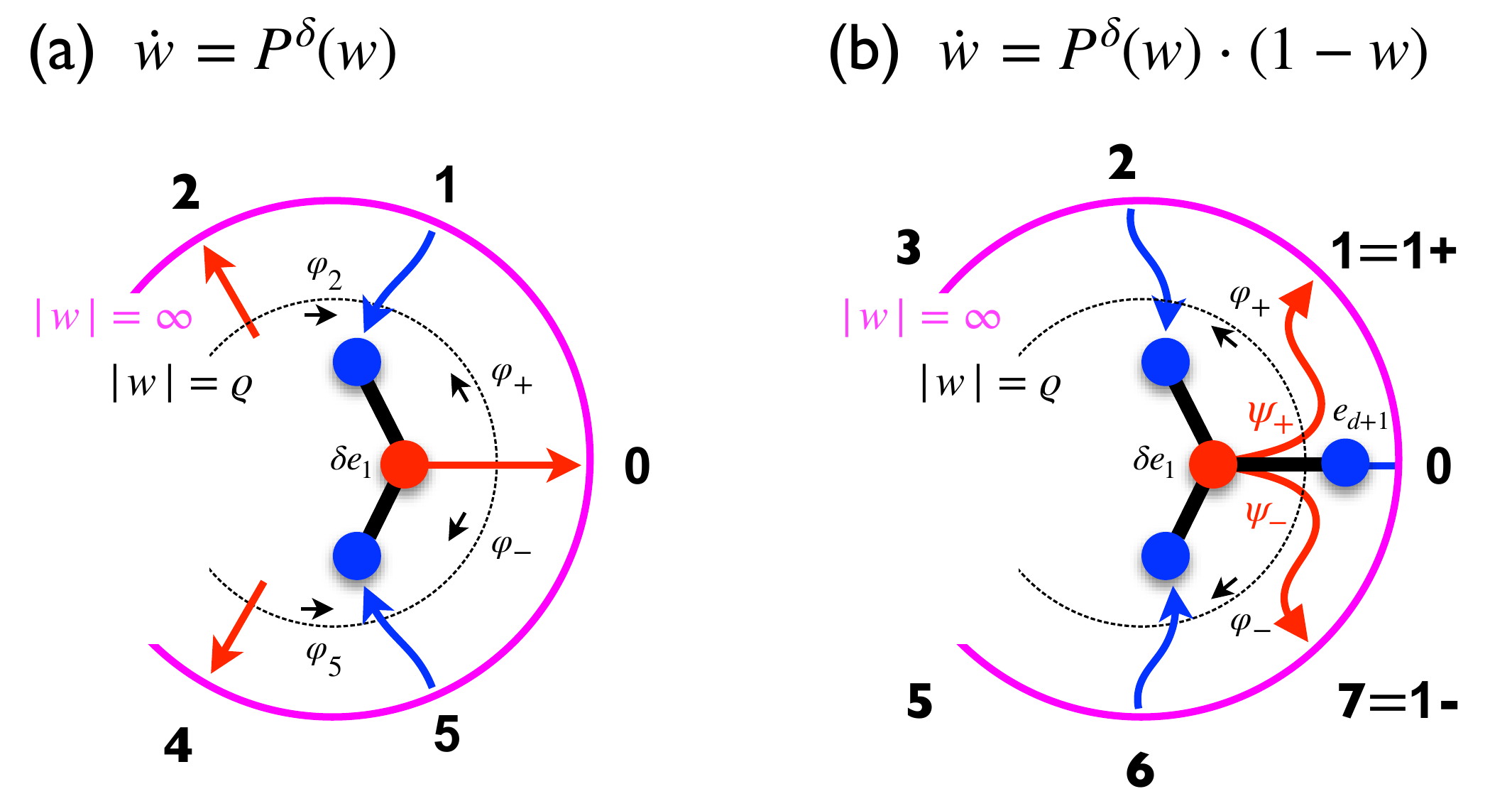}
\caption[Induction over polynomial degree]{\emph{
(a) Phase diagram of ODE \eqref{Pddeltaw} $\dot{w}=P_d^\delta(w)$, for $d=4$.
Tangencies of $P_d^\delta$ to the circle $|w|=\rho$ occur at angles $\phi_j\,,\ j=1,\ldots,2(d-1)=6$.
Indices $+,-$ mark $j=1,\,2(d-1)$.
Intermediate entry and exit intervals are marked by blue and red arrows, respectively.
Trajectories in each entry interval converge to a single blue sink of the downsized tree $\delta\mathcal{T}_d$ in the $\rho$-disk $|w|\leq\rho$, in positive time.
That convergence also applies to the unique blue separatrix entering through the same interval.
Exit intervals show the analogous behavior in negative time, i.e. in red.
Note the red blow-up separatrix $\mathbf{0}$ emanating from the red source vertex $\delta e_1$ in $\delta\mathcal{T}_d$ through the exit interval between $\phi_-=\phi_6$ and $\phi_+=\phi_1$\,.\\
(b) Phase diagram of ODE \eqref{Pdeltaw} for $d=4$ and $\theta=0,\ e_{d+1}=1$, i.e. $\dot{w}=P_d^\delta\cdot(1-w)$.
Up to small perturbations, circle tangencies $\phi_j$\,, entry and exit intervals, and the connection tree $\delta\mathcal{T}_d$ of $P_d^\delta$\,, downsized by $\delta$ to fit in the $\rho$-disk, all persist.
Except for the red separatrix $\mathbf{0}$, the $d-1$ blue and the $d-2$ remaining red separatrices also persist, albeit relabeled, with their continued blue sink and red source targets $\delta e_j$ in $\delta\mathcal{T}_d$\,, respectively.
Only the red separatrix $\mathbf{0}$ gets replaced by the two straddling separatrices $\mathbf{1+}:=\mathbf{1}$ and $\mathbf{1-}:=\mathbf{2d-1}=\mathbf{7}$, which leave the $\rho$-circle at angles $\psi_\pm=\mathcal{O}(\eps)+\mathcal{O}(\rho)$, i.e.~at \,$\phi_-=\phi_6<\psi_-<\psi_+<\phi_1=\phi_+$\,.
The straddling red separatrices $\mathbf{1\pm}$ are heteroclinic from the red source $\delta e_1$ to angles $\alpha_{1}$ and $\alpha_{2d}$ on the purple circle $w=\infty$ of the Poincaré compactification, respectively.
This forces the new separatrix $\mathbf{0}$, now blue, to terminate at the only available sink, $w=e_{d+1}$\,, in the straddled region.
Any trajectory leaving the exit interval between $\psi_-$ and $\psi_+$ has to be heteroclinic from $\delta e_1$ to $e_{d+1}=1$.
In the connection graph $\mathcal{C}_{d+1}^*$ of $P_{d+1}:=P^\delta$, these orbits $\delta e_1\leadsto e_{d+1}$ attach the sink leaf $e_{d+1}$ to the source $\delta e_1$\,, replacing the original red separatrix $\mathbf{0}$ of $\delta e_1$ under $P_d^\delta$\,.\\
General choices $\theta=\theta_\ell$ allow us to substitute any separatrix attached to any source or sink vertex in $\delta\mathcal{T}_d$\,, of (a), by a straddling pair and an attached leaf vertex at $w=\exp(-\mi\theta_\ell)$, in (b).
This proves the polynomial realization theorem \ref{thmTree2Pol}, by induction on the polynomial degree $d$.
}}
\label{fig8}
\end{figure}

As in the example $\delta=0$ of section \ref{Exd}, we also obtain an outer expansion for the straddling separatrices.
Analogously to expansion \eqref{psipml} and the techniques of section \ref{OuterRat}, expansion \eqref{Pdelta=wd} implies
\begin{equation}
\label{psipmleps}
\psi=\psi_\pm(\eps,\rho,\theta_\ell)= -\theta_\ell +\mathcal{O}(\eps)+\mathcal{O}(\rho)\,.
\end{equation}
Both expansions are standard applications of the implicit function theorem.

\subsubsection{Combining inner and outer expansions}\label{InnerOuterPol}
To complete our proof by induction on $d$, we combine inner and outer expansions.
First, consider trajectories of $P_{d+1}:=P^\delta$ in the $2(d-1)$ entry and exit intervals of the $\rho$ circle $|w|=\rho$\,, between the tangencies at $\phi_j$\,.
See \eqref{Pddelta}, \eqref{Pd+1}, \eqref{Pdeltaw}, and \eqref{tang}, \eqref{phijleps}.
In each interval, those trajectories still have to converge to one and the same sink or source vertex of the (downscaled and perturbed) connection tree $\delta\mathcal{T}_d$\,, respectively and in the appropriate time direction.
The inner expansion \eqref{Pdwrho} identifies that target as the unique target of one of the $2(d-1)$ stable and unstable separatrices under unperturbed but downscaled $P_d^\delta$ in \eqref{Pddeltaw}.
The relevant (downscaled) separatrix is defined uniquely by its corresponding crossing interval for $P_d$\,.
See also section \ref{Exd}, and figure \ref{fig8}(a).

By the outer expansion \eqref{psipmleps}, the full system \eqref{Pdeltaw} for $P_{d+1}:=P^\delta$ 
possesses a pair of straddling separatrices, which cross the $\rho$-circle at angles $\psi_\pm=-\theta_\ell+\mathcal{O}(\eps)+\mathcal{O}(\rho)$\,. 
By \eqref{phijleps} that interval is between circle tangencies $\psi=\phi_j$, with $j=-\ell$ and $j=-\ell+1$.
Since $\ell$ can be chosen arbitrarily, let us proceed with our description for $\ell=0$; see figure \ref{fig8}(b).
Odd $\ell$ require reversal of time, and swapped colors.

This identifies the unique source vertex, say $\delta e_1\in\delta E$, of the red straddling separatrices in (perturbed) $\delta\mathcal{T}_d$\,.
The unique blue separatrix between the red straddling separatrices has to converge to $e_{d+1}=1$, as the only available sink in that region.
All other trajectories between the straddling separatrices have to cross the same part of the exit interval and, therefore, are heteroclinic  from $\delta e_1\in\delta\mathcal{T}_d$ to $e_{d+1}=1$.
In other words, $\delta e_1\leadsto e_{d+1}$ attaches $e_{d+1}$ to $\delta e_1$\,, as a leaf of the connection graph $\mathcal{C}_{d+1}^*$\,, in the direction of the original separatrix $\mathbf{0}$ of $P_d$\,.
Since $0\leq \ell < 2(d-1)$ was arbitrary, an appropriate choice of $\ell$ also allows us to attach the vertex $e_{d+1}$, as a leaf, to any vertex $\delta e_j\in\delta\mathcal{T}_d$\,, in the direction of any separatrix of $P_d$ which was attached to $e_j$\,.
In particular, a proper choice of $\ell$ realizes the prescribed graph, $\mathcal{C}_{d+1}^*\equiv\mathcal{T}_{d+1}$\,, for $P_{d+1}:=P^\delta$ with $\delta:=\eps\rho$\, and small $\eps,\rho>0$.
This proves the polynomial realization theorem \ref{thmTree2Pol}, by induction on the polynomial degree $d$. \hfill $\bowtie$

\section{Proof of anti-polynomial theorem \ref{nctree=antipol}}\label{PfAntipol}

With the remarks of section \ref{ResAntipol}, it only remains to prove the realization part \emph{(ii)} of anti-polynomial theorem \ref{nctree=antipol}.

\emph{Proof of theorem} \ref{nctree=antipol}\emph{(ii).} \ 
We sketch a proof of realization claim \emph{(ii)}.
We write $d$ for $d'$, in this proof.
We omit details which are analogous to the techniques developed in section \ref{PfPol}.

Replacing $P$ there by univariate $Q$, let $Q_d:=-\mi Q$ in \eqref{Qbar} realize the dual pair $\mathcal{C}_d^\pm\equiv\mathcal{T}_d^\pm$ of nc-trees with $d$ edges, as phase portraits.
The imaginary prefactor $-\mi$ of the univariate polynomial $Q$, can be chosen without loss of generality; compare \eqref{a=1}.
Again, we downscale $Q_d$ by some $\delta=\eps\rho>0$ as
\begin{equation}
\label{Qddelta}
Q_d^\delta(w):= \delta^d\, Q_d(w/\delta)\,.
\end{equation}
Note $Q_d^0(w)=-\mi w^d$\,.
Our choice for $Q_{d+1}$ remains
\begin{equation}
\label{Qd+1}
Q_{d+1}(w)=Q^\delta(w):=Q_d^\delta(w)\cdot (1-w)\,.
\end{equation} 
Analogously to \eqref{Pddeltaw}--\eqref{P0w}, our three induction-related anti-polynomial ODEs now read
\begin{align}
\label{Qddeltaw}
    \dot{w}=\overline{Q_d^\delta(w)}&\,; \\[1mm]
\label{Qdeltaw}
    \dot{w}=\overline{Q^\delta(w)}&:= \overline{Q_d^\delta(w)}\cdot (1-\overline{w})\,,\quad\mathrm{for}\ \delta>0\,; \\[1mm]
\label{Q0w}
    \dot{w}=\overline{Q^0(w)}&:=\  \,  \mi\,\overline{w}^d \ \ \cdot (1-\overline{w})\,,\quad\mathrm{for}\ \delta=0\,.
\end{align}
We study leaf attachment in the tree portraits $\mathcal{C}_0^\pm$ of anti-polynomial ODE \eqref{Q0w}, first.
We will then complete the proof with a description of the resulting leaf attachments which occur when passing from the nc-trees $\mathcal{C}_d^\pm\equiv\mathcal{T}_d^\pm$ of \eqref{Qddeltaw} to the nc-trees $\mathcal{C}_{d+1}^\pm$ of \eqref{Qdeltaw}, which will realize any prescribed nc-pair $\mathcal{T}_{d+1}^\pm$ \,.

\begin{figure}[t]
\centering \includegraphics[width=0.5\textwidth]{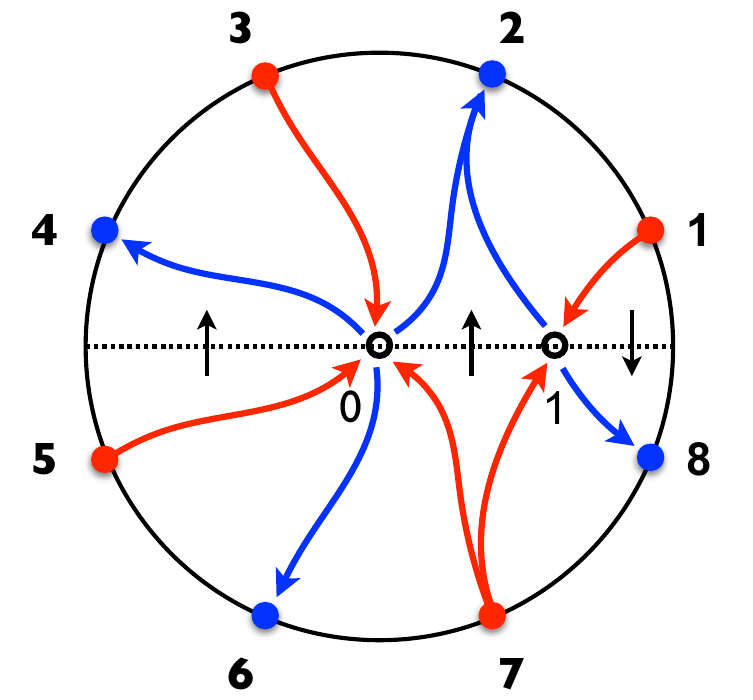}
\caption[Induction over anti-polynomial degree]{\emph{
Illustration of the induction step of theorem \emph{\ref{nctree=antipol}(ii)} on anti-polynomial ODE \eqref{Qbar}, from degree $d=2$ to degree $d+1=3$.
Here $\dot{w}=\mi\,\overline{w}^d(1-\overline{w})$ of degree $d+1$, see \eqref{Q0w}.
Note how complex conjugation $w\mapsto \overline{w}$ reverses time and swaps the colors of separatrices and sources/sinks.\\
Poincaré compactification identifies $d+2$ source angles $\alpha_k=(k-\tfrac{1}{2})\pi/(d+2)$ at $w=\infty$, for odd $1\leq k\leq 2(d+2)$, which we indicate as $\mathbf{k}$.
Even $\mathbf{k}$ define $d+2$ sinks, which alternate with the sources.\\
The equilibrium $w=0$ is of multiplicity $d$ and of hyperbolic type.
On any circle of small radius, $|w|=\rho$, it produces $2(d+1)$ exterior tangencies which separate $2(d+1)$ alternating entry and exit intervals of \eqref{Q0w}.
Each such interval is crossed by a unique separatrix of $w=0$.
The $2(d+1)$ separatrices connect, one-to-one, with the $2(d+1)$ alternating
 sources and sinks at $\alpha_k$\,, for $2\leq k\leq 2d+3$.\\
At the nondegenerate hyperbolic saddle $w=1$, we obtain two red stable separatrices and two blue unstable separatrices: blue outgoing towards the sinks  $\mathbf{2},\ \mathbf{2d+4}=\mathbf{8}$, and red incoming from the sources $\mathbf{1},\ \mathbf{2d+3}=\mathbf{7}$.
In the red portrait of \eqref{Q0w}, this attaches the source $\mathbf{1}$ as a leaf to the vertex $\mathbf{2d+3}$ of the tree of $\dot{w}=\mi\overline{w}^d$.
By time reversibility, the blue portrait of \eqref{Q0w} attaches the sink $\mathbf{2d+4}$ as a leaf to the vertex $\mathbf{2}$ of $\mi\overline{w}^d$, symmetrically.
Rotating the pole $w=1$ appropriately, such leaf attachments can be performed at any other vertex on the Poincaré circle $w=\infty$.
}}
\label{fig9}
\end{figure}

Let $|w|=\rho>0$ denote a small circle and, as in the two-scale expansion \eqref{deltaepsrho}, define the downsizing scaling factor as $\delta:=\eps\rho$\,, for small $\eps$.
See figure \ref{fig9} for an illustration of the Poincaré compactification of ODE \eqref{Q0w}, in case $d=2$.

The equilibrium $w=0$, of degenerate multiplicity $d$, is of hyperbolic type.
In polar coordinates $w=r\exp(\mi\psi)$, and compared to example \ref{Exd}, \eqref{phijtheta}, the $2(d+1)$ exterior tangencies of $\overline{Q^0(w)}$ to the circle $|w|=r=\rho$ now occur at angles \begin{equation}
\label{phijQ0}
\psi=\phi_j(\rho)=\phi_j^*+\mathcal{O}(\rho)\,,\quad \phi_j^*=j\pi/(d+1)\,,
\end{equation}
for $0\leq j< 2(d+1)$.
As before, they define exit and entry intervals on the $\rho$-circle.
The $2(d+1)$ separatrices of $w=0$ cross the circle at the intermediate angles
\begin{equation}
\label{psijQ0}
\psi=\psi_j:=(j-\tfrac{1}{2})\pi/(d+1)+\mathcal{O}(\rho)\,,
\end{equation}
blue outgoing for odd $1\leq j\leq 2(d+1)$, and incoming red for even $j$.

The other finite equilibrium, $w=1$, is a nondegenerate saddle.
The two red incoming separatrices arrive at asymptotic angles $\tfrac{1}{4}\pi$ and $\tfrac{5}{4}\pi$, locally.
Perpendicularly, blue outgoing separatrices emanate from asymptotic angles $\tfrac{3}{4}\pi$ and $\tfrac{7}{4}\pi$.

On the Poincaré circle at $w=\infty$, the equilibria of \eqref{Q0w} are now located at angles
\begin{equation}
\label{alphaQ0}
\psi=\alpha_k:=(k-\tfrac{1}{2})\pi/(d+2)+\mathcal{O}(\rho)\,,
\end{equation}
due to the imaginary prefactor $\mi$.
We obtain $d+2$ red sources $\mathbf{k}$ for odd $1\leq k\leq 2(d+2)$, which alternate with $d+2$ blue sinks for even $\mathbf{k}$.
See figure \ref{fig9}, for $d=2$.

Next, we match the $2(d+3)$ red and blue separatrices of $w=0$ and $w=1$ with the $2(d+2)$ alternating sources and sinks $\mathbf{k}$ at $w=\infty$.
Saddle-saddle loops or heteroclinic orbits are excluded by the simultaneous gradient and Hamiltonian structure.
Complex conjugation $w\mapsto \overline{w}$ reverses time and swaps the colors of separatrices and sources/sinks.
In particular, separatrices cannot cross the real axis.
It is therefore sufficient to study the upper half plane $v=\Im w>0$.
We assume $d$ is even; the arguments for odd $d$ are analogous.

At $w=\infty$ and for even $d$, we observe equal numbers $(d+2)/2$ of upper red sources and upper blue sinks.
The equilibrium $w=0$ contributes $d/2$ upper red separatrices, and $w=0$ contributes another one.
In particular, the numbers of upper red separatrices and red sources coincide.
Indeed, each red source must figure as the $\boldsymbol{\alpha}$-limit set of at least one red separatrix, to separate the domains of attraction of its neighboring sinks. 
This matches the red separatrices and the red sources, one-to-one.
Specifically, $\mathbf{1}\leadsto 1$, and $\mathbf{k}\leadsto 0$, for odd $1< k\leq d+2$.
Since the blue upper separatrices are confined within the red sectors, this also proves $0 \leadsto \mathbf{k}$, for each of the $(d+2)/2$ emanating tangent angles $\phi_j$\,, at odd $1\leq j\leq d+1$ and for even $k=j+1$.
Moreover, we observe an additional blue separatrix $1\leadsto \mathbf{2}$.
Reversibility reflects those results to the lower half plane $v=\Im w<0$, with reversed colors.
Therefore, the red separatrices of the saddle $w=1$ to $\mathbf{1}$ and $\mathbf{2d+3}$, as an edge $[\mathbf{1},\mathbf{2d+3}]$, attach the red source $\mathbf{1}$, as a leaf, to the source vertex $\mathbf{2d+3}$ of the red tree inherited from $\dot{w}=\mi\,\overline{w}^d$.
Reflection attaches the blue leaf $[\mathbf{2},\mathbf{2(d+2)}]$ to the inherited blue tree.
This establishes the leaf attachments of figure \ref{fig9}, as parts of the tree portraits $\mathcal{C}_0^\pm$ for the anti-polynomial ODE \eqref{Q0w}.

More precisely, the attached red leaf is a \emph{short leaf}: the leaf connects two adjacent red sources on the Poincaré circle, which are separated by only a single sink.
And colors may be swapped.
By an elementary induction argument, any nc-tree contains a short leaf.

To perform the induction step and complete the proof, we will substitute the downsized nc-tree $\delta\mathcal{C}_d^-$ of \eqref{Qddelta} for $w=0$, in the above discussion of \eqref{Q0w}.
Passing to the analogue \eqref{Qdeltaw} of \eqref{Q0w} will then realize the desired nc-tree $\mathcal{C}_{d+1}^-\equiv\mathcal{T}_{d+1}$\, with short leaves attached in any desired position. 

We have to start from an arbitrarily given nc-tree, say $\mathcal{T}_{d+1}$\,.
Rotating and reflecting the tree into the proper position, we may assume $[\mathbf{2},\mathbf{2(d+2)}]$ to be a short leaf of $\mathcal{T}_{d+1}$\,. 
By induction hypothesis, let $Q$ in \eqref{Qbar} or, equivalently, downsized $Q^\delta$ in \eqref{Qddelta} realize the nc-tree $\mathcal{C}_d^-\equiv \mathcal{T}_d$\,, defined as $\mathcal{T}_{d+1}$ with the short leaf removed. 
For scaling $\delta=\eps\rho$ consider the $\rho$-circle $|w|=\rho$.
We then obtain expansions for the tangency angles $\phi_j=\phi_j(\eps,\rho)$ and the separatrix angles $\psi_j=\psi_j(\eps,\rho)$, as in \eqref{phijQ0} and \eqref{psijQ0}, but with the remainder terms $\mathcal{O}(\rho)$ replaced by $\mathcal{O}(\rho)+\mathcal{O}(\eps)$.
Outer expansions, for $|w|\geq\rho$, analogously to \eqref{Pdelta=wd}, show
\begin{equation}
\label{Qdelta=wd}
Q^\delta(w)=Q^0(w)(1+\mathcal{O}(\eps))\,.
\end{equation}
We may therefore repeat the matching arguments for separatrices and sinks/sources of the tree portraits $\mathcal{C}_0^\pm$ for $Q^0$ in \eqref{Q0w}, to perform the same matchings for the perturbed separatrices of $\mathcal{C}_\delta^\pm$ for $Q^\delta$ in \eqref{Qdeltaw}.
In \eqref{Q0w}, this substitutes the nc-tree $\mathcal{C}_d^-$ for $w=0$.
Therefore the portrait $\mathcal{C}_{d+1}^-$ of $Q_{d+1}=Q^\delta$ consists of $\mathcal{C}_d^-\equiv\mathcal{T}_d$, with the blue short leaf $[\mathbf{2},\mathbf{2(d+2)}]$ attached at sink $\mathbf{2}$.
This proves $\mathcal{C}_{d+1}^-\equiv\mathcal{T}_{d+1}$\,, completing the induction on the degree $d$, and proves part \emph{(ii)} of theorem \ref{nctree=antipol}. 
\hfill $\bowtie$

Of course, we may reinterpret the purely polynomial results $Q=1$ of theorems \ref{Classthm}--\ref{corCountpol} as an exploration of antiholomorphic vector fields $\dot{w} =\overline{f(w)}=1/\overline{P(w)}$.
Simple zeros $e_j$ of $P$, however, introduce logarithmic potentials $F$, and periods $\mathcal{P}$.
To define the potential $F$ properly, we have to pass to the universal cover of the multi-punctured Riemann sphere -- a problem which does not arise at higher multiplicities.

\section{Discussion}\label{Dis}

We revisit some of our results, and sketch a few possible extensions and open problems.
Specifically, we address \emph{degeneracies} of rational vector fields $\dot{w}=f=P/Q$ of ODE \eqref{ODE} and hybrid holomorphic/anti-holomorphic variants in section \ref{Degens}.
Extensions to \emph{Riemann surfaces} are wide open; see section \ref{Riems}.
Planar \emph{Sturm attractors}, the class of global attractors for scalar parabolic PDEs on the interval with only hyperbolic equilibria $e_j$ of Morse indices $i(e_j)\leq 2$, are all realized in our rational setting; see section \ref{Sturm}.
We commemorate the Yamaguti legacy \cite{YamagutiRomantic}, in section \ref{Romantic}, and conclude on a philosophical note.

\subsection{Degeneracies}\label{Degens}

Complex rational vector fields $f=P/Q$ were assumed nondegenerate; see definition \ref{defnondeg}.
Nondegeneracy led to the gradient-like Morse structure of the regularized flow; see corollary \ref{corMorse} and proposition \ref{propgrad}.
We have established genericity of nondegeneracy, in terms of the coefficients of $P,Q$; see theorem \ref{thmRatGen}.
Violations of nondegeneracy, however, already occur in families $f=f(\lambda,\cdot)$ with a scalar real parameter $\lambda$.
Global descriptions of their unfoldings are a widely open field, in the complex context.
See \cite{Rousseaud=4} for interesting bifurcations, including degree $d=d'+2=4$ and a 2-parameter unfolding.

Multiple zeros or poles violate nondegeneracy, generating elliptic and hyperbolic sectors.
Even the purely polynomial case $\dot{w}=P$ of section \ref{PfPol} provided examples.
We recall how polynomials $P$ of degrees $d\geq 4$ introduce a degenerate pole of multiplicity $d-2\geq 2$ at $w=\infty$.
For $d=2$ recall the quadratic Masuda paradigm of section \ref{Ex2} without poles and with a regular flow at $w=\infty$.
This subsumes the linear case $d=1$.
Figure \ref{fig1}(a) and the cubic rational example of figure \ref{fig5} in section \ref{Ex3} illustrate $d=3$.
For the example $P=w^d(1-w)$ of degree $d+1\geq 3$ with a degenerate equilibrium $w=0$ of multiplicity $d$ see section \ref{Exd} and figure \ref{fig6}.

A general analysis of Lyapunov centers, saddle-saddle connections, multiple zeros, and poles for the rational case $f=P/Q$ is beyond our present scope.
See \cite{DiPol} for a complete characterization of the global dynamics in the polynomial case $f=P$, and \cite{DiRat} for the general rational case $f=P/Q$.
The characterization is presented via nine necessary properties of global flows.
Avoiding induction, meromorphic realization by $f=P/Q$ is then sketched by the construction of a complex structure, from these properties, and subsequent grand uniformization.
The procedure seems to go back to an unpublished manuscript by Adrien Douady et al \cite{Rousseaud=4}.
Genericity and the antiholomorphic case $P=1$ are not addressed specifically.
For some steps towards flows in the presence of essential singularities of $f$ see \cite{MuciEss,LebEss} and the references there.

Lyapunov centers $f(e_j)=0\neq f'(e_j)=\mi\omega$ are surrounded by their local disk of nested periodic orbits.
More generally, real-time periodic orbits occur as iso-periodic local foliations; see section \ref{ODEper}.
They imply degenerate sums \eqref{sumetaJ0} for the real parts of the residues $\eta_j=1/f'(e_j)$ at certain interior zeros $j\in J$ of $P$.
In general, the iso-periodic families need not limit onto Lyapunov centers.
In the presence of poles, they may also limit onto homoclinic saddle loops or heteroclinic cycles \eqref{hetcycle} among several saddles.

We have encountered polynomial examples in sections \ref{Exd} and \ref{PfPol}.
See \eqref{wdtheta}, figure \ref{fig6} and \eqref{Pdeltaw}, figure \ref{fig8}.
Discrete values $\theta=\theta_\ell=\ell\pi/(d-1)\in\mathbb{S}^1=\mathbb{R}/2\pi$ of the angle parameter $\theta$ in $P_{d+1}$ have attached the new hyperbolic equilibrium vertex $e_{d+1}=\exp(-\mi\theta)$ of $\mathcal{T}_{d+1}$ as a leaf to the downsized tree $\delta\mathcal{T}_d$\,.
We recall that $e_{d+1}$ is a sink, for even $\ell$, and a source for odd $\ell$.
Cycling through $\theta\in\mathbb{S}^1$, continuously, rotates $e_{d+1}$ clockwise, and we encounter a sequence of $2(d-1)$ Lyapunov centers.

Consider the simplest case $\dot{w}=P(w):=w^d(1-\exp(\mi\theta)w)$ of section \ref{Exd}.
Then the eigenvalue
\begin{equation}
\label{P'}
P'(e_{d+1})=-\exp(-(d-1)\mi\theta)
\end{equation}
co-rotates clockwise.
We encounter Lyapunov centers $e_{d+1}$ with eigenvalues $P'=(-1)^{\ell+1}\mi$ at the intermediate parameters $\theta=\theta_\ell^*=\tfrac{1}{2}(\theta_{\ell-1}+\theta_\ell)=(\ell-\tfrac{1}{2})\pi/(d-1)$.
At $\theta=\theta_0^*$ with $P'=-\mi$, for example, the family of periodic orbits generated by the Lyapunov center surrounds $e_{d+1}=\exp(\tfrac{1}{2}\pi\mi/(d-1))$ clockwise.
The family is bounded by a homoclinic loop to $w=\infty$, alias a heteroclinic orbit which merges the separatrices $\mathbf{0}$ and $\mathbf{1}$ of the Poincaré compactification into a saddle-saddle connection; compare figure \ref{fig6}(b).
See also the Lyapunov centers $e_1$ and $e_3$ in figure \ref{fig1}(c), (d).
Further increase of $\theta$ reconnects $\mathbf{0}$ to $w=0$.
The sink equilibrium $e_{d+1}$ turns into a source with attached separatrix $\mathbf{1}$, instead, which becomes detached from $w=0$.
Alternatingly, in backward and forward time, this scenario repeats, as $\theta$ cycles through the remaining Lyapunov parameters $\theta_\ell^*$\,.

An analogous description applies to the induction scenario of \eqref{Pdeltaw} and figure \ref{fig8}(b).
We then cycle through the straddle orbit configurations associated to $\theta=\theta_\ell$\,, as described in section \ref{PfPol}.
Again, the rearrangement of connectivity occurs by saddle-saddle loops to $w=\infty$ at the $2(d-1)$ Lyapunov parameters $\theta=\theta_\ell^*+\mathcal{O}(\delta)$.
Structural stability holds between those transitions.
We have sketched this alternative bifurcation approach to leaf attachment in section 4.6 of \cite{FiedlerShilnikov}.

In the rational case $\dot{w}=P/Q$, Lyapunov transitions and saddle loops or cycles define analogous transitions between the dual portrait pairs $\mathcal{C}^\pm$ and their associated Morse connection graphs $\mathcal{C}$.
For any fixed degree $d=d'+2$, we may therefore view any degeneracy of real codimension 1 as an edge between the nondegenerate graphs, as vertices.
The resulting ``metagraphs'' of heteroclinic connectivity seem wide open, at this stage.
For $d=4$ and codimension up to 2 see \cite{Rousseaud=4}.

Returning to the hybrid holomorphic/anti-holomorphic viewpoint, of course, our whole present study can be viewed as a study of products of polynomials, $\dot{w}=P\overline{Q}$\,; see the original regularization \eqref{ODEreg}.
In general, however, it remains to explore and include the effects of multiple poles and multiple zeros.
Swapping the roles of $P$ and $Q$ then also subsumes the fully ``anti-rational'' case $\dot{w}=\overline{P/Q}$ as $\dot{w}=Q\overline{P}$.
The global dynamics of meromorphic, anti-meromorphic, and hybrid vector fields $f$ on given Riemann surfaces, even of compact type, seems fairly open, as do higher complex dimensions.
For a recent survey of at least well-posedness for some hybrid quadratic Schrödinger PDEs including the variant $\overline{w}^2$ of the purely quadratic Fujita nonlinearity $w^2$, see \cite{LiuR} and the references there.

\subsection{Riemann surfaces}\label{Riems}

In our context of complex scalar ODEs, Riemann surfaces $\mathcal{R}$ arise in (at least) two different ways: as ODEs $\dot{w}=f(w)$ \emph{on} Riemann surfaces $w\in\mathcal{R}$, and as Riemann surfaces \emph{of} solutions $t\mapsto w(t)$.
See \cite{Forster, Jost, Lamotke} for a general background.

See sections 3 and 6 of \cite{FiedlerShilnikov} for the solution aspect of the polynomial case $f=P$.
The notion \eqref{calP} of periods $\mathcal{P}$ generated by the residues $\eta_j=1/f'(e_j)$ plays a central role there.
Solutions of complex entire vector fields are another example.
Consider solutions $w(t)$ of $\dot{w}=f(w)=\exp(w)$, for example.
Separation of variables \eqref{sov} shows $w(t)=-\log(-t+\exp(-w_0))\ \mathrm{mod}\ 2\pi\mi\mathbb{Z}$.
At $w=1/z=\infty$, we run into an essential singularity of the vector field $\dot{z}=-z^2 \exp(1/z)$.
In particular, Poincaré compactification fails.

Concerning $\dot{w}=f(w)$ \emph{on} Riemann surfaces $w\in\mathcal{R}$, we have only studied rational $f=P/Q$ here, i.e.~the case of the Riemann sphere $f:\widehat{\mathbb{S}}\rightarrow \widehat{\mathbb{S}}$.
Meromorphic vector fields $f:\mathcal{R}\rightarrow\widehat{\mathbb{S}}$ seem like a natural generalization.
Second order ODEs
\begin{equation}
\label{ODEw"}
\ddot{w}+\tfrac{1}{2}f(w)=0
\end{equation}
for polynomial nonlinearities $f=P$ are an example.
Indeed, preservation of complex energy $c$ reduces this to scalar ODEs
\begin{equation}
\label{ODEw'2}
\dot{w}^2+F(w)=c\,,
\end{equation}
where $F':=P$ defines any primitive function $F$ of $f=P$.
In other words, we obtain an ODE on compact Riemann surfaces $w\in\mathcal{R}$, alias complex algebraic curves, here of hyperelliptic type $z^2+F(w)=c$\,.
The cubic cases $F=-4w^3+g_2\,,\ c=-g_3$ lead to the classical meromorphic \emph{Weierstrass elliptic functions} 
\begin{equation}
\label{wp}
w(t)=\wp(t;g_2,g_3)\,.
\end{equation}
These are doubly periodic solutions of the ODE \eqref{ODEw'2} on 2-tori $w\in\mathcal{R}\cong\mathbb{T}_\tau^2=\mathbb{C}/\Lambda_\tau$ of certain lattices $\Lambda_\tau:=\mathbb{Z}+\tau\mathbb{Z},\ \Im\tau>0$.
Explicitly, 
\begin{equation}
\label{T2R}
(\wp,\wp'):\ \mathbb{T}_\tau^2\rightarrow\mathcal{R}
\end{equation} 
provides a biholomorphic equivalence, which pulls the flow \eqref{ODEw'2} back to $\mathbb{T}_\tau^2$ as the parallel flow $\dot{t}=1$.
See section 3 in \cite{FiedlerFila}, for example, and the detailed references there.

More generally, for example, it may be interesting to explore $\dot w=f(w)$ for meromorphic vector fields $f$ on the $2$-torus $w\in\mathbb{T}_\tau^2$.
The elliptic function $f$ then turns out to be rational in the Weierstrass function $\wp$ and its derivative $\wp'$.
Generalizations of realization theorem \ref{thmPortrait2Rat} point at real-time recurrence like blow-up involving Cherry flows on $\mathbb{T}_\tau^2$ \cite{MuciRiem}.

These few remarks on even the simplest case of \eqref{ODEw'2} for cubic $F$ may be sufficient to indicate the intriguing complexities awaiting us, in real time, for meromorphic vector fields $f$ on hyperelliptic or, more generally, on prescribed compact Riemann surfaces $\mathcal{R}$ -- not to speak of holomorphic or meromorphic $f$ in the non-compact hyperbolic cases.

\subsection{Planar PDE Sturm attractors}\label{Sturm}

Planar \emph{Sturm attractors} $\mathcal{A}$ of parabolic PDEs in one space dimension provide a nontrivial class of Morse systems beyond \eqref{PDEw}.
In real time $t\geq0$, we consider general dissipative semilinear parabolic equations
\begin{equation}
\label{w_t}
w_t=w_{xx}+f(x,w,w_x)
\end{equation}
on the unit interval $0<x<1$, with Neumann boundary conditions $u_x=0$ at the endpoints $x=0,1$.
For an early discussion of some \emph{systems} of such reaction-advection-diffusion models, see \cite{YamagutiPDE}.
For some recent work on the \emph{scalar} case, and background references on the terminology used, see \cite{firoSFB,firoFusco}.
For a recent survey on \eqref{PDEw} under periodic boundary conditions see \cite{firoS1}.

The scalar PDE \eqref{w_t} is gradient-like \cite{Zelenyak, Matano, Lappicy}.
The nonlinearity $f\in C^2$ is assumed \emph{dissipative}, i.e., a large ball in a suitable Sobolev space $w(t,\cdot)\in X\subseteq C^1$ of spatial profiles attracts all solutions $w(t,x)$ in forward time.
Similarly to example \ref{cubic}, we think of $w=\infty$ as a source ``equilibrium'', in the dissipative case.

Then there exists a \emph{global attractor} $\mathcal{A}\subset X$ which consists of all eternal, i.e.~global, solutions $w=w(t,x),\ t\in\mathbb{R}$, which remain uniformly bounded in, both, positive and negative real times $t$.
The global attractor $\mathcal{A}$ consists of all equilibria $v\in E$, and the heteroclinic orbits among them.
We also assume hyperbolicity of all equilibria $v$, i.e.~of the time independent solutions $w(t,x)=v(x)$ of \eqref{w_t}.
Then the set $E$ of equilibria is finite.
As before, let $i(v)$ denote the Morse index of $v\in E$, alias the unstable dimension, i.e.~the  dimension of the local unstable manifold.
Due to a spatial nodal property of Sturm type in the spirit of \cite{MatanoLap}, the gradient-like flow then becomes Morse, with automatic transversality of stable and unstable manifolds.
Note $\dim\mathcal{A}=\max i(v)$.
The global attractor embeds into a dissipative flow on Euclidean space $\mathbb{R}^n$ of the same dimension; see \cite{Jolly, Brunovsky}.
The celebrated $\lambda$-Lemma implies transitivity of the heteroclinic relation $\leadsto$ on $E$.
Again due to the nodal property, a \emph{cascading principle} asserts, conversely, that any heteroclinic orbit is generated by a sequence of heteroclinic orbits between equilibria of decreasing, but adjacent, Morse indices \cite{Cascading}.
Moreover, the heteroclinic orbit between any pair of equilibria at adjacent Morse levels is unique.
In particular, the \emph{connection graph} $\mathcal{C}_f$ of heteroclinic connectivity between adjacent Morse levels determines all heteroclinicity.

In the planar case $\dim\mathcal{A}\leq 2$, of course, we obtain a Morse flow on $\mathbb{S}^2$, as in section \ref{S2Morse}.
In particular, rational realization theorem \ref{thmPortrait2Rat} implies the following corollary.

\begin{cor}
Consider the connection graph $\mathcal{C}_f$ under \eqref{w_t}, for any planar PDE global attractor $\mathcal{A}$ of Sturm type, with hyperbolic equilibria.
Then $\mathcal{C}_f\subsetneq\mathcal{C}=\mathcal{C}^+\cup \mathcal{C}^-$ arises in the connection graph $\mathcal{C}$ for the regularization \eqref{ODEreg} of a nondegenerate rational ODE \eqref{ODE}, with a source at $w=\infty$.
To obtain $\mathcal{C}_f$ from $\mathcal{C}$, only the source at infinity and its unstable separatrices have to be removed.
In particular, the PDE semiflow on $\mathcal{A}$ is $C^0$ orbit equivalent to the correspondingly restricted, regularized flow of the rational model in the Riemann sphere $\widehat{\mathbb{C}}$.
\end{cor} 

Unlike rational vector fields, however, the Sturm case $\mathcal{C}_f$ does not realize \emph{all} connection graphs, with just a source removed.
A characterization of the planar connection graphs $\mathcal{C}_f$ for the hyperbolic Sturm case \eqref{w_t} has been established in \cite{Sturm2, Sturm1, Sturm3}.
In higher dimensions, the graph defines a regular cell complex (except for the cell at infinity).
See also the decomposition as a signed Thom-Smale complex in \cite{ThomSmale}.
This excludes planar portraits $\mathcal{C}^\pm$ with loops, for example, as well as cell boundaries with an extra spike-edge pointing inward.
Moreover, the connection graph possesses a bipolar orientation, in the sense of \cite{Rosenstiehl}.
For Hamiltonian nonlinearities $f=f(w)$, further restrictions apply \cite{firoS1}.
A general geometric characterization of higher-dimensional Sturm global attractors $\mathcal{A}$ remains wide open.

\subsection{And what does it all mean?}\label{Romantic}
Yamaguti's views were much broader than comments like the above could fathom.
In \cite{YamagutiRomantic}, he offers us a deeply philosophical perspective on past and present \hbox{mathematics --} with a glimpse of the future.
He identifies a \emph{``classical''} direction with protagonists like David Hilbert, John von Neumann, and others. 
Mathematics, in that sense, is a purely axiomatic construct, consistent and complete, supposedly.
\emph{Consistency} requires absence of self-contradictions.
\emph{Completeness} requires that any statement can be asserted to be either true, or else false, in finitely many steps.
Kurt Gödel defeated that program, for axiomatic systems containing at least the axioms of the positive integers \cite{Godel}.
Yamaguti, in contrast, criticizes how the purely abstract, axiomatic viewpoint rejects any \emph{meaning} of mathematics, outside its axiomatic realm.
He fascinatingly contrasts this ``classical'' attitude with a \emph{``romantic''} direction, exemplified by Henri Poincaré, Norbert Wiener, and others, where mathematical \emph{meaning} transcends the confines of a pure mind game. 
Most importantly from an applied perspective, ``romanticism'' opens the door to the world of Science and, through such filters, the real world itself.
Yamaguti explores this connection with techniques of semiotics and linguistics, inspired by René Thom and others.

We add a few remarks, following up on \cite{Truth, Wunder}.
At least since Plato, reflections concerning that amazing ``real world'', which includes ourselves as such a small part, have remained a primary motivation of philosophy.
An example is the full title \emph{``Philosophiae naturalis principia mathematica''} of Newton's famed 1687 \emph{Principia} \cite{Newton}, at a juncture to modern Science and arguably the birth place of nonlinear dynamics.
See \cite{dynameis} for the long and convoluted history of \begin{otherlanguage}{greek}δυνάμεις\end{otherlanguage}\, (dynámeis), from the Aristotelian concept towards present day virtual realities.
In \cite{Thomas} Q16 De veritate, Thomas Aquinas highlights the same precarious juncture, when he defines (absolute) \emph{Truth} as
\begin{center}
\textbf{ad{\ae}quatio rei et intellectus.}
\end{center}
The adequacy or correspondence -- but not equality -- is between the ``thing'', \emph{``res"}, and its representations in the \emph{``intellect''}.
A century before, Averroes (aka Ibn Rushd) had required consistency: \emph{``Truth is not self-contradictory''} \cite{Ave}.

In Chinese and Japanese logographic script, I am told, three characters become involved:
\begin{enumerate}[(i)]
  \item \begin{CJK}{UTF8}{min}真\end{CJK} (zh\={e}n, shin), i.e.~``true/truth", absolute, universal, crystalline, without ambiguity;
  \item \begin{CJK}{UTF8}{min}実\end{CJK} (shí, jitsu), i.e.~``reality", also a (hard) ``fruit", concrete and specific;
  \item \begin{CJK}{UTF8}{min}理\end{CJK} (lĭ, ri), i.e.~``reason'', also ``order'' or ``principle'', as in Law of Nature.
\end{enumerate}
For example, they combine as \begin{CJK}{UTF8}{min}真理\end{CJK} (zh\={e}n-lĭ, shin-ri), for underlying absolute Truth in an abstract sense, or as
\begin{CJK}{UTF8}{min}真実\end{CJK} (zh\={e}n-shí, shin-jitsu), for apparent truth as opposed to falsehood and lie, e.g.~in legal context.
With ``truth'' \begin{CJK}{UTF8}{min}真\end{CJK} in hand, nothing else matters.
``Reason'' \begin{CJK}{UTF8}{min}理\end{CJK} permeates \cite{YamagutiRomantic}, in various combinations like axiom, theorem, sciences, psychology, logics, and understanding.
It represents a complex concept in Neo-Confucian philosophy, although of much older origin.
In the work of influential Ming philosopher-poet-general \begin{CJK}{UTF8}{min}王陽明\end{CJK}(Wang Yangming),
the categories ``reality'' \begin{CJK}{UTF8}{min}実\end{CJK} and even absolute ``truth'' \begin{CJK}{UTF8}{min}真\end{CJK} emanate from the light of innate reason \begin{CJK}{UTF8}{min}理\end{CJK}, human \emph{and} universal.
The Yamaguti category of ``classical'' seems to be concerned with ``truth'' \begin{CJK}{UTF8}{min}真\end{CJK} only: absolute and universal, but detached.
The ``romantic'', in contrast, cultivates ``reason'' \begin{CJK}{UTF8}{min}理\end{CJK} by (self-)reflection towards meaning and understanding of ``reality'' \begin{CJK}{UTF8}{min}実\end{CJK} and ``truth'' \begin{CJK}{UTF8}{min}真\end{CJK}.
Yamaguti's vision aims for a confluence of both these currents of mathematical thought, which  separated so deplorably in the 20th century.

Mathematical truth however, once let loose, develops a dynamics of her own~--
so much so, in fact, that we might suspect her to carry some innate Truth, after all.
The Thomist definition requires an ``ad{\ae}quatio'' between mathematics, as ``intellect'', and some ``things".
But what are those ``things", in mathematics?
Concerning the lofty ``glass bead game'' \cite{Hesse} of mathematical truth~-- detached, universal, and absolute~-- how could it apply to anything outside herself, ever?
And here mathematics only stands for a particularly simple, a particularly limited but comparatively well-defined, mode of thought.
Law, philosophy, language, religion are so much more complex, by comparison.

Heinrich Hertz describes the Thomist adequacy dynamically, as a commutative diagram: 
\begin{quotation}
\noindent
{\emph{``We produce inner images [sketches] or symbols [models] of the external objects, such that the consequences of the images, by abstract thought, always shall be the images of the consequences, in nature, of the objects which the images describe''} \cite{Hertz}.}
\end{quotation}
Over-interpretation may lead to a kind of incantation by mathematical voodoo.
Faster-than-light tachyons $v>c$ are an example: over-stretched  standard relativistic terms like $\sqrt{1-(v/c)^2}$ lead to imaginary momenta or rest masses, and other quantities without physical interpretation.
In a similar vein, the ongoing quest for the universal quantum computer, beyond the old idea of analog computers, is probably the broadest experimental test ever of the universal accessibility -- and relevance -- of the vast complex Schrödinger phase space model of quantum mechanics.

Even if the experimental \emph{``res" }side of the Thomist ad{\ae}quatio seems relatively clear, scientifically speaking, the theoretical \emph{``intellectus"} part requires some actual ``incarnation", some mindful representation, some word, some language, preferably some \begin{otherlanguage}{greek}λόγος\end{otherlanguage}{ }(logos).
As noted by Galileo Galilei, mathematics steps in here, the designated ``universal language of science" \cite{Galilei}.
Indeed, the mathematical formulas of Isaac Newton, James Maxwell, Albert Einstein, Ernst Schrödinger, Werner Heisenberg, and many others, define a language in which those ominous ``Laws of Nature'' appear to be carved in stone.
The lamentable and, to outsiders, inscrutably hermetic symbolism of mathematics reflects that universality.
Logographic scripts, like Chinese or Japanese characters, Egyptian hieroglyphs, or computer icons, are a similarly ingenious and universal semiotic device: content and meaning are symbolized, processed, and conveyed beyond ephemeral phonetic particulars.

In the present paper, specifically,  we have dealt with complex and real time.
``Romantically'' speaking, however, what are these?
Augustine of Hippo develops a concept of time which is very closely related to Dedekind's construction of $\mathbb{R}$ by sections in $\mathbb{Q}$, but sighs:
\begin{quotation}
\noindent
\emph{And I confess to thee, O Lord, that I am still ignorant as to what time is. And again I confess to thee, O Lord, that I know that I am speaking all these things in time, and that I
have already spoken of time a long time, ...
Or, is it possible that I do not know how I can express what I do know? Alas for me! I do not even know the extent of my own ignorance} \cite{Confess}
\end{quotation} 
Such earlier doubts did not discourage him to lucidly interpret the Book Genesis in terms of \emph{co-creation of space-time}, decades later:
\begin{quotation}
\noindent
\emph{... then certainly the world was made, not in time, but simultaneously with time. For that which is made in time is made both after and before some time,—after that which is past, before that which is future. But none could then be past, for there was no creature by whose movements its duration could be measured. However, simultaneously with time the world was made, if in the world’s creation change and motion were created ...
} \cite{Civitate} XI.6.
\end{quotation}
It took another 1500 years before differential geometry, general relativity, and Big Bang cosmology arrived at PDE models. Not to mention the mires of string theory \cite{Penrose}.

From a ``classical'' mathematical viewpoint, of course, we may just talk about dynamics, in the guise of flows $\Phi^t$, as the study of (at least local, and mostly differentiable) group actions $t\in(\mathbb{R},+)$, for continuous time $t$, or of $t=n\in(\mathbb{Z},+)$, in the time-discrete case of ``self-referential'' iterations.
Applications to real time, in the real world, abound.
We identify the flow $\Phi^t$ as the law of nature, or model, which carries us from a Hertzian ``image'' at a starting time $t_0$ to the image of a later or earlier(!) time $t_0+t$.

Extensions to complex time $t\in\mathbb{C}$, as in \eqref{flow}, seem natural, from the axiomatic point of view.
Descriptions by power series, likewise, induce natural ``classical'' extensions to complex time arguments $t$, undeterred by Voodoo and devoid of any associated meaning.

Under the somewhat pretentious title of ``Summa technologiae'' \cite{Lem}, Stanis{\l}aw Lem has described that ``classical'' attitude in chapter V.5 on ``Methodological madness'':
\begin{quotation}\noindent
\emph{Let us imagine a mad tailor who makes all sorts of clothes. 
He does not know anything about people, birds, or plants. 
He is not interested in the world; he does not examine it. ...
Some clothes are spherical, without any openings for heads or feet; others have tubes sewn into them that he calls “sleeves” or “legs.” ...
The tailor is only concerned about one thing: he aims to be consistent. 
He produces symmetrical and asymmetrical dresses, large and small, stretchy and rigid. 
When he starts producing a new item, he makes certain assumptions. 
They are not always the same. 
But once fixed he follows them, and he expects them not to lead to any contradiction. ... 
He deposits the finished clothes at a huge storage. 
If we had access there, we would discover that some of the clothes fit an octopus, others trees, butterflies, or humans. 
We would find clothes for the centaur and the unicorn, and for creatures beyond imagination.
The great majority of his clothes would not find any application. 
Anyone will admit that the tailor’s Sisyphean labor is pure madness.\\
Mathematics works in the same way.}
\end{quotation}

So, what could our ``arxival'' storage be good for?
There are consequences of our results, for meaningful real-time flows.
Blow-up is an example, as well as the relation between linear and quadratic real vector fields via Möbius transformations \eqref{Mobius}.
Complex $w=u+\mi v$, rather than complex time, can be seen as a special form of coupled real systems, separately for $u$ and $v$, with appropriate real meaning attached to $u$ and $v$.
We have already mentioned \cite{Yanagida} for a PDE variant of this interpretation.
Simple examples then may show unexpected relations among the real-time interpretations of different equations.
Consider odd real analytic force laws $f(u)=ug(u^2)$, for example.
Then the superficially different second-order real $v$-pendulum \eqref{pendulumv} in
\begin{align}
\label{pendulumu}
\ddot u\ +\ u\,g(u^2)&=0\,,\\
\label{pendulumv}
\ddot v+v\,g(-v^2)&=0\,,
\end{align}
is closely related to the $u$-pendulum \eqref{pendulumu}, once we replace $u$ there by complex coordinates $w=u+\mi v$.
Simple examples are the widely studied Duffing oscillators, or the transition from the ``soft spring'' mathematical pendulum \eqref{pendulumu} to the ``hard spring'' pendulum \eqref{pendulumv} with a $\sinh v$ force law.
Real solutions $u$ of \eqref{pendulumu} with initial conditions $u(0)=u_0=0$ and $u'(0)=p_0\in\mathbb{R}$ also generate real solutions $v(t):=\mi u(\mi t)$ of \eqref{pendulumv} with nonlinearity $-g$, in their imaginary time direction.
Similarly, initial conditions $u_0\in\mathbb{R},\ p_0=0$ and $v(t):=u(\mi t)$ relate odd pendula $\ddot u\pm f(u)=0$.
Alternatively, we may view the two ``pendula'' as equilibria for reaction-diffusion equations.
Then the relation is induced, mathematically, by a passage from real to imaginary space variable $x$ for Dirichlet and Neumann boundary, respectively -- a quite ominous transition within conventional physics.

Real discretization demands yet another use of complex time.
In \cite{FiedlerClaudia}, we have raised a 1,000\,\euro\ prize question concerning the existence of complex entire homoclinic loops $\Gamma$ for entire vector fields $f$ in several complex dimensions.
See sections 1.8 and 7 there, for further discussion, including some classical results and the relevance of this question for super-exponentially small separatrix splittings under time-periodic forcing or time-discretization.
Poles of real $\Gamma$ in complex time induce real transverse homoclinics and chaotic shift dynamics, under discretization.
See also \cite{FiedlerScheurle} for the related problem of ``invisible chaos''.
Yamaguti, again, was among the pioneers to explore Li-Yorke or Sharkovsky chaos \cite{LiYorke}, as induced by discretization of ODEs; see \cite{YamagutiMatano} and others.

Physics itself is a long-standing source of ``romantic'' complex analysis applications.
For a short summary see section 5.3 of \cite{Marsden}.
In electrostatics, for example, the blue/cyan and orange circle families of figure \ref{fig4} are known as field lines and potential levels of an electrical dipole with charges $\pm1$ located at $w=\pm1$.
Similarly they illustrate flow lines and the potential of a static, incompressible, irrotational planar fluid flow with source and sink at $w=\pm1$.
Section \ref{ResAntipol}, likewise, identifies anti-polynomial vector fields \eqref{ODEF} as divergence-free gradients in \eqref{ODEwbar}.
Real parts $u=u(x,y)$ of holomorphic functions $w=w(x+\mi y)$ are harmonic, i.e.~they satisfy the Laplace equation $\Delta_{x,y} u=0 $.
In particular, this property is therefore invariant under biholomorphic equivalences in $(x,y)$.
The Riemann mapping theorem and its Teichmüller variants then open the door to general geometries.
The resulting force fields $\nabla u$ describe particle dynamics in the limit of high viscosity.
Semilinear elliptic real PDEs like $\Delta u+f(u)=0$ in two space dimensions $(x,y)\in\mathbb{R}^2$ lead to mathematical interpretations with a two-dimensional ``time'' group $(\mathbb{R}^2,+)$ given by spatial shifts in $(x,y)$.
See for example \cite{Vishik} and the references there.
Similarly to complex time $t=r+\mi s\in(\mathbb{C},+)$, spatial shifts in the two real ``time'' directions $x$ and $y$ commute.
Alternatively, real analyticity of solutions defines local extensions to shifts $(x,y)$ in $\mathbb{C}^2$.

The linear Schrödinger equation makes use of complex state spaces for $\psi$.
Motivated by gauge invariance, nonlinear versions involve hybrid terms like $|\psi|^2\psi$.
See for example \cite{Sulem}.
For applied aspects of the nonconservative variant \eqref{PDEpsi}, see \cite{JaquetteHet}.
We also mention the complex Ginzburg-Landau equation as a rich source of physics \cite{CGL}.
The phase space of the 3d Navier-Stokes equation has been complexified, in \cite{LiSinai}, for a purely mathematical exploration of blow-up under complex initial conditions.

We have seen how mathematical modeling is omnipresent in Science.
It lends ``romantic'' meaning to ``classical'' mathematics, e.g., by the Hertz variant of the Thomist ``ad{\ae}quatio'' \cite{Thomas,Hertz}.
Conversely, new ``classical'' mathematics emerges by abstraction, when applications are stripped of ``romantic'' interpretations. 
Very close to such two-way ad{\ae}quatio, however, Albert Einstein wonders:
\begin{quotation}
\noindent
\emph{``How is it possible that mathematics, clearly a product of human thought independent of all experience, fits the objects of reality so excellently? Is it at all conceivable that human reason, without experience, and based on pure thought, fathom the properties of real things through mere thinking?''} \cite{Einstein}
\end{quotation}
\emph{``The unreasonable effectiveness of mathematics in the natural sciences''} \cite{Wigner} weighs in here.
A deep mystery of the \begin{otherlanguage}{greek}λόγος\end{otherlanguage} \, (lógos), the Word of Creation, indeed. See \cite{Truth, Wunder} for further speculation.

In his analysis \cite{YamagutiRomantic}, Yamaguti has envisioned a future confluence and lively interaction among the ``classical'' pure and the ``romantic'' applied currents of mathematical thought.
We have attempted to contribute a small sketch towards his grand vision.



\bigskip

\begin{thebibliography}{9999)999}

{\footnotesize{

\bibitem[APMR22]{MuciEss}
A.~Alvarez-Parrilla and J.~Muciño-Raymundo.
Dynamics of singular complex analytic vector fields with essential singularities. I and II.
\emph{Conform. Geom. Dyn.} \textbf{21} (2017), 126--224, and \emph{J.\ Sing.} \textbf{24} (2022), 1--78.

\bibitem[Arn88]{ArnoldODE}
V.I.~Arnold.
\emph{Geometrical Methods in the Theory of Ordinary Differential Equations.}
Springer-Verlag, Berlin 1988.

\bibitem[Aug396]{Confess}
Augustine of Hippo.
\emph{Confessiones.} (Latin) 396.

\bibitem[Aug426]{Civitate}
Augustine of Hippo.
\emph{De civitate dei.} (Latin) 426.

\bibitem[Ave1180]{Ave} Averroes (Ibn Rushd).
\emph{The Incoherence of the Incoherence.}
S.~van den Bergh (transl.), Cambridge University Press 1954.

\bibitem[BCN02]{critval}
A.F. Beardon, T.K. Carne, T.W. Ng.
The critical values of a polynomial.
\emph{Constr.~Approx.} \textbf{18} (2002), 343--354;
\url{https://doi.org/10.1007/s00365-002-0506-1}

\bibitem[BraDi10]{DiPol}
B.~Branner and K.~Dias.
Classification of complex polynomial vector fields in one complex variable.
\emph{J.\ Difference Eqs.\ Appl.} \textbf{16} (2010), 463--517;
\url{https://doi.org/10.1080/10236190903251746} 

\bibitem[Br90]{Brunovsky}
P.~Brunovsk\'y.
The attractor of the scalar reaction diffusion equation is a smooth graph. 
\emph{J. Dyn. Differential Eqs.} \textbf{2} (1990), 293--323. 

\bibitem[BrFie88]{brfi88}
P.~Brunovsk\'y and B.~Fiedler.
 Connecting orbits in scalar reaction diffusion equations.
 \emph{Dynamics Reported} \textbf{1} (1988), 57--89.

\bibitem[BrFie89]{brfi89}
P.~Brunovsk\'y and B.~Fiedler.
Connecting orbits in scalar reaction diffusion equations {II}: The complete solution.
 \emph{J.~Diff.~Eqs.} \textbf{81} (1989), 106--135.

\bibitem[COS16]{COS}
C.-H. Cho, H. Okamoto, M. Sh\={o}ji. 
A blow-up problem for a nonlinear heat equation in the complex plane of time. 
\emph{Japan J. Ind. Appl. Math.} \textbf{33} (2016), 145--166.

\bibitem[Con78]{Conley}
C.C.~Conley.
\emph{Isolated Invariant Sets and the Morse Index.} 
AMS, Providence R.I. 1978. 

\bibitem[DLYor24]{Yorke-a}
R. De Leo and J.A. Yorke.
Streams and graphs of dynamical systems.
\emph{Qual. Theory Dyn. Syst.} \textbf{24} (2024);
\url{https://doi.org/10.1007/s12346-024-01112-x}

\bibitem[DLYor25]{Yorke-b}
R. De Leo and J.A. Yorke.
Streams, graphs and global attractors of dynamical systems on locally compact spaces. (2025);
\url{https://arxiv.org/abs/2503.02262}

\bibitem[DiGa21]{DiRat}
K.~Dias and A.~Garijo.
On the separatrix graph of a rational vector field on the Riemann sphere.
\emph{J.\ Diff.\ Eqs.} \textbf{282} (2021), 541--565;
\url{https://doi.org/10.1016/j.jde.2021.02.021}
 
\bibitem[DNZ23]{CGL}
G.K.~Duong, N.~Nouaili, H.~Zaag.
Construction of Blowup Solutions for the Complex Ginzburg-Landau Equation with Critical Parameters. 
Mem.~Am.~Math.~Soc. \textbf{1411}, Providence RI 2023. 

\bibitem[Ein21]{Einstein}
A.~Einstein.
Geometrie und Erfahrung. (German)
\emph{Sitzungsber. Preuss. Akad. Wiss.} (1921), 123--130.

\bibitem[FKW24]{Fasondini24}
M. Fasondini, J.R. King, J.A.C. Weideman.
Complex-plane singularity dynamics for blow up in a nonlinear heat equation: analysis and computation. 
\emph{Nonlinearity} \textbf{37} (2024);
\url{https://doi.org/10.1088/1361-6544/ad700b}

\bibitem[Fie21]{Truth}
B.~Fiedler. 
Absolute Truth – a toxic chimera?
In \emph{The Impact of Academic Research
on Character Formation, Ethical Education, and the Communication of Values
in Late Modern Pluralistic Societies.} 
W.~Schweiker, M.~Welker, J.~Witte, S.~Pickard (eds.), 
Evangelische Verlagsanstalt, Leipzig 2021, 65--87.

\bibitem[Fie23]{FiedlerClaudia}
B.~Fiedler.
Real chaos and complex time. (2023); 
\url{https://arxiv.org/abs/2310.08136}

\bibitem[Fie24]{FiedlerShilnikov}
B.~Fiedler.
Scalar polynomial vector fields in real and complex time. 
\emph{Reg. Chaotic Dyn.} \textbf{30} (2025), 188--225.
\url{https://doi.org/10.1134/S1560354725020030}

\bibitem[Fie25]{Wunder} 
B.~Fiedler.
Wunder -- Widersprüche -- Wirklichkeiten.
In \emph{Die unsichtbare Welt. Beiträge zur Weltwahrnehmung in den Wissenschaften.} (German) 
K. Böhmer, D.~Evers, P.-G.~Reinhard (eds.),
Evangelische Verlagsanstalt, Leipzig 2025, 183--213.


\bibitem[FieRo96]{Cascading}
B.~Fiedler and C.~Rocha.
Heteroclinic orbits of semilinear parabolic equations.
\emph{J.~Differ.\ Eqs.} \textbf{125} (1996), 239--281.

\bibitem[FieRo08]{Sturm2}
 B.~Fiedler and C.~Rocha.
 Connectivity and design of planar
global attractors of Sturm type. II: Connection graphs. 
\emph{J. Diff.~Eqs.}~\textbf{244} (2008), 1255--1286.

\bibitem[FieRo09a]{Sturm1}
B.~Fiedler and C.~Rocha.
Connectivity and design of planar
global attractors of Sturm type. I: Bipolar orientations and
Hamiltonian paths.
\emph{Crelle J.~Reine Angew.~Math} \textbf{635} (2009), 71--96.

\bibitem[FieRo09b]{Sturm3}
B.~Fiedler and C.~Rocha.
Connectivity and design of planar global attractors of Sturm type.
III: Small and Platonic examples.
\emph{J.~Dyn.~Diff.~Eqs.} (2009); 
\url{doi.org/10.1007/s10884-009-9149-2}

\bibitem[FieRo20]{ThomSmale} 
B.~Fiedler and C.~Rocha.
Boundary orders and geometry of the signed Thom-Smale complex for Sturm global attractors.
\emph{J.~Dyn.~Diff,~Eqs.} (2020), 32pp.
\url{doi.org/10.1007/s10884-020-09836-5}

\bibitem[FieRo23]{firoSFB}
B.~Fiedler and C.~Rocha.
Design of Sturm global attractors 1: Meanders with three noses, and reversibility.
\emph{Chaos} \textbf{33}, 083127 (2023); \url{https://doi.org/10.1063/5.0147634}

\bibitem[FieRo24]{firoFusco}
B.~Fiedler and C.~Rocha.
Design of Sturm global attractors 2: Time-reversible Chafee-Infante lattices of 3-nose meanders.
\emph{São Paulo J. Math. Sciences.} (2024);
\url{https://doi.org/10.1007/s40863-023-00385-5}

\bibitem[FSV98]{Vishik}
B.~Fiedler, A.~Scheel, M.I.~Vishik. 
Large patterns of elliptic systems in infinite cylinders. 
\emph{J.~Math.\ Pures Appl.} \textbf{77} (1998), 879--907.

\bibitem[FieSch96]{FiedlerScheurle}
B. Fiedler and J. Scheurle.
\emph{Discretization of Homoclinic Orbits, Rapid Forcing and “Invisible” Chaos.}
Mem. Am. Math. Soc. \textbf{570}, Providence RI 1996. 

\bibitem[FieStu25]{FiedlerFila}
B.~Fiedler and H. Stuke.
Real eternal PDE solutions are not complex entire: a quadratic parabolic example.
\emph{J. Ell. Par. Eqs.} (2025), 53pp;
\url{doi.org/10.1007/s41808-024-00309-0}



\bibitem[FMR95]{Rosenstiehl}
H. de Fraysseix, P.O. de Mendez, P. Rosenstiehl.
Bipolar orientations revisited.
\emph{Discrete Appl. Math.} \textbf{56} (1995),
157--179.

\bibitem[For81]{Forster}
O. Forster.
\emph{Lectures on Riemann Surfaces.}
Springer-Verlag, New York 1981.

\bibitem[Gala04]{Galaktionov}
V.A.~Galaktionov.
\emph{Geometric Sturmian Theory of Nonlinear Parabolic Equations and Applications.} 
Chapman\&Hall/CRC, Boca Raton FL 2004. 

\bibitem[Gal1623]{Galilei}
G.~Galilei.
\emph{Il Saggiatore.} (Italian)
Rome 1623.

\bibitem[Gö1931]{Godel}
K.~Gödel.
Über formal unentscheidbare Sätze der Principia Mathematica und verwandter Systeme I. (German)
\emph{Monatsh.~Math.~Phys.} \textbf{38} (1931), 173--198.

\bibitem[GNSY13]{Yanagida}
J.-S. Guo, H. Ninomiya, M. Shimojo, E. Yanagida.
Convergence and blow-up of solutions for a complex-valued heat equation with a quadratic nonlinearity. 
\emph{Trans. Am. Math. Soc.} \textbf{365} (2013), 2447--2467. 

\bibitem[HMO02]{Hale}
J.K.~Hale, L.T.~Magalhães, W.M.~Oliva.
\emph{An Introduction to Infinite Dimensional Dynamical Systems -- Geometric Theory.}
Springer-Verlag, New York 1984. 

\bibitem[Hartm02]{Hartman}
Ph. Hartman.
\emph{Ordinary Differential Equations.}
SIAM, Providence RI 2002.

\bibitem[Her1894]{Hertz}
H.~Hertz.
\emph{Die Principien der Mechanik in neuem Zusammenhange dargestellt.} (German)
Leipzig 1894.

\bibitem[Hes43]{Hesse}
H.~Hesse.
\emph{Das Glasperlenspiel. Versuch einer Lebensbeschreibung des Magister Ludi Josef Knecht samt Knechts hinterlassenen Schriften.} (German)
Fretz \& Wasmuth, Zürich 1943.

\bibitem[IlYa08]{Ilyashenko}
Y. Ilyashenko and S. Yakovenko.
\emph{Lectures on Analytic Differential Equations.}
AMS, Providence RI 2008.

\bibitem[Jaq21]{Jaquetteqp}
J. Jaquette.
Quasiperiodicity and blowup in integrable subsystems of nonconservative nonlinear Schrödinger equations. 
\emph{J. Dyn. Differ. Eqs.} \textbf{36} (2024), 1--25. 
\url{https://doi.org/10.1007/s10884-021-10112-3}

\bibitem[JLT22a]{JaquetteHet}
J.~Jaquette, J.-P.~Lessard, A.~Takayasu.
Global dynamics in nonconservative nonlinear Schrödinger equations. 
\emph{Adv. Math.} \textbf{398} (2022), 108234. 

\bibitem[JLT22b]{JaquetteStuke}
J.~Jaquette, J.-P.~Lessard, A.~Takayasu.
Singularities and heteroclinic connections in complex-valued evolutionary equations with a quadratic nonlinearity. 
\emph{Commun. Nonlinear Sci. Numer. Simul.} \textbf{107} (2022), 106188.

\bibitem[Jol89]{Jolly}
M.S.~Jolly.
Explicit construction of an inertial manifold for a reaction diffusion equation. 
\emph{J. Differential Eqs.} \textbf{78} (1989), 220--261. 

\bibitem[Jost06]{Jost}
J. Jost.
\emph{Compact Riemann Surfaces. An Introduction to Contemporary Mathematics.}
Springer-Verlag, Berlin 2006.

\bibitem[KMM04]{Mischaikow}
T.~Kaczynski, K.~Mischaikow, M.~Mrozek.
\emph{Computational Homology.} 
Springer-Verlag, New York 2004. 

\bibitem[KaiLe25]{LebEss}
N.~Kainz and D.~Lebiedz.
Separatrix configurations in holomorphic flows.
(2025), 20pp.
\url{https://doi.org/10.48550/arXiv.2505.14594}

\bibitem[KliRou21]{Rousseaud=4}
M.~Klimeš and Ch.~Rousseau.
Remarks on rational vector fields on $\mathbb{CP}^1$.
\emph{J.\ Dyn.\ Control Syst.} \textbf{27} (2021), 293--320. 

\bibitem[Lam09]{Lamotke}
K. Lamotke.
\emph{Riemannsche Flächen.} (German)
Springer-Verlag, Heidelberg 2009.

\bibitem[Lap23]{Lappicy}
P. Lappicy and E. Beatriz. 
An energy formula for fully nonlinear degenerate parabolic equations in one spatial dimension.
\emph{Math. Ann.} (2023); \url{https://doi.org/10.1007/s00208-023-02740-5}

\bibitem[Lem64]{Lem}
S.~Lem. 
\emph{Summa technologiae.} (Polish),
after J. Zylinska (transl.), 2013.
Kraków 1964.

\bibitem[LiSin08]{LiSinai}
D. Li and Y.G. Sinai.
Blow ups of complex solutions of the 3D Navier-Stokes system and renormalization group method.
\emph{J. Eur. Math. Soc.} \textbf{10} (2008), 267--313;
\url{https://doi.org/10.4171/JEMS/111}

\bibitem[LiYor75]{LiYorke}
T-Y. Li and J. A. Yorke. 
Period three implies chaos. 
\emph{Amer. Math., Monthly} \textbf{82} (1975), 985--992.
\url{https://doi.org/10.2307/2318254}

\bibitem[Liu25]{LiuR}
R. Liu.
Local well-posedness of the periodic nonlinear Schrödinger equation with a quadratic nonlinearity $\bar{u}^2$ in negative Sobolev spaces.
\emph{J.~Dyn.~Differential Eqs.} \textbf{37} (2025), 509--538;
\url{https://doi.org/10.1007/s10884-023-10295-x}

\bibitem[MH99]{Marsden}
J.E. Marsden and M.J. Hoffman.
\emph{Basic Complex Analysis.}
Freeman, New York 1999.

\bibitem[Mas82]{Masuda1}
K. Masuda.
Blow-up of solutions of some nonlinear diffusion equations.  
\emph{North-Holland Math. Stud.} \textbf{81} (1982), 119--131. 

\bibitem[Mas84]{Masuda2}
K. Masuda.
Analytic solutions of some nonlinear diffusion equations. 
\emph{Math. Z.} \textbf{187} (1984), 61--73. 

\bibitem[Mat82]{MatanoLap}
H.~Matano.
Nonincrease of the lap-number of a solution for a one-dimensional semilinear parabolic equation. 
\emph{J. Fac. Sci. Univ. Tokyo I A} \textbf{29} (1982), 401--441. 

\bibitem[Mat88]{Matano}
H.~Matano.
Asymptotic behavior of solutions of semilinear heat equations on
$(S^1$.
In \emph{Nonlinear Diffusion Equations and their Equilibrium States
  {II}}. W.-M.~Ni, L.A.~Peletier, J.~Serrin (eds.), Springer-Verlag,
New York 1988.

\bibitem[Muc02]{MuciRiem}
J.~Muciño-Raymundo.
Complex structures adapted to smooth vector fields. 
\emph{Math.\ Ann.} \textbf{322} (2002), 229--265. 

\bibitem[New1687]{Newton}
I.~Newton.
\emph{Philosophiae Naturalis Principia Mathematica.} (Latin)
London 1687.

\bibitem[Noy98]{Noy}
M.~Noy.
Enumeration of noncrossing trees on a circle. 
\emph{Discr.~Math.} \textbf{180} (1998), 301--313. 

\bibitem[Pa69]{Palis}
J.~Palis.
 On Morse-Smale dynamical systems.
 \emph{Topology} \textbf{8} (1969), 385--404.

\bibitem[PaSm70]{PalisSmale}
J.~Palis and S.~Smale.
 Structural stability theorems.
 In \emph{Global Analysis},  S. Chern, S. Smale (eds.). Proc. Symp. in
 Pure Math.~vol.~XIV.~AMS, Providence 1970. 

\bibitem[PdM82]{PalisdeMelo}
J.~Palis and W.~de~Melo.
\emph{Geometric Theory of Dynamical Systems. An Introduction.}
Springer-Verlag, New York 1983. 

\bibitem[Pei62]{Peixoto1}
M.M. Peixoto.
Structural stability on two-dimensional manifolds.
\emph{Topology} \textbf{2} (1963), 101--120. 

\bibitem[Pei63]{Peixoto2}
M.M. Peixoto.
Structural stability on two-dimensional manifolds. A further remark. 
\emph{Topology} \textbf{2} (1963), 179--180.

\bibitem[Pen16]{Penrose}
R. Penrose.
\emph{Fashion, Faith, and Fantasy in the New Physics of the Universe.}
Princeton University Press, 2016.

\bibitem[Pil92]{Pilyugin}
S.Yu.~Pilyugin.
\emph{Introduction to Structurally Stable Systems of Differential Equations.} 
Birkhäuser, Basel 1992. 

\bibitem[QS19]{Quittner}
P. Quittner and Ph. Souplet.
\emph{Superlinear Parabolic Problems. Blow-Up, Global Existence and Steady States.} 2nd ed. Birkhäuser, Cham 2019. 

\bibitem[Rel40]{Rellich}
F. Rellich.
Elliptische Funktionen und die ganzen Lösungen von $y''=f(y)$.
\emph{Math. Z.} \textbf{47} (1940), 153--160.
\url{https://doi.org/10.1007/bf01180954} 

\bibitem[RFLN25]{firoS1}
C.~Rocha, B.~Fiedler, A.~López-Nieto:
A classification of global attractors for $\mathbb{S}^1$-equivariant parabolic equations.
(2025), 36pp.
\url{https://doi.org/10.48550/arXiv.2507.10051}

\bibitem[Slo24a]{oeispol}
L.P.A. Sloane (ed.).
The Online encyclopedia of integer sequences.
A002995 (2024). \url{https://oeis.org/A002995}

\bibitem[Slo24b]{oeisanti}
L.P.A. Sloane (ed.).
The Online encyclopedia of integer sequences.
A296533 (2024). \url{https://oeis.org/A296533}

\bibitem[Sot73]{Sotomayor}
J. Sotomayor.
Generic one-parameter families of vector fields on two-dimensional manifolds.
\emph{Publ. Math. I.H.E.S.} \textbf{43} (1973), 5--46.

\bibitem[Sto00]{countpol2}
A. Stoimenow.
On the number of chord diagrams.
\emph{Disc. Math.} \textbf{218} (2000), 209--233.

\bibitem[Stu17]{Stukediss}
H. Stuke.
\emph{Blow-up in Complex Time.}
Dissertation Thesis, Freie Universität Berlin 2017. 
\url{http://dx.doi.org/10.17169/refubium-11743}

\bibitem[Stu18]{Stukearxiv}
H. Stuke.
\emph{Complex time blow-up of the nonlinear heat equation.}
(2018); 
\url{https://arxiv.org/abs/1812.10707} 

\bibitem[SuSu99]{Sulem}
C.~and P.-L.~Sulem.
\emph{The Nonlinear Schrödinger Equation. Self-Focusing and Wave Collapse.}
Springer-Verlag, New York 1999. 

\bibitem[Tho1274]{Thomas}
Thomas Aquinas.
\emph{Summa Theologica.} (Latin)
1265-1274.

\bibitem[TLJO22]{JaquetteMasuda}
A.~Takayasu, J.-P.~Lessard, J.~Jaquette, H.~Okamoto.
Rigorous numerics for nonlinear heat equations in the complex plane of time.
\emph{Numer. Math.} \textbf{151} (2022), 693--752. 

\bibitem[Ush80]{Ushiki1}
S. Ushiki. 
On unstable manifolds of analytic diffeomorphisms of the plane.
\emph{RIMS Kyoto Kokyuroku} \textbf{403} (1980), 1--7.

\bibitem[Ush81]{Ushiki2}
S. Ushiki. 
Unstable manifolds of analytic dynamical systems.
\emph{J. Math. Kyoto Univ.} \textbf{21} (1981), 763--785.

\bibitem[Wal72]{countpol1} 
D.W. Walkup.
The number of plane trees.
\emph{Mathematika} \textbf{19} (1972), 200--204.

\bibitem[Web24]{dynameis}
J.~Weber.
\emph{Dynameis. Bausteine zu einer Geschichte der Virtualität.} (German)
de Gruyter, Berlin 2024;
\url{https://doi.org/10.1515/9783111322520}

\bibitem[Wig60]{Wigner}
E.P.~Wigner.
The unreasonable effectiveness of mathematics in the natural sciences.
\emph{Comm.~Pure Appl.~Math.} \textbf{13} (1960), 1-14.

\bibitem[Wit41]{Wittich1}
H. Wittich.
Ganze Lösungen der Differentialgleichung $w''=f(w)$.
\emph{Math. Z.} \textbf{47} (1941), 422--426.
\url{https://doi.org/10.1007/BF01180973}

\bibitem[Wit50]{Wittich2}
H. Wittich.
Ganze transzendente Lösungen algebraischer Differentialgleichungen.
\emph{Math. Ann.} \textbf{122} (1950), 37--46.
\url{https://doi.org/10.1007/BF01342967}

\bibitem[Yam75]{YamagutiPDE}
M.~Yamaguti.
A certain semilinear system of partial differential equations.
\emph{Dynamical Systems},
 Proc.~Symp.~Univ.~Warwick 1973/74, 
 \emph{Lect. Notes Math.} \textbf{468}, 78--79;
 Springer-Verlag, Heidelberg 1975.

\bibitem[Yam97]{YamagutiRomantic}
M.~Yamaguti.
What is lacking in the philosophy of complexities. (Japanese)
{\emph{J.~Jap.~Soc.}} \emph{Fuzzy Th.~Syst.} \textbf{9} (1997), 603--613. 
\url{https://doi.org/10.3156/jfuzzy.9.5_603}

\bibitem[YamMat79]{YamagutiMatano}
M.~Yamaguti and H.~Matano.
Euler’s finite difference scheme and chaos.
\emph{Proc.~Jap.~Acad., Ser.} \textbf{A 55} (1979), 78--80. 

\bibitem[Zel68]{Zelenyak}
T.I.~Zelenyak.
Stabilization of solutions of boundary value problems for a second order parabolic equation with one space variable.
\emph{Diff.~Eqs.} \textbf{4} (1968), 17--22.
}}

\end{thebibliography}
\end{document}